\pdfoutput=1
\documentclass[english,letterpaper,reqno]{amsart}
\usepackage[T1]{fontenc}
\usepackage[latin9]{inputenc}
\usepackage{setspace}
\usepackage{textcomp}
\usepackage{enumitem}
\usepackage{amsmath, amssymb,bm, cases, mathtools, thmtools}
\usepackage{graphicx}
\usepackage{esint} %
\usepackage{xargs}[2008/03/08]
\usepackage{datetime}

\usepackage[usenames,dvipsnames]{xcolor}

\usepackage{euscript,microtype}

\usepackage[%
    minnames=4,maxnames=99,maxcitenames=4,
    style=alphabetic,
    doi=false,url=false,
    firstinits=true,hyperref,natbib,backend=bibtex]{biblatex}
\renewbibmacro{in:}{%
  \ifentrytype{article}{}{\printtext{\bibstring{in}\intitlepunct}}}
\AtEveryBibitem{\clearfield{issn}}
\AtEveryBibitem{\clearfield{isbn}}
\AtEveryCitekey{\clearfield{issn}}
\AtEveryCitekey{\clearfield{isbn}}
\AtEveryBibitem{\ifentrytype{book}{\clearfield{pages}}{}}
\DeclareFieldFormat[book,article,unpublished,misc]{title}{\textit{#1}}
\DeclareFieldFormat[article]{journaltitle}{{#1}}

\bibliography{sample_dist_paper,dans_references}

\usepackage[colorlinks,citecolor=blue,urlcolor=blue,linkcolor=RawSienna,hypertexnames=true]{hyperref}
\usepackage{hypernat}

\usepackage[capitalize]{cleveref}

\crefname{lem}{Lemma}{Lemmas}
\crefname{cor}{Corollary}{Corollaries}
\crefname{thm}{Theorem}{Theorems}
\crefname{assumption}{Assumption}{Assumptions}

\numberwithin{equation}{section}

\declaretheorem[style=plain,numberwithin=section,name=Theorem]{thm}
\declaretheorem[style=plain,sibling=thm,name=Lemma]{lem}

\declaretheorem[style=definition,sibling=thm,name=Definition]{defn}
\declaretheorem[style=definition,name=Assumption]{assumption}

\declaretheorem[style=remark,qed=$\triangleleft$,sibling=thm,name=Remark]{rem}

\numberwithin{thm}{section}

\newlist{casenv}{enumerate}{4}
 \setlist[casenv]{leftmargin=*,align=left,widest={iiii}}
 \setlist[casenv,1]{label={{\itshape\ \casename} \arabic*.},ref=\arabic*}
 \setlist[casenv,2]{label={{\itshape\ \casename} \roman*.},ref=\roman*}
 \setlist[casenv,3]{label={{\itshape\ \casename\ \alph*.}},ref=\alph*}
 \setlist[casenv,4]{label={{\itshape\ \casename} \arabic*.},ref=\arabic*}

\def\[#1\]{\begin{align}#1\end{align}}
\def\*[#1\]{\begin{align*}#1\end{align*}}

\newcommand{\LATER}[1]{\error}
\newcommand{\fLATER}[1]{\error}
\newcommand{\TBD}[1]{\error}
\newcommand{\fTBD}[1]{}
\newcommand{\PROBLEM}[1]{\error}
\newcommand{\fPROBLEM}[1]{\error}
\newcommand{\NA}[1]{#1}

\usepackage{babel}

 \providecommand{\casename}{Case}

\begin{document}

\newcommand{\ErdosRenyi}{Erd\H{o}s--R\'enyi--Gilbert}

\newcommand{\defnphrase}[1]{\emph{#1}}

\global\long\def\defas{\vcentcolon=}

\global\long\def\st{\,:\,}

\global\long\def\dist{\ \sim\ }

\global\long\def\given{\mid}

\global\long\def\distiid{\overset{iid}{\dist}}

\global\long\def\distind{\overset{ind}{\dist}}

\global\long\def\Naturals{\mathbb{N}}

\global\long\def\Rationals{\mathbb{Q}}

\global\long\def\Reals{\mathbb{R}}

\global\long\def\BorelSets{\mathcal{B}}

\global\long\def\Nats{\mathbb{N}}

\global\long\def\Ints{\mathbb{Z}}

\global\long\def\NNInts{\Ints_{+}}

\global\long\def\NNExtInts{\overline{\Ints}_{+}}

\global\long\def\Cantor{\2^{\mathbb{N}}}

\global\long\def\NNReals{\Reals_{+}}

\global\long\def\as{\textrm{ a.s.}}

\global\long\def\epi{\textrm{epi}}

\global\long\def\intr{\textrm{int}}

\global\long\def\conv{\textrm{conv}}

\global\long\def\cone{\textrm{cone}}

\global\long\def\aff{\textrm{aff}}

\global\long\def\cone{\textrm{cone}}

\global\long\def\dom{\textrm{dom}}

\global\long\def\cl{\textrm{cl}}

\global\long\def\ri{\textrm{ri}}

\global\long\def\grad{\nabla}

\global\long\def\imp{\Rightarrow}

\global\long\def\downto{\!\downarrow\!}

\global\long\def\upto{\!\uparrow\!}  \global\long\def\AND{\wedge}

\global\long\def\OR{\vee}

\global\long\def\NOT{\neg}

\global\long\def\PowerSet{\mathcal{P}}

\global\long\def\Measures{\mathcal{M}}

\global\long\def\ProbMeasures{\mathcal{M}_{1}}

\global\long\def\equaldist{\overset{d}{=}}

\global\long\def\equalprob{\overset{p}{=}}

\global\long\def\inv{^{-1}}

\global\long\def\norm#1{\lVert#1 \rVert}

\global\long\def\event#1{\left\lbrace #1 \right\rbrace }

\global\long\def\tuple#1{\langle#1 \rangle}

\global\long\def\bspace{\Omega}

\global\long\def\bsa{\mathcal{A}}

\global\long\def\borelspace{(\bspace,\bsa)}

\global\long\def\card#1{\##1}

\global\long\def\iid{i.i.d.\ }

\global\long\def\iff{iff\ }

\global\long\def\gprocess#1#2{(#1)_{#2}}

\global\long\def\nprocess#1#2#3{\gprocess{#1_{#3}}{#3 \in#2}}

\global\long\def\process#1#2{\nprocess{#1}{#2}n}

\newcommand{\EE}{\mathbb E}

\global\long\def\expect#1{\EE [#1]}

\global\long\def\var#1{\mbox{var}\left[#1\right]}

\global\long\def\equalas{\overset{\mathrm{a.s.}}{=}}

\global\long\def\abs#1{\lvert#1 \rvert}

\global\long\def\norm#1{\lVert#1\rVert}

\global\long\def\inedge#1{e_{#1}^{\mbox{in}}}

\global\long\def\outedge#1{e_{#1}^{\mbox{out}}}

\global\long\def\intd{\mathrm{d}}

\global\long\def\suchthat{\mid}

\global\long\def\Pr{\mbox{P}}

\global\long\def\convPr{\xrightarrow{\,p\,}}

\global\long\def\floor#1{\lfloor#1\rfloor}

\global\long\def\asympLim#1{,\ #1\rightarrow\infty}

\global\long\def\diri{\mathrm{Diri}}

\global\long\def\categ{\mathrm{Cat}}

\global\long\def\betaDist{\mathrm{Beta}}

\global\long\def\bern{\mathrm{Bernoulli}}

\global\long\def\bernoulli#1{\mbox{Bernoulli}(#1)}

\global\long\def\binDist{\mathrm{Bin}}

\global\long\def\uniDist{\mathrm{Uni}}

\global\long\def\poiDist{\mathrm{Poi}}

\global\long\def\gammaDist{\mathrm{Gamma}}

\global\long\def\PP{\Pi}

\global\long\def\PPnu{\Pi_{\nu}}

\global\long\def\vertexset#1{v\left(#1\right)}

\global\long\def\edgeset#1{e\left(#1\right)}

\global\long\def\PPbelowth#1{\Pi_{\nu,\le#1}}

\global\long\def\PPaboveth#1{\Pi_{\nu,>#1}}

\global\long\def\linverse{^{-1}}

\global\long\def\threshold{T_{\nu}}

\global\long\def\popthreshold{T_{\nu,\mbox{pop}}}

\global\long\def\popgraph{P_{\nu}}

\global\long\def\upperthreshold{T_{\nu,u}}

\global\long\def\numDegNu#1{N_{\nu,#1}}

\global\long\def\PP{\Pi}

\global\long\def\PPnu{\Pi_{\nu}}

\global\long\def\vertexset#1{v\left(#1\right)}

\global\long\def\edgeset#1{e\left(#1\right)}

\global\long\def\PPbelowth#1{\Pi_{\nu,\le#1}}

\global\long\def\PPaboveth#1{\Pi_{\nu,>#1}}

\global\long\def\linverse{^{-1}}

\global\long\def\threshold{T_{\nu}}

\global\long\def\popthreshold{T_{\nu,\mbox{pop}}}

\global\long\def\popgraph{P_{\nu}}

\global\long\def\upperthreshold{T_{\nu,u}}

\global\long\def\numDegNu#1{N_{\nu,#1}}

\global\long\def\pointspace{\mathbb{M}}

\global\long\def\degNuFn{D_{\nu}}

\newcommand{\Lebesgue}{\Lambda}

\newcommand{\NatSubs}[1]{\tilde \Nats_{#1}}

\newcommand{\StarF}{S}
\newcommand{\IsoF}{I}

\title[Kallenberg Exchangeable Graphs] 
{The class of random graphs arising from exchangeable random measures}

\author[V.~Veitch]{Victor Veitch}
\address{University of Toronto\\Department of Statistical Sciences\\Sidney Smith Hall\\100 St George Street\\Toronto, Ontario\\M5S 3G3\\Canada}

\author[D.~M.~Roy]{Daniel M.~Roy}

\begin{abstract}
We introduce a class of random graphs that we argue meets many 
of the desiderata one would demand of a model to serve as the foundation 
for a statistical analysis of real-world networks. The class of random graphs 
is defined by a probabilistic symmetry:  invariance of the distribution of each graph to an arbitrary relabelings of its vertices. 
In particular, following Caron and Fox, we interpret a symmetric simple point process on $\mathbb{R}_+^2$ as 
the edge set of a random graph, and formalize the probabilistic symmetry as joint exchangeability of the point process.
We give a representation theorem for the class of random graphs satisfying this symmetry 
via a straightforward specialization of Kallenberg's representation theorem for jointly exchangeable random measures on $\mathbb{R}_+^2$. 
The distribution of every such random graph is characterized by three (potentially random) components: 
a nonnegative real $I \in \mathbb{R}_+$, 
an integrable function $S: \mathbb{R}_+ \to \mathbb{R}_+$,
and a symmetric measurable function $W: \mathbb{R}_+^2 \to [0,1]$ that satisfies several weak integrability conditions. 
We call the triple $(I,S,W)$ a graphex, in analogy to graphons, 
which characterize the (dense) exchangeable graphs on $\Nats$. 
Indeed, the model we introduce here contains the exchangeable graphs as a special case, as well as the "sparse exchangeable" model of Caron and Fox.
We study the structure of these random graphs,
and show that they can give rise to interesting structure, including sparse graph sequences. 
We give explicit equations for expectations of certain graph
statistics, as well as the limiting degree distribution. We also show that
certain families of graphexes give rise to random graphs that, asymptotically,
contain an arbitrarily large fraction of the vertices in a single connected
component.
\end{abstract}

\maketitle

\begin{center}
  \begin{minipage}{.80\linewidth}
    \setcounter{tocdepth}{1}
    \tableofcontents
  \end{minipage}
\end{center}

\section{Introduction}

Random graph models are a key tool for understanding the structure of real-world networks,
especially through data.
In particular, 
a random graph model can serve as the foundation for a statistical analysis:
observed link structure is modeled as a realization from the random graph model, whose parameters are in some unknown configuration. The goal is to then infer the configuration of the parameters, and in doing so,  understand properties of the network that gave rise to the observed link structure.

The quality of the inferences we can make 
depends in part on the fidelity of the model,
but building realistic models of networks is challenging:
the models must be simple enough to be tractable,
yet flexible enough to accurately represent a wide range of phenomena.
In the setting of densely connected networks,
the well-known exchangeable graph model provides a
tractable yet general framework.
However, the vast majority of real-world networks are sparsely 
connected---two nodes chosen at
random are very unlikely to be directly connected by a link.
Accordingly, for some configuration of their parameters, 
realistic random graph models for networks must be
sparse, exhibiting only a vanishing fraction of all possible edges as
they become large. 
At the same time, the link structure of real-world networks is rich:
e.g., in social networks, 
phenomena such as
homophily (informally, friends of friends are more likely to be friends),
``small-world'' connectivity (two randomly chosen individuals are likely to be connected by a short path of friendship), 
and 
power law degree distributions (the number of friends an individual may have varies across many orders of magnitude)
are common~\citep{Newman:2009,Durrett:2006}.
It is a remarkable gap in modern statistical practice that there is no 
general framework for the statistical analysis of real-world networks. 
 There is no shortage of proposals for random graph models of real-world networks;
 however, these models tend to be ad hoc,
 exhibiting certain properties of real-world networks by design, 
 but behaving pathologically in other aspects.
 It is difficult to assess the statistical
 applicability of such models. 

One approach to identifying large but tractable families of random graphs
is to consider the family of all random graphs satisfying 
a small number of natural assumptions. In this paper, we define
a class of random graph models in terms of a single invariance principle: that the distribution of a graph should be invariant to an arbitrary relabeling of its vertices.
 From this assumption, we derive and study a general class of
random graphs suitable for modeling network structures. We show that these random graphs admit a simple, 
tractable specification and give rise to complex structures of the kinds observed in real world networks.
Moreover, our derivation is closely analogous to an approach that has been used to define broadly useful
statistical models in other settings. For instance, the classical \iid setting 
and 
the graphon setting for densely connected networks are 
both derived from analogous invariance assumptions~\cite{Orbanz_Roy_Exchangeable_Struct}. 
Indeed, we show that the exchangeable graph models are a special case of the models we derive here. 
These observations suggest that the models we identify in this paper may be broadly useful for the statistical
analysis of real-world networks.

To explain our approach we begin by reviewing a closely related approach
used to define models for
the statistical analysis of densely connected networks. In this setting,
networks are modeled as random graphs represented by their adjacency
matrices; an observed $n\times n$ adjacency matrix is modeled as the
leading size-$n$ principal submatrix of some infinite array of random
variables. The infinite structure automatically provides consistent
models for datasets of different size.  
The foundational structural
assumption by which the dense graph framework is defined is a probabilistic
symmetry: \emph{joint exchangeability of the infinite array}. This is the requirement that
the distribution of the infinite array is invariant under joint
permutations of the indices of the array; intuitively, this means that
the labeling of the vertices of a graph does not carry information
about its structure.

The statistical framework can be 
derived using the Aldous--Hoover
representation theorem for jointly exchangeable arrays. Specialized to
the case of infinite adjacency matrices, this theorem asserts that the
adjacency matrix of a random graph on $\Nats$ is jointly exchangeable \iff its
distribution can be written as a mixture over a certain privileged
family of distributions (namely, the ergodic measures).  Each member of this
family is specified in terms of a symmetric, measurable function $W :
[0,1]^2 \to [0,1]$, now known as a \emph{graphon}.  It follows that
the space of probability distributions on $n \times n$ observations of
a densely connected networks can be parameterized by the space of
graphons. A particular consequence of the theorem is that the expected
number of links among every $n$ individuals is ${n \choose 2} \| W \|_1 $;
i.e., the graph is either empty or dense.  As stated plainly in
\cite{Orbanz_Roy_Exchangeable_Struct}, these models are thus
misspecified as statistical models for real-world networks.

The derivation of the dense graph framework is a particular instance
of a general recipe for constructing statistical models: a
probabilistic symmetry is assumed on some infinite random structure
and an associated representation theorem characterizes the 
ergodic measures,
 forming the foundation of a
framework for statistical analysis.  The first main contribution of
the present paper is the analogous representation theorem for the
sparse (and dense) graph setting, which we arrive at by a straightforward 
adaptation of a result of 
Kallenberg~\cite{Kallenberg_Random_Meas_Plane,Kallenberg:2005}.  
Our inspiration
comes from recent paper of Caron and Fox~\cite{Caron_Fox_CRM_Graphs}
that exploits a connection between random measures and random graphs
to exhibit a class of sparse random graphs. In their paper, they observe
that their random graphs satisfy a natural analogue of joint
exchangeability when considered as a point process and make use of an
associated representation theorem to study the model.  The present
paper reverses this chain of reasoning, beginning with the symmetry on
point processes and elucidating the full family of random graphs that
arise from the associated representation theorem. In the graph
context, joint exchangeability of point processes retains the
interpretation that the labels of vertices carry no information about the
structure of the graph.

Following Caron and Fox, we represent random graphs as an infinite
simple point processes on $\NNReals^2$ with finite random graphs given
by truncating the support of the point process to a finite set (see
\cref{fig:random_graph_measure_corresp}). The representation theorem
associated to joint exchangeability of point processes is known by the
work of Kallenberg
\cite{Kallenberg_Random_Meas_Plane,Kallenberg:2005}.
We arrive at our
representation theorem by a straightforward translation of this result
into the random graph setting.  The random graphs picked out by our
representation theorem have three possible components: isolated edges,
infinite stars, and a final piece that provides the interesting
graph structure. The basic object for the distributions of these
random graphs is a triple $(\IsoF,\StarF,W)$ where $\IsoF\in\NNReals$, $S:\NNReals\to\NNReals$ is integrable, 
and $W:\NNReals^{2}\to[0,1]$ is a symmetric measurable function satisfying certain weak integrability conditions. 
(See \cref{lfgraphex}; $W$ integrable is sufficient but not necessary.)
We call the triple a \emph{graphex}. In this paper we focus on random graphs 
without isolated edges or infinite stars, and so we take $\IsoF = \StarF=0$;
when there is no risk of confusion, we will use the term \emph{graphex} to refer
to the function $W$ alone with the understanding that the triple is then of the form $(0,0,W)$.
The distribution of every such random graph, which we call a Kallenberg exchangeable graph,
is characterized by some (possibly random) graphex.
Graphexes are the analogues of graphons and the space of
distributions on (sparse) graphs can be parameterized by the space of
graphexes.

\begin{figure}
   \protect \vspace*{2em}
  \centering{}\includegraphics[width=0.5\linewidth]{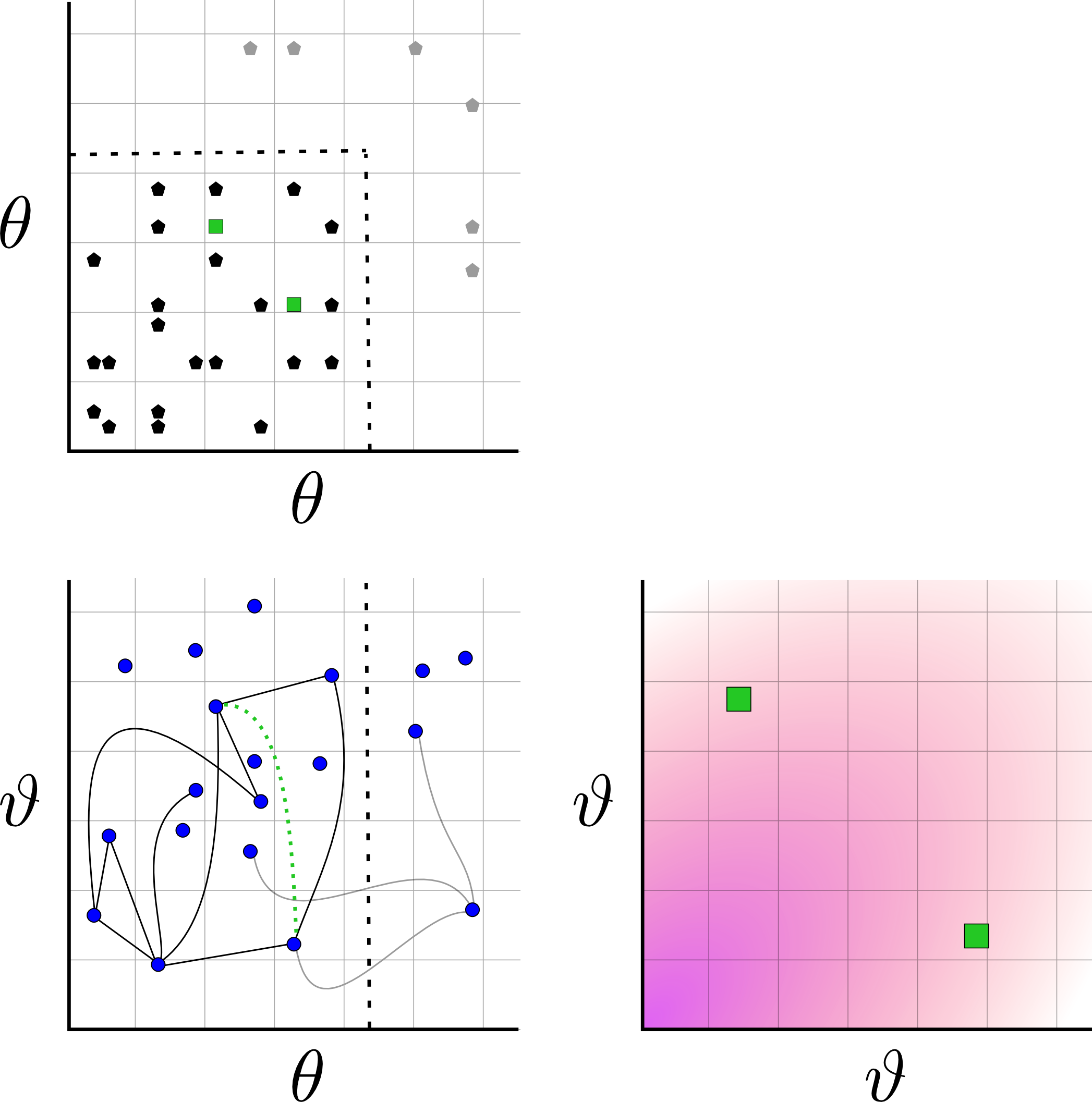}
\caption{\label{fig:graphex_model}
    (Kallenberg exchangeable graph) Random graphs arising from
    exchangeable random measures are characterized by three (potentially random) components: 
    a non-negative real $\IsoF \in \NNReals$, an integrable function $\StarF: \NNReals \to \NNReals$, 
    and a symmetric measurable function $W: \NNReals^2 \to [0,1]$ satisfying some weak integrability conditions.  We call the triple $(\IsoF,\StarF,W)$ a graphex.
    The most interesting structure arises from $W$.  
    A particular $W$ is illustrated by the magenta heatmap (lower right).
    Given $W$, an infinite random
    graph with a vertex set in $\boldsymbol{\theta}$ is generated in
    this model
    according to:\protect \\
    \hspace*{2em}1. Sample a (latent) unit rate Poisson process $\PP$ on $\boldsymbol{\theta}\times\boldsymbol{\vartheta}$.\protect \\
    \hspace*{2em}2. For each pair of points
    $(\theta_{i},\vartheta_{i}),(\theta_{j},\vartheta_{j})\in\PP$
    include\protect \\
    \hspace*{3.2em}edge $(\theta_{i},\theta_{j})$ with probability $W(\vartheta_{i},\vartheta_{j})$.\protect \\
    \hspace*{2em}3. Include $\theta_{i}$ as a vertex whenever
    $\theta_{i}$
    participates\protect \\
    \hspace*{3.2em}in at least one edge.\protect \\
    Finite subgraphs are given by restricting the space
    $\boldsymbol{\theta}$ to be less than some finite value. The lower
    left panel of the figure shows a realization of a latent Poisson
    process with a realization of the edge structure superimposed. A
    finite subgraph (black edges) is given by taking only points with
    $\theta<4.2$. The edge $(3.2,2.1)$ (green, dotted squares) is
    included with probability $W(1.1,4.7)=W(4.7,1.1)$; this is shown
    in the middle panel. Edges that include a point of $\PP$ with
    $\theta>4.2$ (grey, transparent) are not included in the subgraph.
    Vertices, such as $2.7$, that participate only in edges with a
    terminus that has $\theta>4.2$ are not included in the
    subgraph. The upper left panel shows the pictured graph as a
    realization of a random measure on
    $\boldsymbol{\theta}\times\boldsymbol{\theta}$ space.
    \protect \\
}
\end{figure}

It remains to explain the construction of the random graph associated with a graphex.
Let
$\boldsymbol{\theta}=\NNReals$ be the space of labels of the graph,
$\boldsymbol{\vartheta}=\NNReals$ be the space of latent parameters, and
$\PP$ be a unit rate Poisson process on
$\boldsymbol{\theta}\times\boldsymbol{\vartheta}$.  Intuitively, the
random graph is given by independently randomly including each pair of
points in $\PP$ as an edge of graph with a probability determined by
the graphex $W$.  A point of the Poisson process is included as a vertex
of the graph if and only if it participates in at least one
edge. The construction of the random graph is
explained in \cref{fig:graphex_model}.
Formally, treating the collection of edges
$\{(\theta_{i},\theta_{j})\}$ as the basic random object of interest
the generative model given $W$ and $\PP$ is:
\begin{align}
  (\theta_{i},\theta_{j})\given W, \PP & \distind
  \bern(W(\vartheta_{i},\vartheta_{j})).
\end{align}
Finite size graphs are given by restricting to only edges
$(\theta_i,\theta_j)$ such that $\theta_i,\theta_j<\nu$ and including
vertices only if they participate in at least one such edge.  These
distributions are consistent for datasets of different sizes and admit
sparse graphs, 
allowing for the realistic modeling of
physical networks. Moreover, in a sense we make precise in \cref{sec:graphon_models},
the exchangeable graphs derived from the Aldous--Hoover theory are contained
as a subfamily of the Kallenberg exchangeable graphs, and correspond those graphs generated by graphexes of the form $(0,0,W)$ where $W$ is compactly supported, and therefore equal to the dilation of some graphon. Thus the KEG
framework is a generalization of the exchangeable graph framework to the
sparse graph regime.

Let $G_{\nu}$ be the random graph given by truncating the label space
$\boldsymbol{\theta}$ to $[0,\nu]$ (see \cref{fig:graphex_model}); we
call the random graph model $(G_\nu)_{\nu\in\NNReals}$ the
\emph{Kallenberg exchangeable graph} (KEG) associated with $W$.  The
bulk of the present paper is devoted to deriving properties of these
graphs in terms of the graphex $W$.  For simplicity of presentation we
ignore self edges here, giving full statements in the body of the
paper. Let $\mu_{W}(x)=\int_{\NNReals}W(x,y)\intd x$.
\begin{enumerate}

\item Given a point $(\theta,\vartheta)$ in the latent Poisson
  process, the degree of the vertex labeled $\theta$ is Poisson
  distributed with mean $\nu \mu_{W}(\vartheta)$.

\item The expected number of edges $e_{\nu}=\abs{\edgeset{G_{\nu}}}$
  is
  \[
  \expect{e_{\nu}}=\frac{1}{2}\nu^{2}\iint_{\NNReals^{2}}W(x,y)\intd
  x\intd y.
  \]

\item The expected number of vertices
  $v_{\nu}=\abs{\vertexset{G_{\nu}}}$ is
  \[
  \expect{v_{\nu}}=\nu\int_{\NNReals} (1-e^{-\nu\mu_{W}(x)} )\intd x.
  \]

\item Subject to some technical constraints, the scaling limit of the
  asymptotic degree distribution has an explicit expression in terms
  of $W$. Let $k_{\nu}$ be some non-decreasing function of $\nu$ and
  let $D_{\nu}$ be the degree of a randomly selected vertex of
  $G_{\nu}$, then
  \[
  P(D_{\nu}\ge k_{\nu}\given
  G_{\nu})\convPr\lim_{\nu\to\infty}\frac{\sum_{k=k_{\nu}}^{\infty}\frac{\nu^{k}}{k!}\int\mu_{W}(x)^{k}e^{-\nu\mu_{W}(x)}\intd
    x}{\int_{\NNReals}(1-e^{-\nu\mu_{W}(x)})\intd x}.
  \]
  This result establishes that the random graph construction in this
  paper can give rise to sparse graphs.

\item Certain choices of $W$ admit highly connected graphs. Suppose
  $W(x,y)=f(x)f(y)$, let $C_{1}(G_{\nu})$ be the largest connected
  component of $G_{\nu}$, and let $\epsilon>0$, then
  \[
  \lim_{\nu\to\infty}P(\abs{C_{1}(G_{\nu})}>(1-\epsilon)\abs{\vertexset{G_{\nu}}})=1.
  \]
  This means that the sparse structure can arise in an interesting
  way: it is not simply a consequence of having a collection of
  disjoint dense graphs.
\end{enumerate}

We begin by giving background on random graph modeling and the use of
probabilistic symmetry in \cref{sec:Background}.  In
\cref{sec:Examples}, we give a number of illustrative examples of
Kallenberg exchangeable graphs to make the construction concrete. In
\cref{sec:Representation-Theorem}, we establish the representation
theorem and give a formal characterization of the models we
derive. In \cref{sec:Expectations}, we derive the first moments of
several graph statistics of $G_{\nu}$ using point process techniques,
allowing self edges. An expression for asymptotic degree distribution
of these graphs in terms of the graphex is derived in
\cref{sec:Degree-distribution}.  Finally, in \cref{sec:Connectivity},
we study the structure of the Kallenberg exchangeable graphs generated by
graphexes of the form $W(x,y)=f(x)f(y)1[x\neq y]$ with the goal of establishing the
asymptotic connectivity structure. Several other interesting features
of these random graphs are uncovered in the course of establishing
this result. In particular, we show that degree power law
distributions and ``small-world'' phenomena arise naturally in this
framework.

\section{Background\label{sec:Background}}
\newcommand{\FGraphs}{\mathcal G}
\newcommand{\id}{\mathrm {id}}
\newcommand{\SSM}{\mathcal N_1}

In order to relate the Kallenberg exchangeable graph model to a diverse range of existing random graph models,
it will be useful to have a general definition for the term `random graph model'.
In this paper, a random graph model is an indexed family of 
graph-valued random variables $G_{s,\phi}$,
where 
$s$ specifies the ``size'' of the graph and takes values in a totally ordered set $S$,
and where $\phi \in \Phi$ determines some distributional properties (and so could play the role of a parameter in a statistical model).  
We will write $\mu_{s,\phi}$ for the distribution of $G_{s,\phi}$.\footnote{In a statistical setting, the family of distributions $\mu_{s,\phi}$ would be the natural structure to call a model.  Here we adopt the language of graph theorists.}
Our definition is deliberately vague about the meaning of `graph-valued' as different models will naturally be described in terms of different concrete spaces.

For example, the well-known \ErdosRenyi\ model is the family of simple random graphs $G_{n,p}$ on $n \in \Nats$ vertices, where each edge appears independently with probability $p \in [0,1]$.  Concretely, we can think of $G_{n,p}$ as a random $n \times n$ adjacency matrix, or equivalently, as a symmetric $n \times n$ 
array of $0/1$-valued (i.e., binary) random variables whose diagonal is zero.
In a statistical setting, we might model the network of friendships among $n$ individuals as a realization of $G_{n,p}$ for some unknown $p$. 
In this case, the goal of statistical analysis would be to make inferences about the parameter $p$ given some particular observed dataset in the form of an adjacency matrix.

The \ErdosRenyi\ model can be seen as special case of the more general random graph model that arises from the graphon theory or from the Aldous--Hoover representation theorem.  In this case, the size again determines the number of vertices, but the parameter is a graphon, i.e., a symmetric, measurable function $W : [0,1]^2 \to [0,1]$.  
(The \ErdosRenyi\ model corresponds with constant graphons $W(x,y) = p$ for some $p \in [0,1]$.)
This class of random graphs are known as the \defnphrase{exchangeable graphs},
although we will sometimes refer to them as the (dense) exchangeable graphs to distinguish them from the Kallenberg exchangeable graphs.

In the exchangeable graph model, the size parameter is the number of vertices. This is the typical
approach to indexing random graph models. 
In contrast, the size parameter of a Kallenberg exchangeable graph model is a non-negative real $\nu$
that is proportional to the square root of the \emph{expected} number of \emph{edges}.

\subsection{Desiderata for random graph models}

For the purpose of modeling real-world networks, 
one of the key properties of a random graph model is the relationship between the number of edges and vertices.
Consider a random graph model $G_{s,\phi}$, fix a parameter $\phi$, and let $s_n \upto \infty$ be some diverging sequence of sizes.  
For a graph $G$, let $\abs{\edgeset{G}}$ and $\abs{\vertexset{G}}$ denote the number of edges and vertices, respectively.
To avoid pathologies, we will assume that $\abs{\vertexset{G}} \to \infty$ as $n \to \infty$.
Then the sequence $(G_{s_n,\phi})$ is \defnphrase{sparse} 
or \defnphrase{not dense}
if, with probability one,
\[
\frac 
{ \sqrt { \abs{\edgeset{G}} } }
{ \abs{\vertexset{G}} } 
  \to 0 \qquad \text{ as $n \to \infty$. }
\]
This condition states that, asymptotically, graphs with $v$ vertices have $o(v^2)$ edges.
More generally, it is interesting to identify whether there is a (potentially random) exponent $k$ such that, asymptotically, there are $\Theta(v^k)$ edges.

For statistical applications,
it is desirable to
impose a desideratum in addition to sparsity.
The prototypical statistical network analysis has the following structure: 
an observed network $g_s$ is modeled as a realization of a random graph $G_{s,\phi}$ 
for some size $s$ and
for some unknown parameter $\phi$;
the goal is to infer the parameter $\phi$.  
In some random graph models, the sequence $G_{s_1,\phi},G_{s_2,\phi},\dotsc$ of graphs 
is a model of the dynamics by which a network grows and evolves.
In the statistical problems motivating this paper, however,
the size parameter $s$ is akin to sample size in the sense that collecting more data
corresponds to choosing larger values of $s$.  
It is therefore natural to demand that the distributions associated with different sizes are ``consistent'' with one another
in the sense that moving from $G_{s,\phi}$ to $G_{t,\phi}$, for $t > s$, can be understood as collecting additional data.

One way to formalize this notion of consistency is to demand that the distributions of the random graphs $G_{s,\phi}$ be \emph{projective}.
Projectivity is defined in terms of a projective system, i.e., 
a family of measurable maps $ ( f_{s,t};\, s \le t \in S)$ where $f_{s,t}$ maps graphs of size $t$ to graphs of size $s \le t$, $f_{s,s}$ is the identity, and $f_{r,t} = f_{r,s} \circ f_{s,t}$ for all $r \le s \le t$.
A random graph model is \defnphrase{projective} if, for some projective system $( f_{s,t};\, s \le t \in S)$, 
it holds that $G_{s,\phi} \equaldist f_{s,t}(G_{t,\phi})$ for every $s < t \in S$ and parameter $\phi$.

Intuitively, this is simply the requirement that a data set
of size $t$ can be understood as a data set of size $s < t$ augmented with
some additional observations.
Indeed, if a random graph model $(G_{s,\phi})$ is projective with respect to a 
projective system $( f_{s,t};\, s \le t \in S)$,
then it is possible to construct the random variables $G_{s,\phi}$ in such a way that
the identity $G_{s,\phi} = f_{s,t}(G_{t,\phi})$ holds almost surely, and not only in distribution.  In view of this, the connection with the idea of $s$ as sample size is clear.  The graphs $G_{s_j,\phi}$ for an increasing sequence $s_1,s_2,\dotsc$ of sizes are nested.

Both the (dense) exchangeable graph model 
and the Kallenberg exchangeable graph model are projective.  
(See \cref{fig:Graphon_model,fig:graphex_model} for illustrations).
The (dens) exchangeable graph model is projective with respect to the maps $f_{m,n}$ that take an $n \times n$ adjacency matrix to its principal leading $m \times m$ submatrix.  In other words, dropping the last $n-m$ rows and columns from $G_{n,W}$ produces an array with the same distribution as $G_{m,W}$.
The Kallenberg exchangeable graph model is projective with respect to the maps $f_{s,t}$ that take a measure on $[0,t]^2$ to its restriction on $[0,s]^2$.  In other words,
$G_{s,W} \equaldist G_{t,W}(\cdot \cap [0,s]^2)$ for all $s,t \in \NNReals$.

The projectivity of the KEG
model sets it apart from random graph models 
that achieve sparsity by percolating dense random graph models such as the exchangeable graph model,
i.e., a sparse graph model is produced by randomly deleting each edge in a dense graph model independently with a probability that grows with the number of vertices.
Examples of such models abound \citep{Bollobas:Janson:Riordan:2007,Bollobas:Riordan:2007,Borgs:Chayes:Cohn:Zhao:2014:sgc1,Borgs:Chayes:Cohn:Zhao:2014:sgc2},
and in some cases consistent estimators have been developed~\citep{Wolfe:Olhede:2013:a,Borgs:Chayes:Cohn:Ganguly:2015,Borgs:Chayes:Smith:2015}.
Each of these random graph models is parametrized by a size $n$ 
that determines the number of vertices,
and, for every size $n$, these random graph models are also jointly exchangeable. 
It then follows from the Aldous--Hoover and graphon theory, 
as well as the fact that they are not dense, that 
these random graph models are not projective.

While dropping projectivity allowed for sparse random graph models,
the lack of projectivity
complicates the statistical applicability of these models.
At the very least, the interpretation of the aforementioned consistency results is not straightforward.
Indeed, these models are usually understood to
generate the size $n$ graphs independently of each other. Even an
adaptation of these models designed to impose some consistency between
datasets of different size seems inappropriate for modeling
data observation as, for instance, every time a new vertex is observed some
fraction of the edges already in the graph will be randomly deleted.

\subsection{Models from symmetries}

Up until this point, we have focused on very general desiderata
for random graph models. Merely requiring
sparsity and projectivity, however, does not alone lead to a tractable
class of models. 
Indeed, without any restrictions on the model, data will convey no information as to the process that gave rise to it.
To enable statistical inference, it is necessary to make some
structural assumptions on the parametrization of the random graph model.
At the same time, we want a flexible model to serve as the foundation 
of a broadly applicable framework for the statistical analysis of network data, and so we want to  impose 
as few assumptions as possible.

A general approach towards identifying large tractable families of distributions
is to consider the class of all distributions satisfying a particular invariance.
The structure of such invariant classes can be understood in general terms using very general results on ergodic decompositions,
or, in some cases, via explicit characterizations given by 
so-called representation theorems.
Both (dense) exchangeable graphs and KEGs are examples of such families, 
but to clarify the idea of defining a class of models by an invariance principle,
we will review a fundamental class of examples: the exchangeable sequences.
(The following development owes much to \citep{Orbanz_Roy_Exchangeable_Struct}, where the reader can find more details.)

Consider the classical setting of statistical inference:
a sequence of real-valued measurements $x_1,\dotsc,x_n$ are made of a system in some unknown configuration,
and this sequence is modeled as a realization from some unknown distribution $\mu_n \in \ProbMeasures(\Reals^n)$.
If, in principle, we could have made any number of measurements, 
then there exists a sequence of distributions $\mu_{1}, \mu_{2},\dotsc$ that are projective with respect to the maps $f_{m,n}$ that take length-$n$ sequences to their length-$m$ prefixes.
It follows from general results in probability theory that there exists an infinite sequence $X_1,X_2,\dotsc$ of 
random variables such that $\mu_{n}$ is the distribution of $(X_1,\dotsc,X_n)$. 
Therefore, we are modeling observed length-$n$ sequences $(x_1,\dotsc,x_n)$ as realizations of prefixes $(X_1,\dotsc,X_n)$ of the infinite random sequence $(X_1,X_2,\dotsc)$.
Let $\mu$ be the unknown distribution of the infinite sequence. %

Without making any further assumptions, it would seem that $\mu$ is an unknown element of the 
space $\ProbMeasures(\Reals^\infty)$ of all distributions on infinite sequences of real numbers.
However, a finite prefix of a realization drawn from an arbitrary element $\mu \in \ProbMeasures(\Reals^\infty)$ 
does not convey any information about the generating process $\mu$.
However, if we assume that the infinite sequence of random variables $X_1,X_2,\dotsc$ is \defnphrase{exchangeable}, i.e., 
\[
(X_{1},\dotsc,X_n ) \equaldist (X_{\sigma(1)},\dotsc,X_{\sigma(n)})
\]
for every $n \in \Nats$ and every permutation $\sigma$ of $[n]=\{1,\dotsc,n\}$, 
then, 
by de~Finetti's representation theorem~\citep{DeFinetti:1930,DeFinetti:1937,Hewitt:Savage:1955},
the random variables $X_1,X_2,\dotsc$ are conditionally i.i.d., i.e.,
there exists a probability measure $\mathcal P$ on the space $\ProbMeasures (\Reals)$ of probability measures on $\Reals$
such that
\begin{align}
  M & \dist \mathcal P \\
  X_{1},X_{2},\dotsc \given M & \distiid M.
\end{align}
We can express the distribution $\mu$ in terms of $\mathcal P$:
For a distribution $m$ on $\Reals$, let $m^\infty$ be the distribution of an infinite i.i.d.-$m$ sequence.  Then
\[\label{mixtureiid}
\mu(B) = \int_{\ProbMeasures(\Reals)} m^{\infty}(B) \, \mathcal P(\intd m), 
 \qquad \text{ for measurable $B \subseteq \Reals^\infty$.}
\]
The distribution $\mu$ is uniquely determined by $\mathcal P$, and vice versa.  
From \cref{mixtureiid}, we can see that the space of distributions of exchangeable sequences is a convex set.  
It is known that  every such distribution can be written as a unique mixture of the infinite product measures of the form $m^\infty$, which are the extreme points.  These extreme points are precisely the \defnphrase{ergodic measures}.

The statistical utility of exchangeability 
is obvious: it follows from the disintegration theorem~\citep[][Thm.~4.4]{Kallenberg:2001}
and the law of large numbers that
\[
M(A) = \lim_{n\to\infty} \frac 1 n \sum_{j=1}^n 1(X_j \in A)\ \text{a.s.}
\]
On the other hand, even an infinite realization $(x_1,x_2,\dotsc)$ gives no information about $\mathcal P$.  
For this reason, in a statistical setting, 
in addition to assuming that $\nu$ is an element in the space of distributions of exchangeable sequences,
we assume that $\nu$ is ergodic, i.e., 
$\nu$ is an unknown element in the space of distributions of i.i.d.\ sequences.
Since every $\nu$ has the form $m^\infty$ for some probability measure $m$ on $\Reals$, 
it follows that the natural parameter space is the space $\ProbMeasures(\Reals)$, and our model is $\mu_{n,\phi} = \phi^n$.

The statistical utility of exchangeability is not merely a matter of
theoretical convenience; the vast majority of statistical practice falls
under the remit of this framework. Inference of
the kind taught in introductory statistics courses is recovered by
restricting $\mathcal P$ to have support only on families of models
with finite dimensional parameterizations, e.g., the normal
distributions. The case where $\mathcal P$ has support on distributions without
finite dimensional parameterizations are so called non-parametric models,
of which there are many practical examples.

It is worth emphasizing that although de Finetti's representation theorem
is often characterized as a
justification for the use of independence in Bayesian modeling,
for our purposes the deeper point is that assuming a probabilistic
symmetry characterizes the primitive of random sequence models ($M$, a
probability distribution on $\Reals$) and gives a simple
generative recipe for the data in terms of this primitive. It is this later
perspective that is paralleled in the derivation of the KEG model.

\subsection{Models for graphs from symmetries}

We have seen how the assumption that an idealized infinite sequence of observations is exchangeable leads to a considerable simplification
of the space of distributions under consideration.  Moreover, it is clear that finite samples can be used to make inferences about the generating process.
We now turn to related results for networks. In particular, we derive the traditional exchangeable graph model from  exchangeability and then connect it to the Kallenberg exchangeable graph model.

Consider a partial observation of a network:
an array of measurements $x_{i,j}$, for $1 \le i,j \le n$, are made between $n$ entities numbered from $1$ to $n$.  We write $x_{i,j} = 1$ if a link exists between $i$ and $j$, and write $x_{i,j} = 0$ otherwise. We will assume the relationship is symmetric, i.e., $x_{i,j} = x_{j,i}$ and that no entity links to itself, i.e., $x_{i,i} = 0$.  In other words, our data is a simple graph over $n$ vertices, and we can model it as a realization from some distribution $\mu_n \in \ProbMeasures(\{0,1\}^{n \times n})$ concentrating on symmetric arrays with zeros along the diagonal.
If, in principle, we could have collected data on any number of entities, 
then there exists a sequence of distributions $\mu_{1}, \mu_{2},\dotsc$ that are projective with respect to the maps $f_{m,n}$ that take $n \times n$ arrays to their leading $m \times m$  subarrays.
Again, from general results in probability theory, there exists an infinite array of random variables
$X_{i,j}$, for $i,j \in \Nats$, 
such that $\mu_{n}$ is the distribution of $(X_{i,j};\, i,j \le n)$. 
Therefore, we model observed $n \times n$ adjacency matrices $(x_1,\dotsc,x_n)$ as realizations of prefixes $(X_{i,j};\, i,j \le n)$ of the infinite adjacency matrix $(X_{i,j};\, i,j \in \Nats)$.
Let $\mu$ be the distribution of the infinite array matrix.

Let us now consider probabilistic symmetries on this infinite idealized network observation.
The class of exchangeable sequences has a literal---if  na\"ive---counterpart in the graph setting: 
the class of edge-exchangeable graphs. 
The assumption that the edges are exchangeable is the assumption that
\[\label{edgeexch}
(X_{i,j};\, i,j \le n ) \equaldist (X_{\sigma(i,j)};\, i,j \le n),
\]
for every $n \in \Nats$ and every permutation $\sigma$ of $[n] \times [n]$ that is symmetric, i.e., $\sigma(i,j) = (i',j')$ if and only if $\sigma(j,i) = (j',i')$. 
This assumption is too severe, however, because it is simply exchangeability of a sequence in disguise.  

To see this, let $\Nats_2$ be the set of pairs $(i,j) \in \Nats^2$ such that $i < j$
let  $\iota : \Nats \to \Nats_2$ be an arbitrary bijection, and  define $Y_n = X_{\iota(n)}$.  
Then \cref{edgeexch} implies that the sequence of random variables $Y_1,Y_2,\dotsc$ are exchangeable and so they are conditionally i.i.d.  But then the edges $X_{\iota(n)}$, for $n \in \Nats$, are also conditionally i.i.d.
Therefore, there exists a random variable $p$ in $[0,1]$ such that, conditioned on $p$, 
the edges $X_{i,j}$
are i.i.d.\ and each edge appears 
with probability $p$.
This is none other than the \ErdosRenyi\ model with a random edge probability.  The class of ergodic measures in this case is precisely the \ErdosRenyi\ model.

\begin{figure}
  \begin{centering}
    \includegraphics[width=0.5\linewidth]{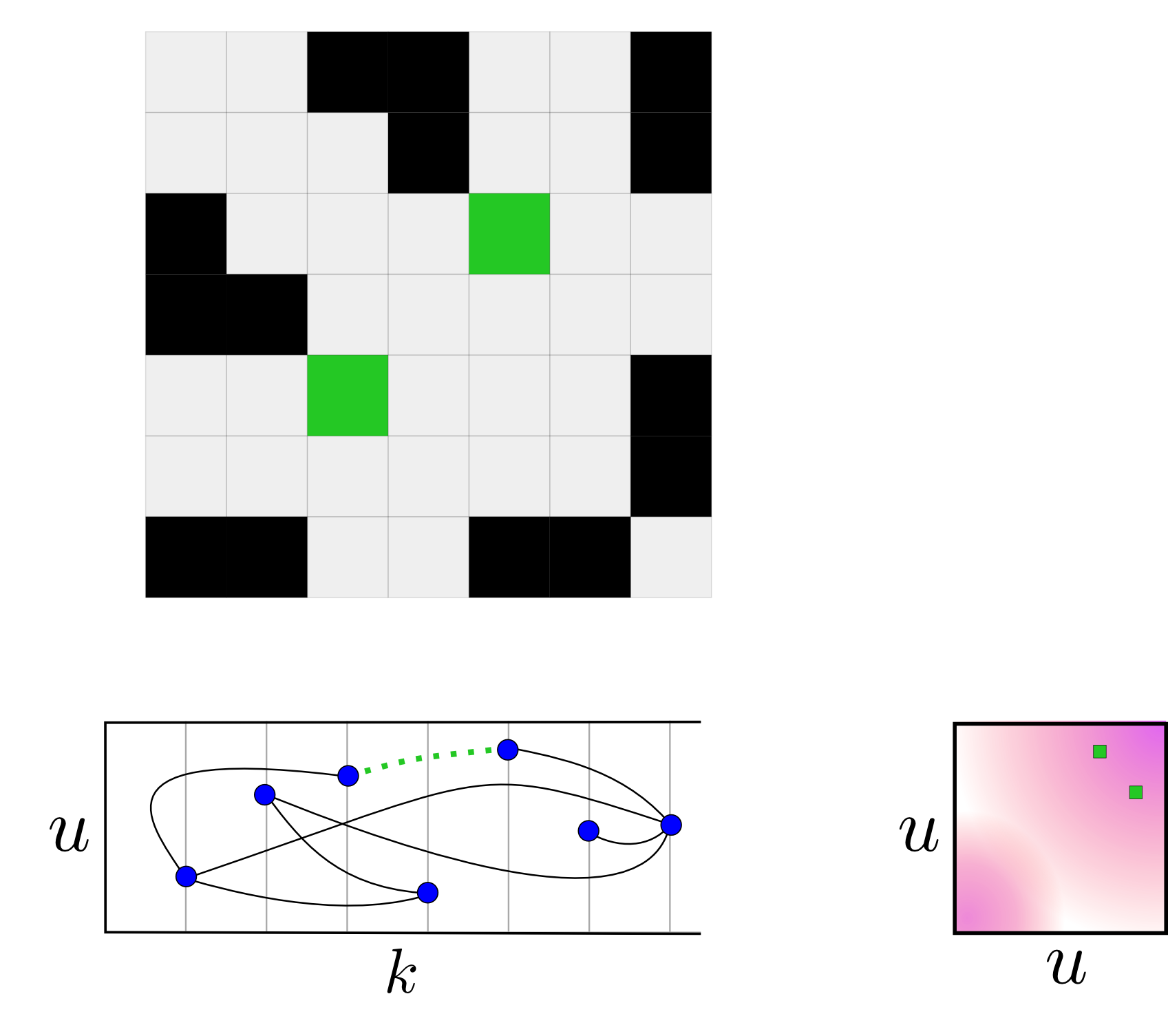}
    \par\end{centering}

  \caption{\label{fig:Graphon_model}Graphon random graph model. In the
    jointly exchangeable array setting a random graph model is
    characterized by a (potentially random) 
    symmetric measurable function $W: [0,1]^2 \to [0,1]$ called a graphon.
    An example graphon is depicted as a magenta heatmap (lower
    right). Conditional on $W$, a random graph of size $n$ is generated by
    independently assigning to each vertex $k \in \{1,\dotsc,n\}$ a latent random variable
    $U_{k}\dist\uniDist (0,1)$ (values along vertical axis) and
    including each edge $(k,l)$ independently with probability $W(U_{k},U_{l})$.  
    For
    example, edge $(3,5)$ (green, dotted) is present with probability
    $W(0.72,0.9)$; the green boxes in the right square represent the
    values of $(u_{3},u_{5})$ and $(u_{5},u_{3})$. The upper left
    panel shows the graph realization as an adjacency matrix.}
\end{figure}

The natural analogue of exchangeability in the graph setting is to assume that the labels of the vertices are exchangeable.
Informally, this is the assumption that the vertex labels carry no information.
Given that we are representing an observed adjacency matrix as a prefix of  an idealized infinite symmetric binary array,
vertex-exchangeability is formalized as the requirement that distribution of the array is invariant under simultaneous permutation of its rows and columns.
More carefully, an array of random variables $X_{i,j}$ is \emph{jointly exchangeable} when
\[
(X_{i,j};\, i,j \le n)\equaldist(X_{\sigma(i),\sigma(j)};\, i,j\le n)
\]
for every $n \in \Nats$ and every permutation $\sigma$ of $[n]$.
A characterization of infinite jointly exchangeable adjacency matrices
can be easily derived from the Aldous--Hover representation theorem 
for general jointly exchangeable arrays~\cite{Aldous:1981,Hoover:1979}. 
In particular, every ergodic measures is characterized by 
a symmetric measurable function $W:[0,1]^{2}\to[0,1]$, whose diagonal is zero.
This same object was later rediscovered independently
by graph theorists as the limit object in a theory of limits of dense
graphs~\cite{Lovasz:Szegedy:2006,Lovasz:Szegedy:2007,Lovasz:2013:A}.
In this context it was named a graphon, which is the nomenclature we
use here. 
The relationship between the graphon as the defining object
for distributions of jointly exchangeable arrays and as the limit object of
dense graph theory is explained by \cite{Diaconis:Janson:2007}.
More concretely, the generative model for vertex-exchangeable graphs
is (see \cref{fig:Graphon_model})
\begin{align}
  W & \dist  \mu\\
  \{ U_{i}\}  & \distiid \uniDist [0,1]\\
  (X_{ij})\given W,U_{i},U_{j} & \distind \bern (W(U_{i},U_{j})),
\end{align}
where $\mu$ is a measure on the space of symmetric functions from the
unit square to the unit interval with zero diagonal. 
The fact that projective and
jointly exchangeable adjacency matrices cannot be sparse is a simple
consequence of this generative model and the law of large
numbers. In particular, any nondiagonal entry is one with probability $\| W \|_1$. This framework is
the exchangeable graph model, whose nomenclature is now self explanatory. Comparing the generative model
for the exchangeable graph model with the KEG generative model
(see \cref{fig:graphex_model}) makes it clear that the distinction
that allows for more general graphs in the KEG setting is that the
latent variables associated with each vertex are not independent, 
and the sizes of the graphs are random.

It is possible to construct a sparse and projective random graph model 
if we drop the requirement that the arrays of each size $n \in \Nats$ be exchangeable.
For example, the preferential
attachment model of \cite{Barabasi:Albert:1999} can be understood in these terms, although historically it was developed independently of these concerns for the special purpose of giving a mechanism of graph \emph{growth} that leads to power law behavior in the degree distribution. Ad hoc models of this kind tend to fail to capture certain
key elements of real-world network structure. For instance, as shown by \cite{Berger:Borgs:Chayes:Saberi:2014},
the limiting local structure of preferential attachment graphs is a tree, and so these networks would be pathological models of social networks, which exhibit homophily.

\subsection{Random graphs as random measures}

\begin{figure}[t]

  \begin{centering}
    \includegraphics[width=0.5\linewidth]{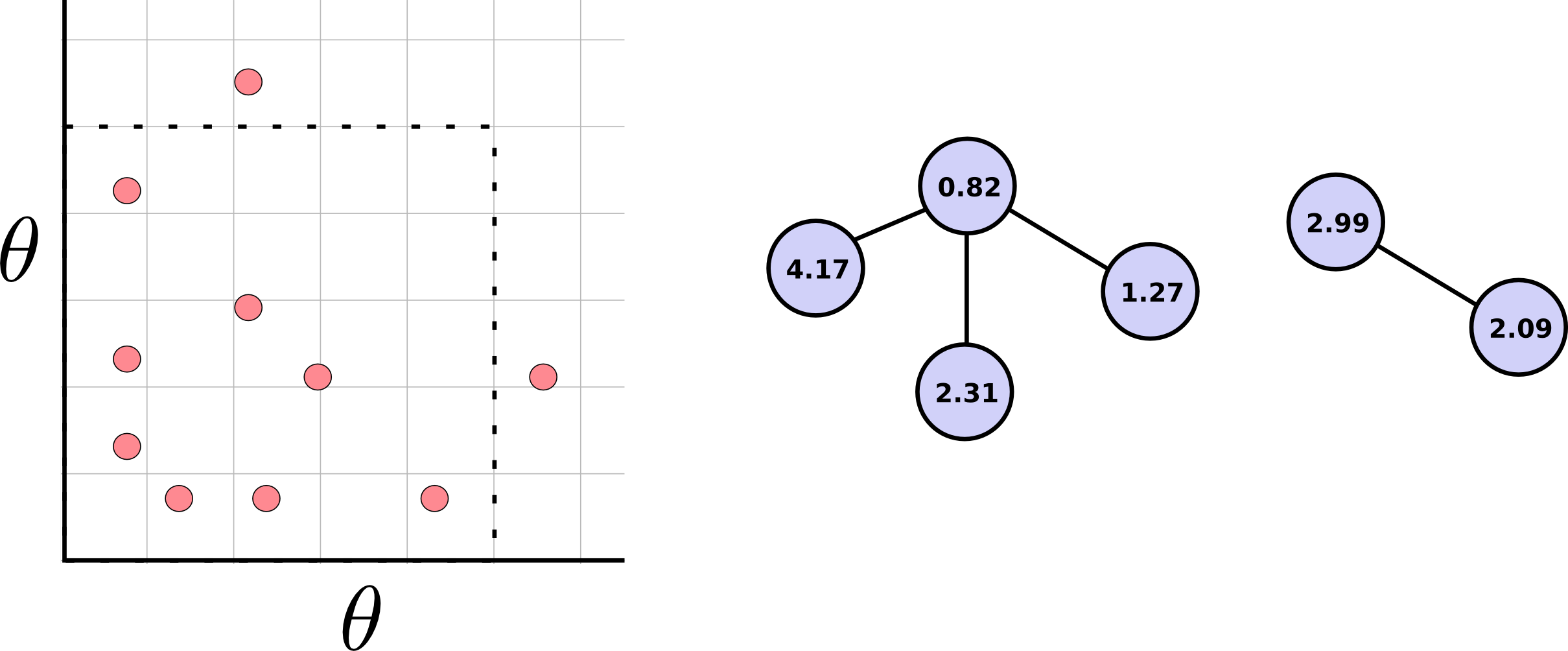}
    \par\end{centering}

  \caption{\label{fig:random_graph_measure_corresp}Random graphs as
    point processes.  Random point processes on $\NNReals^{2}$
    correspond to infinite random graphs, with finite subgraphs given
    by restricting the point process to a finite square. Points of the
    process correspond to graph edges and the vertex structure is
    deduced from the edge structure.  Pictured is a realization of a
    point process and the realization of the random graph that
    corresponds to truncating at $\theta=5$. }
\end{figure}

The key ingredient for generalizing the exchangeable graph model is a
correspondence
between random graphs and symmetric simple point processes due to
Caron and Fox~\citep{Caron_Fox_CRM_Graphs} (see
\cref{fig:random_graph_measure_corresp}).
Again, restricting ourselves to simple
graphs for simplicity of presentation, 
the edge set of a random graph is a random finite or countable collection of 
tuples $(x,y)\in\NNReals^{2}$, and the vertex set is the set of those
real numbers $x$ such that $x$ participates in at least one edge.
Concretely, the random graph is represented by a simple point process $G$ on $\NNReals^2$
containing a point
$(x,y)$ iff there is an edge $(x,y)$ in the random graph.

It will be mathematically convenient to represent simple point processes by simple random measures, i.e., purely atomic random measures whose atoms all have mass one. In this case, each atom in the simple random measure represents a point of the point process.
Having made this choice, the idealized infinite observation in this setting is the infinite point process $G$,
and finite observations are the restrictions $G_t = G( \cdot \, \cap [0,t]^2)$, 
for $t \in \NNReals$, of the infinite point process $G$ to the 
bounded square subsets $[0,t]^2 \subset \NNReals^2$ containing the origin.
The distribution of these restrictions of $G$ are automatically projective
with respect to the maps $f_{s,t}$ that takes a measure on $[0,t]^2$ to its restriction on $[0,s]^2$.
In contrast to the exchangeable graph model,
the KEG model has a continuously indexed size parameter and
the number of vertices in each finite restriction $G_t$
is itself a random quantity.

It is important to note that the graph corresponding to 
the restriction $G_s$ to $[0,s]^2$
has as its vertex set only those vertices $x \in [0,s]$ that appear
in some edge $(x,y)$ where $y \in [0,s]$.
In particular, there will, in general, be vertices in $[0,s]$ that appear
for the first time in a restriction $[0,t]$, for $t > s$.
This is an essential property of this representation, and is
the way that the seeming equivalence between exchangeability and density can be
relaxed. The point labeled 2.7 in \cref{fig:graphex_model} provides a concrete
example of this phenomena.

As observed by Caron and Fox, when random graphs are represented as point processes,
vertex-exchangeability corresponds to joint exchangeability for \emph{random measures}.
  Formally, a random measure $\xi$ on $\NNReals^2$ is \emph{jointly exchangeable} when
\[
\xi\equaldist\xi\circ(f\otimes f)^{-1}
\]
for every measure preserving transformation $f:\NNReals\to\NNReals$,
where $\otimes$ is the tensor product.  
This probabilistic symmetry was introduced by Aldous, who also
conjectured a concrete representation theorem~\citep[][Conj. 15.15]{Aldous:1983},
later established rigorously by 
Kallenberg~\cite{Kallenberg_Random_Meas_Plane,Kallenberg:2005}.  We will refer to the representation theorem as the Kallenberg representation theorem.

We now describe the Kallenberg exchangeable graph model plainly:  
It is the random graph model 
that arises from the symmetry of joint exchangeability of symmetric simple point processes on $\NNReals^2$, when these structures are interpreted as the edge sets of random graphs.
We give a representation theorem for these structures
via a straightforward application of Kallenberg's representation theorem 
in the specific context of symmetric simple point processes on $\NNReals^2$.
From this result, we see that every ergodic measure is determined by a triple $(\IsoF,\StarF,W)$, which we call a graphex.
From a statistical standpoint, the graphexes are the natural parameters,
and every random graph is seen to arise via the corresponding generative process (\cref{fig:graphex_model}).
The KEG model is projective, exchangeable,  and admits sparse graphs, 
thereby providing a statistical framework for network analysis that avoids some of 
the pitfalls of other random graph models. 
Both the traditional exchangeable graph model 
and the Caron--Fox model are special cases, and so the KEG model can be seen as a generalization and unification of these models.%

\section{Examples\label{sec:Examples}}

The aim of this section is to work through the details of several
informative examples to build intuition for the structure of the
Kallenberg exchangeable graph models we consider here. 
We focus on those graphexes where $\IsoF=\StarF=0$, and so we will
refer
to $W$ as the graphex without any risk of confusion.
We are particularly
interested in the sparsity of these graph
models. \cref{thm:expected_edges} establishes that (ignoring self
edges) for all random graphs $G_\nu$ generated by graphex $W$ it holds
that $\expect{e_{\nu}}=\frac{1}{2}\nu^{2}\norm W_{1}$; i.e., the
expected number of edges scales as $\nu^{2}$ in all cases. Intuitively
then we expect the sparsity of a random graph model to be determined
by $\expect{v_{\nu}}=\nu\int_{\NNReals} 1-e^{-\nu\mu_{W}(x)}\intd x$
(from \cref{thm:expected_vertices}, ignoring self edges). This suggests
that the slower $\mu_{W}(x)=\int_{\NNReals}W(x,y)\intd y$ decays the
sparser the graph will be, an intuition that is borne out by the
examples of this section.

\subsection{Graphon models}
\label{sec:graphon_models}

The above argument suggests that the most densest graphs will correspond to those $W$ 
that are compactly supported. Let $\widetilde{W}:[0,1]^2\to[0,1]$ be a graphon and consider the graphex
given by the dilation
\[
W(x,y)=\begin{cases}
  \widetilde{W}(x/c,y/c) & x\le c,y\le c\\
  0 & \mbox{otherwise.}
\end{cases}
\]
In this case, points $(\theta,\vartheta)\in\PP$ of the latent Poisson
process will fail to connect to an edge if $\vartheta>c$, and so such points
 they never
participate in the graph and can be discarded. This means that for finite size graph
$G_{\nu}$ given by restricting $\theta\le\nu$ the relevant underlying
process is the unit rate Poisson process on $[0,\nu]\times[0,c]$.
The generative model for the graph can be expressed as:
\begin{align}
  N_{\nu} & \dist \poiDist(c\, \nu)\\
  \{ \theta_{i}\} \given N_{\nu} & \distiid \uniDist [0,\nu]\\
  \{ \vartheta_{i}\} \given N_{\nu} & \distiid \uniDist [0,1]\\
  (\theta_{i},\theta_{j}) \given\widetilde{W},\vartheta_{i},\vartheta_{j}
  & \distind \bern (\widetilde{W}(\vartheta_{i},\vartheta_{j})).
\end{align}
A little thought shows that this is just a trivial modification of the
graphon model. Instead of indexing the family of graphs by the number
of vertices ($\Nats$) we now index them by the continuous parameter $\nu$
and have $\poiDist(c\, \nu)$ candidate vertices at each stage.  The vertices
now have \iid uniform labels instead of the integer labels of the
traditional graphon model and vertices are only included if they
connect to at least one edge. The critical components of the graphon
model structure are unchanged: the primitive is still the graphon
$\widetilde{W}:[0,1]^{2}\to[0,1]$, the conditional independence of the
edges is the same, the latent variables are independent, and these
graphs are necessarily asymptotically dense (or empty). This is the
sense in which the graphon model is a special case of the graphex
model derived in this paper.

In fact, these are the only dense KEGs arising from (integrable) graphexes:
\cref{thm:dense_iff_compact} shows that $G$ is dense \iff the
generating (integrable) graphex has compact support.

\subsection{Slow Decay}

We next consider a graphex with tails that go to 0 slowly:
\begin{align}
  W(x,y) & = \begin{cases}
    0 & x=y,\\
    (x+1)^{-2}(y+1)^{-2} & \mbox{otherwise,}
  \end{cases}
\end{align}
where the condition $W(x,x)=0\ \forall x\in\NNReals$ forbids self
edges. In this case $\mu_{W}(x)=\frac{1}{3}(x+1)^{-2}$ and by
\cref{thm:expected_vertices}
\begin{align}
  \expect{v_{\nu}} & = \nu(\sqrt{\pi}\sqrt{\nu/3}\mbox{erf}(\sqrt{\nu/3})+e^{-\nu/3}-1)\\
  & \sim \sqrt{\frac{\pi}{3}}\nu^{3/2}\asympLim{\nu}.
\end{align}
By \cref{thm:expected_deg_k} the number of vertices with degree $k$ has
expectation:
\begin{align}
  \expect{\numDegNu k} & = \frac{\nu^{k+1}}{k!}(\frac{1}{3})^{k}\int_{1}^{\infty}x^{-2k}e^{-\frac{1}{3}\nu x^{-2}}\intd x\\
  & = \frac{\nu^{k+1}}{k!}(\frac{1}{3})^{k}\int_{0}^{1}x^{2(k-1)}e^{-\frac{1}{3}\nu x^{2}}\intd x\\
  & = \frac{\Gamma(-\frac{1}{2}+k)-\Gamma(-\frac{1}{2}+k,\frac{\nu}{3})}{2\sqrt{3}k!}\nu^{3/2}\\
  & \sim
  \frac{\Gamma(-\frac{1}{2}+k)}{2\sqrt{3}k!}\nu^{3/2}\asympLim{\nu}.
\end{align}
By \cref{thm:deg_dist_theorem} it follows that the degree $D_{\nu}$ of
a uniformly selected vertex of $G_{\nu}$ satisfies
\[
\Pr(D_{\nu} = k\given
G_{\nu})\convPr\frac{\Gamma(-\frac{1}{2}+k)}{2\sqrt{\pi}k!}\asympLim{\nu},
\]
so in particular a randomly selected vertex of $G_{\nu}$ will have
finite degree even in the infinite graph limit. For large $k$
\[
\frac{\Gamma(-\frac{1}{2}+k)}{2\sqrt{\pi}k!}\sim
k^{-\frac{3}{2}}\asympLim k,
\]
so this is an example of a random graph model with power-law degree
distribution.  Note that, in the limit, while the degree of a randomly
chosen vertex is finite almost surely, it is infinite in expectation.

\subsection{Fast Decay}

Next we consider a graphex with quickly decaying tails. Let
\[
W(x,y)=\begin{cases}
  0 & x=y\\
  e^{-x}e^{-y} & \mbox{otherwise.}
\end{cases}
\]
Then $\mu(x)=e^{-x}$ and so by \cref{thm:expected_vertices}
\begin{align}
  \expect{v_{\nu}} & = \nu\int_{\NNReals}1-e^{-\nu e^{-x}}\intd x\\
  & = \nu\int_{0}^{1}\frac{1}{x}(1-e^{-\nu x})\intd x\\
  & = \nu(\gamma+\Gamma(0,\nu)+\log(\nu))\\
  & \sim \nu\log\nu\asympLim{\nu}.
\end{align}
As expected, the rapidly decaying graphex gives rise to a graph that
is much more dense than one from the slowly decaying graphex.

By \cref{thm:expected_deg_k} the number of vertices with degree $k$ has
expectation:
\begin{align}
  \expect{\numDegNu k} & = \frac{\nu^{k+1}}{k!}\int_{0}^{\infty}e^{-kx}e^{-\nu e^{-x}}\intd x\\
  & = \frac{\nu}{k!}(\Gamma(k)-\Gamma(k,\nu))\\
  & \sim \frac{\nu}{k} \asympLim{\nu}.
\end{align}
so that for fixed $k$ only a vanishing fraction of the vertices will have
degree $k$ as $\nu\to\infty$. More precisely, since
$\sum_{k=1}^{\nu^{\beta}}\frac{\nu}{k}\sim\beta\nu\log\nu
\asympLim{\nu}$ we have by \cref{thm:deg_dist_theorem} that for
$0<\beta<1$
\[
P(D_{\nu}\le\nu^{\beta})\convPr\beta\asympLim{\nu}
\]
where $D_{\nu}$ is a random vertex of $G_{\nu}$.

\subsection{Caron and Fox}

As already alluded to, the family of random graph models considered by
Caron and Fox in \cite{Caron_Fox_CRM_Graphs} is a special case of the
one considered here. Indeed, in their paper they prove their model
satisfies joint exchangeability when considered as a random measure
and use Kallenberg's representation theorem to derive some model
properties. Nevertheless, the connection is opaque
because their model is constructed from products of
completely random measures and 
they cast their model in terms of L\'evy process intensities. If the
$\boldsymbol{\theta}\times\boldsymbol{\theta}$ measure they had studied
had been a product of completely random measures, that model would have
corresponded to a graphex of the form $W(x,y)=f(x)f(y)$.  Instead,
they actually consider a measure on
$\boldsymbol{\theta}\times\boldsymbol{\theta}$ given by using the
product of completely random measures as a base measure for a Cox
process. This gives rise to a directed multigraph which is then
transformed into a simple graph by including edge
$\{\theta_{i},\theta_{j}\}$ if and only if there is at least one
directed edge between $\theta_{i}$ and $\theta_{j}$. A little algebra
shows this model corresponds to the graphex 
\[
W(x,y)=\begin{cases}
  1-\exp(-g(x)g(y)) & x=y\\
  1-\exp(-2g(x)g(y)) & x\neq y
\end{cases}
\]
where $g(x):\NNReals\to\NNReals$. Caron and Fox derive this expression in their
paper, and give $g$ in terms of the intensity of the defining L\'evy process.

\section{Representation Theorem for Random Graphs represented by Exchangeable
  Symmetric Simple Point Processes}
\label{sec:Representation-Theorem}

We now turn to giving formal statements of our construction and
proving the representation theorem at the heart of the paper. In fact,
this mostly amounts to translating Kallenberg's representation theorem
for jointly exchangeable random measures on $\NNReals$ to the random graph setting.

The central objects of study here are undirected, unweighted graphs whose vertices are labeled with values in $\NNReals$.
For a graph $G$, we will write \defnphrase{$\vertexset G$} and \defnphrase{$\edgeset G$}
to denote the set of vertices and edges, respectively.  
We begin by formalizing the idea of a graph represented by a measure.

\begin{defn}
An \defnphrase{adjacency measure} is a locally finite symmetric simple measure on $\NNReals^2$.
The \defnphrase{$\nu$-truncation} of an adjacency measure $\xi$ is the adjacency measure $\xi( \cdot\, \cap [0,\nu]^2 )$ obtained by restricting $\xi$ to $[0,\nu]^2$.
\end{defn}

\begin{defn}
Let $G$ be a simple graph, possibly with loops,
whose edge set $e(G)$ is a locally finite subset of $\NNReals^2$.
Then the \defnphrase{adjacency measure of $G$} is the adjacency measure $\sum_{(x,y) \in e(G)} \delta_{(x,y)}$.
\end{defn}

Note that the adjacency measures of a graphs $G$ and $G'$ coincide if and only if their edge sets do.  In particular, vertices that do not participate in an edge are ``forgotten''.
We will be interested in the smallest graph corresponding to an adjacency measure $\xi$,
which is necessarily the graph with the same edge set and no isolated vertices.
(See \cref{fig:random_graph_measure_corresp} for an illustration.)

\begin{defn}
Let $\xi = \sum_{i < \kappa} \delta_{e_i}$ be an adjacency measure, where $\kappa \in \NNInts \cup \{\infty\}$ and $e_1,e_2,\dotsc$ is a sequence of distinct elements of $\NNReals^2$.  
Then the \defnphrase{simple graph associated with $\xi$}
is the graph $G$ whose edge set is $\{ e_i : i < \kappa \}$ and whose vertex set is
$\{ x : \exists i < \kappa\, \exists y \in \NNReals \, e_i = (x,y) \}$. 
\end{defn}

\begin{rem}
  This correspondence extends to directed weighted graphs in an
  obvious way by dropping the requirement that the adjacency measure
  be symmetric and allowing the adjacency measure to assign a mass 
  other than one to each of its atoms; i.e., a directed weighted adjacency measure 
  is a locally finite purely atomic measure, and so would have the form
  $\xi=\sum_{ij}\omega_{ij}\delta_{(\theta_{i},\theta_{j})}$.
\end{rem}

A random adjacency measure is an (a.s.\ locally finite) symmetric simple point process.  
We will represent random graphs by their random adjacency measures, noting that only nonisolated vertices 
are captured by this representation.

Informally, we are interested in those simple random graphs embedded in $\NNReals$ whose
distributions are invariant to every relabeling of the vertices of
the random graph.
We can formalize this notion of invariance in terms of a symmetry of the corresponding adjacency measure.  We begin with a definition of exchangeability for random measures due to Aldous:
\begin{defn}
  A random measure $\xi$ on $\NNReals^{2}$ is said to be
  \defnphrase{jointly exchangeable} if, for every measure preserving
  transformation $f$ on $\NNReals$, we have
  \[
  \xi\circ(f\otimes f)^{-1}\equaldist\xi.
  \]
\end{defn}
The following result, due to Kallenberg, characterizes the space of exchangeable measures on $\NNReals^2$ as well as its extreme points:  
Let $\Lebesgue$
  denote Lebesgue measure on $\NNReals$ and let $\Lebesgue_{D}$ denote Lebesgue
  measure on the diagonal of $\NNReals^2$.
\begin{thm}[{Kallenberg
    \cite{Kallenberg:2005,Kallenberg_Random_Meas_Plane}}]
  \label{jointexchmeasure}
  A random measure $\xi$ on $\NNReals^{2}$ is jointly exchangeable
  \iff almost surely
  \begin{align}\label{KEGgen}
    \xi = &  \hspace*{3.6mm}   \sum_{i,j}f(\alpha,\vartheta_{i},\vartheta_{j},\zeta_{\{i,j\}})\delta_{\theta_{i},\theta_{j}} \\
    & +  \sum_{j,k}(g(\alpha,\vartheta_{j},\chi_{jk})\delta_{\theta_{j},\sigma_{jk}}+g^{\prime}(\alpha,\vartheta_{j},\chi_{jk})\delta_{\sigma_{jk},\theta_{j}})\\
    & +  \sum_{k}(l(\alpha,\eta_{k})\delta_{\rho_{k},\rho'_{k}}+l^{\prime}(\alpha,\eta_{k})\delta_{\rho'_{k},\rho_{k}}) \\
    & +  \sum_{j}(h(\alpha,\vartheta_{j})(\delta_{\theta_{j}}\otimes\Lebesgue)+h^{\prime}(\alpha,\vartheta_{j})(\Lebesgue\otimes\delta_{\theta_{j}})) +\beta\Lebesgue_{D}+\gamma\Lebesgue^{2}, \label{lebterms}
  \end{align}
  for some measurable function $f\ge0$ on $\NNReals^{4}$,
  $g,g^{\prime}\ge0$ on $\NNReals^{3}$ and
  $h,h^{\prime},l,l^{\prime}\ge0$ on $\NNReals^{2}$, 
  some collection of independent uniformly distributed random variables $(\zeta_{\{i,j\}})$ on $[0,1]$,
  some independent unit rate Poisson processes $\{ (\theta_{j},\vartheta_{j})\} $ and
  $\{ (\sigma_{ij},\chi_{ij})\} _{j}$, for $i\in\Nats$, on $\NNReals^{2}$
  and $\{ (\rho_{j},\rho_{j}^{\prime},\eta_{j})\} $ on $\NNReals^{3}$,
  and some independent set of random variables $\alpha,\beta,\gamma\ge0$.
   The latter can be chosen to be non-random iff $\xi$ is extreme.
\end{thm}
The task is to translate this into a statement about random graphs, or more specifically, their adjacency measures.
Because adjacency measures are purely atomic,
all terms with a Lebesgue component (\cref{lebterms}) must have measure zero. The remaining purely atomic terms underlying a jointly exchangeable random measure have the following
interpretation for adjacency measures:
\begin{enumerate}
\item
  $\sum_{i,j}f(\alpha,\vartheta_{i},\vartheta_{j},\zeta_{\{i,j\}})\delta_{\theta_{i},\theta_{j}}$:
  this term contributes most of the interesting structure for the
  random graph models. The random measure $\xi$ will be symmetric and simple if and only if 
  $f$ is a.e. $\{0,1\}$-valued and
  symmetric in its second and third arguments, for a.e.\ fixed first and fourth argument.  
  (It is clear that this can easily be strengthened to hold everywhere.) 
  This leads to the correspondence illustrated in
  \cref{fig:graphex_model}. (General $f$ could be used to model directed,
  weighted graphs in an obvious way.) The tuples
  $(\theta_{i},\theta_{j})$ are possible edges of the graph and the
  points $\theta_{i}$ are candidate vertices.
\item
  $\sum_{j,k}(g(\alpha,\vartheta_{j},\chi_{jk})\delta_{\theta_{j},\sigma_{jk}}+g^{\prime}(\alpha,\vartheta_{j},\chi_{jk})\delta_{\sigma_{jk},\theta_{j}})$:
  this term contributes stars. To see this, note that each candidate vertex $\theta_{j}$ has an
  associated Poisson process $\{\sigma_{jk}\}$.  The points are a.s.\ distinct: i.e.,
  $\{\theta_{l}\}\cap\{\sigma_{jk}\}=\emptyset$ and
  $\{\sigma_{jk}\}\cap\{\sigma_{lk}\}$ for $j\neq l$ with probability
  one. This means the candidate vertices $\{\sigma_{jk}\}$ will only ever
  participate in edges with $\theta_{j}$, hence the star
  structure. 
  The random measure $\xi$ will be a.s.\ symmetric and simple iff
  $g=g^{\prime}$ and $g$ is $\{0,1\}$-valued.
\item
  $\sum_{k}(l(\alpha,\eta_{k})\delta_{\rho_{k},\rho'_{k}}+l^{\prime}(\alpha,\eta_{k})\delta_{\rho'_{k},\rho_{k}})$:
  this term contributes isolated edges.  To see this, note that, with probability one,
  $\{\rho_{k}\}\cap\{\rho_{k}^{\prime}\}=\emptyset$ and these
  candidate vertices do not coincide with any other candidate vertices
  (e.g., $\{\rho_{k}\}\cap\{\theta_{l}\}=\emptyset$). This means that
  if $(\rho_{i},\rho_{j})$ is an edge of the graph then with
  probability 1 $(\rho_{i},x)$ will not be an edge for any
  $x\in\NNReals$.  Again, the random measure $\xi$ will be a.s.\ symmetric and simple iff
  $l = l'$ and $l$ is $\{0,1\}$-valued.
\end{enumerate}

The following theorem characterizes the space of exchangeable adjacency measures
as well as its extreme points: 
\begin{thm}[{Random graph representation}]
  \label{thm:graphex_rep_theorem}
  Let $\xi$ be a random adjacency measure.  
  Then $\xi$ is jointly exchangeable \iff almost surely
  \begin{align}
    \xi = 
    & \hspace*{3.6mm} 
            \sum_{i,j}1[W(\alpha,\vartheta_{i},\vartheta_{j}) \le \zeta_{\{i,j\}}]\delta_{\theta_{i},\theta_{j}}\\
    & +  \sum_{j,k} 1[ \chi_{jk} \le \StarF(\alpha,\vartheta_{j}) ](\delta_{\theta_{j},\sigma_{jk}}+\delta_{\sigma_{jk},\theta_{j}})\\
    & +   \sum_{k} 1[\eta_{k} \le \IsoF(\alpha) ](\delta_{\rho_{k},\rho'_{k}}+\delta_{\rho'_{k},\rho_{k}}),
  \end{align}
  for 
  some measurable function
     $\StarF : \NNReals^2 \to \NNReals$,
     $\IsoF : \NNReals \to \NNReals$,
     $W : \NNReals^3 \to [0,1]$, where $W(a,\cdot,\cdot)$ is symmetric for every $a \in \NNReals$;
  some collection of independent uniformly distributed random variables $(\zeta_{\{i,j\}})$ in $[0,1]$;
  some independent unit rate Poisson processes $\{ (\theta_{j},\vartheta_{j})\} $ and
  $\{ (\sigma_{ij},\chi_{ij})\} _{j}$, for $i\in\Nats$, on $\NNReals^{2}$
  and $\{ (\rho_{j},\rho_{j}^{\prime},\eta_{j})\} $ on $\NNReals^{3}$;
  and an independent random variable $\alpha\ge0$.
   The latter can be chosen to be non-random iff $\xi$ is extreme.
\end{thm}
  The second term of this measure corresponds to stars centered at the
  points $\{ \theta_{j}\} $ and the third term corresponds to isolated edges
  that do not connect to the rest of the graph.
\begin{proof}
  Most of this result is immediate from the text preceding the
  theorem.  One direction of the correspondence is immediate: the
  random measure $\xi$ is obviously jointly exchangeable.

  In the other direction, let $f$, $\alpha$,
  $\{\theta_i,\vartheta_i\}$, and $\{\zeta_{\{i,j\}}\}$ be as in
  \cref{jointexchmeasure}, and let
  \[
  \xi_{\{i,j\}} \defas
  f(\alpha,\vartheta_{i},\vartheta_{j},\zeta_{\{i,j\}}),
  \]
  which is well-defined because $f$ is symmetric in its second and
  third arguments.  Define $W : \NNReals^3 \to \NNReals$ by 
  \[
  W(a,t,t') 
  = \Lebesgue \{ z \in [0,1] \st f(a,t,t',z) = 1 \}
  = \Lebesgue f(a,t,t',\cdot\,) ,
  \]
  and write $W_a$ for $W(a,\cdot,\cdot)$.  Note that $W_a$ is
  symmetric.  Let $\mathcal F \defas
  \sigma(\alpha,\{(\vartheta_i,\theta_i)\}_{i \in \Nats} )$.  Then the
  random variables $\xi_{\{i,j\}}$, for $\{i,j\} \in \NatSubs 2$,
  are
  independent given $\mathcal F$ and satisfy
  \[
  \EE [ \xi_{\{i,j\}} | \mathcal F ] \equalas W_\alpha(\vartheta_i,
  \vartheta_j).
  \]
  Let $\{\zeta'_{\{i,j\}}\}$ be an i.i.d.\ uniform array on $\NatSubs
  2$, independent from $\mathcal F$, and define, for $\{i,j\} \in
  \NatSubs 2$,
  \[
  \xi'_{\{i,j\}} = 1(W_\alpha(\vartheta_i,\vartheta_j) \le
  \zeta'_{\{i,j\}}).
  \]
  Then it is clear that
  \[
  ( \alpha, ((\theta_i,\vartheta_i)_{i \in \Nats}),
  (\xi'_{\{i,j\}})_{\{i,j\} \in \NatSubs 2}) \equaldist ( \alpha,
  ((\theta_i,\vartheta_i)_{i \in \Nats}), (\xi_{\{i,j\}})_{\{i,j\} \in
    \NatSubs 2})
  \]
  and so, by a transfer argument \citet[][Cor~6.11]{Kallenberg:2001},
  there exists an \iid uniform array $\{\zeta''_{\{i,j\}}\}$ on
  $\NatSubs 2$ independent also from $\mathcal F$ such that
  \[
  \xi_{\{i,j\}} \equalas 1(W_\alpha(\vartheta_i,\vartheta_j) \le
  \zeta''_{\{i,j\}}).
  \]
Similarly, letting $g$ and $l$ be as in \cref{jointexchmeasure},
define
\[
\StarF(a,t) \defas \Lebesgue \{ z \in \NNReals \st g(a,t,z) = 1 \} = \Lebesgue g(a,t,\cdot\,)
\]
and
\[
\IsoF(a) \defas \Lebesgue \{ z \in \NNReals \st l(a,z) = 1 \} = \Lebesgue l(a,\cdot\,).
\]
A similar argument to above can be used to show that the terms involving $\StarF$ and $\IsoF$ agree with their 
counterparts in \cref{jointexchmeasure}.
\end{proof}

From the representation theorem, we learn
that the extreme members, from which all other can be recovered as mixtures,
are naturally defined in terms of a triple $(\IsoF,\StarF,W)$,
where $I \in \NNReals$ and  $S : \NNReals \to \NNReals$ and $W :\NNReals^2 \to \NNReals$ are measurable, and $W$ is symmetric.

In general, an exchangeable simple point process $\xi$ of the form above 
may not be finite when restricted to a finite region $[0,t]^2$.  We want finite restrictions of the adjacency measure to correspond to finite size observations, and so we must isolate conditions on the triple $(\IsoF, \StarF, W)$ so that the random measure is a.s. finite on bounded sets.
The following result, due to Kallenberg, gives necessary and sufficient conditions for a jointly exchangeable measure to be a.s.\ locally finite.
\begin{thm}[local summability {\cite[][Prop.~9.25]{Kallenberg:2005}}]
Let $\xi$ be as in \cref{jointexchmeasure}, 
write $\hat f = f \wedge 1$, 
and let
\[
f_1 = \Lebesgue_{23}^2 \hat f, \qquad
f_2 = \Lebesgue_{13}^2 \hat f, \qquad
g_1 = \Lebesgue_{2} \hat g,
\]
where $\Lebesgue_{23}^2$ denotes two-dimensional Lebesgue measure in the second and third coordinates, 
and similarly for $\Lebesgue_{13}^2$ and $\Lebesgue_2$.
For fixed $\alpha$, the random measure $\xi$ is  a.s.\ locally finite iff these five conditions are fulfilled:
\begin{enumerate}[label=\normalfont(\roman*)]
\item $\Lebesgue(\hat l + \hat h + \hat h') < \infty$,
\item $\Lebesgue(\hat g_1 + \hat g_1') < \infty$,
\item $\Lebesgue\{f_i = \infty \} = 0$ and $\Lebesgue \{f_i > 1\} < \infty$ for $i=1,2$,
\item $\Lebesgue^2[\hat f; f_1 \vee f_2 \le 1] < \infty$,
\item $\Lebesgue \hat l' + \Lebesgue_D \Lebesgue \hat f < \infty$.
\end{enumerate}
\end{thm}

(Note that we have corrected a typo in part (iv), where the integral was taking w.r.t.\ $\Lebesgue$ not $\Lebesgue^2$.)
The consequences for adjacency measures is as follows:

\begin{thm}[locally finite graphex]\label{lfgraphex}
Let $\xi$ be as in \cref{thm:graphex_rep_theorem} for fixed $\alpha$, and drop the first coordinate from 
the definitions of $\IsoF$, $\StarF$, and $W$.
Let $\mu_W(t) = \Lebesgue W(t,\cdot) = \int_{\NNReals} \!W(t,t')\, \intd t'$.
The random measure $\xi$ is  a.s.\ locally finite iff these four conditions are fulfilled:
\begin{enumerate}[label=\normalfont(\roman*)]
\item $\IsoF < \infty$, 
\item $\Lebesgue \StarF = \int_{\NNReals} \StarF(t)\, \intd t < \infty$,
\item $\Lebesgue\{ \mu_W = \infty \} = 0$ and $\Lebesgue \{\mu_W > 1\} < \infty$, 
\item $\Lebesgue^2[W; \mu_W \vee \mu_W \le 1] = \int_{\NNReals^2} W(x,y)\,1[ \mu_W(x) \le 1 ]\, 1[\mu_W(y) \le 1] \intd x \intd y < \infty$,
\item $\int_{\NNReals} W(x,x) \, \intd x < \infty$.
\end{enumerate}
In particular, $\xi$ is a.s.\ locally finite if $S$ and $W$ are integrable and $I < \infty$.
\end{thm}

\begin{rem}
An example showing that there are nonintegrable $W$ admitting a.s.\ locally finite exchangeable adjacency measures 
is the function $W(x,y) = 1[xy \le 1]$. Its marginal is $\mu_W(x) = \frac 1 x$, which obviously satisfies (iii).  Moreover, $W=0$ a.e.\ on the set $\{ (x,y) : \mu_W(x) \wedge \mu_W(y) \le 1 \} = \{(x,y) : x, y \ge 1\}$, satisfying (iv).
\end{rem}

These conditions leads us to the following definition:
\begin{defn}
  A \defnphrase{graphex} is a triple $(\IsoF,\StarF,W)$, where $\IsoF \ge 0$ is a non-negative real,
  $\StarF : \NNReals \to \NNReals$ is integrable,
  and $W : \NNReals^2 \to [0,1]$ is symmetric, and satisfies parts (iii)--(v) of \cref{lfgraphex}.
\end{defn}
\NA{In situations where there is no risk of confusion, we will abuse nomenclature and use the term graphex to refer to
  the $W$ component alone, with the understanding that the corresponding triple is $(0,0,W)$.}

The name graphex is chosen in analogy to graphon, the limit object
in the dense graph setting, and graphing, the limit objects in the bounded
degree graph setting~\cite{Lovasz:2013:A}.

The marginal $\mu_W$ of the graphex component $W$ arises in the characterization of a.s.\ finite undirected graph point processes.  This function will turn out to be an important quantity in a number of different contexts.
\begin{defn}
  The \defnphrase{graphex marginal} is
  $\mu_{W}(x)=\int_{\NNReals}W(x,y)\intd y$.
\end{defn}

\cref{{thm:graphex_rep_theorem}} gives us a precise picture of the structure of random graphs corresponding to jointly exchangeable simple point processes:
First, the potential vertices are the points of a collection of Poisson processes.  For the graph component corresponding to $W$, 
there is a Poisson process on $\boldsymbol{\theta}\times\boldsymbol{\vartheta} = \NNReals^2$, 
and each pair of vertices 
$(\theta_{i},\vartheta_{i}),(\theta_{j},\vartheta_{j})$ of the process
are connected independently with probability 
$W(\vartheta_{i},\vartheta_{j})$.
For each vertex $(\theta_{i},\vartheta_{i})$ in this component, there is a corresponding Poisson process on $\NNReals$ with rate $\StarF(\vartheta_{i})$.  Every point of this Poisson process connects to the vertex $(\theta_{i},\vartheta_{i})$ and no other point.
Finally, a Poisson process on $\NNReals^2$ with rate $\IsoF$ produces pairs $(x,y) \in \NNReals^2$ of vertices that are connected to each other but no other vertices.

We now define the class of Kallenberg exchangeable graphs:

\begin{defn}
  A \defnphrase{Kallenberg exchangeable graph (KEG) associated with
  graphex $(\IsoF,\StarF,W)$} is the random graph $G$ associated with an exchangeable adjacency measure
                 $\xi$ of the form given in \cref{KEGgen}.
  The Kallenberg exchangeable graph model is the family of $\nu$-truncations $G_\nu = \xi( \cdot \, \cap [0,\nu]^2)$, for $\nu \in \NNReals$.
  When the graphex is clear from
  context, we will simply refer to $G$ as the Kallenberg exchangeable graph.
\end{defn}

The first term of \cref{KEGgen} gives
essentially all of the interesting graph structure, and so for the rest of
the paper, we will restrict attention to models that take
$\StarF=\IsoF=0$.  Before doing so, we note that the natural analogue of \ErdosRenyi\ graphs in the KEG model corresponds to graphs for which $\IsoF \ge 0$, $\StarF= 0$, and $W$ is constant on a set of the form $[0,c]^2$ and 0 otherwise.  In this case, if $W$ is not identically zero, then later results will imply that the truncated graph sequence is dense.

Consider now the structure arising from $W$ alone.  Because $\IsoF=\StarF=0$, we will refer to $W$ as the graphex without any risk of confusion.
Let $\PP$ be a unit rate Poisson process on
$\boldsymbol{\theta}\times\boldsymbol{\vartheta}$ as in
\cref{thm:graphex_rep_theorem}.  A Kallenberg exchangeable graph $G$
associated with $W$ has vertex set
\[
v(G)=\{
\theta_{i}\suchthat(\theta_{i},\vartheta_{i})\in\PP\AND\exists\theta_{j}\in\PP\st
W(\vartheta_{i},\vartheta_{j})>\zeta_{\{ i,j\} }\}
\]
and edge set
\[
e(G)=\{ \{ \theta_{i},\theta_{j}\}
\suchthat(\theta_{i},\vartheta_{i}),(\theta_{j},\vartheta_{j})\in\PP\AND
W(\vartheta_{i},\vartheta_{j})>\zeta_{\{ i,j\} }\} .
\]
 
\begin{rem}
  A graphex with $W(\vartheta,\vartheta)=0\mbox{ for all
  }\vartheta\in\NNReals$ generates a KEG with no self edges.
\end{rem}

\begin{rem}
  Notice that if $G$ is a KEG associated to $W$ and $G_{\nu}$ is $G$
  restricted to $\left[0,\nu\right]$ then $G_{\nu}$ is \emph{not} the
  same as the induced subgraph of $G$ given by restricting to vertices of
  $G$ with labels $\le\nu$. The reason for this is that the induced
  subgraph includes an (infinite) collection of vertices that do not
  connect to any edges. However, it is true that $G_{\nu}\upto G$ in
  the sense that $v(G_{\nu})\upto v(G)$ and $e(G_{\nu})\upto e(G)$ as
  $\nu\upto\infty$.
\end{rem}

\begin{rem}
 The model can be extended to weighted graphs by replacing the
 indicator term $1[W(\alpha,\vartheta_{i},\vartheta_{j}) \le \zeta_{\{i,j\}}]$ by a general random variable parameterized by
 $W(\alpha,\vartheta_i,\vartheta_j)$.  The model can be extended to directed graphs
 by mimicking the 4-graphon approach used by
 \cite{Cai:Ackerman:Freer:2015} to extend the exchangeable graph model to
 directed graphs.
\end{rem}

\begin{defn}
  We will often refer to $\PP$ as the \defnphrase{latent Poisson
    process}.  For a point of the latent Poisson process
  $(\theta_{i},\vartheta_{i})\in\PP$ the \defnphrase{label} of the
  point is $\theta_{i}$ and the \defnphrase{latent value} is
  $\vartheta_{i}$.
\end{defn}

We close this section with a word of warning about point process
notation:

\begin{rem}
  Point processes are central to our construction. For a point process
  $\boldsymbol{P}$ we will often refer to points
  $p_{i}\in\boldsymbol{P}$ where the index $i$ is given by some
  unspecified measurable function of $\boldsymbol{P}$.  For example, if
  $\boldsymbol{P}$ is a Poisson process then the points could be
  indexed by the ordering of their Euclidean distances to the
  origin. This is convenient for writing summations across the point
  process and for unambiguously associating dimensions when the points
  are multidimensional (e.g., $p_{i}=(a_{i},b_{i})$ then we understand
  $a_{i}$ and $b_{i}$ are part of the same tuple in $\boldsymbol{P}$).
  However, there is a small subtlety here: any choice of indexing
  function will be informative about the value of the point of the
  process. For example, if the points of a Poisson process are indexed
  by their distance to the origin then the value of the index is
  informative about the value of the point. 
  As a result, some care must be taken when making statements of (conditional) independence.
\end{rem}

\section{Expected Number of Edges and Vertices\label{sec:Expectations}}
\global\long\def\indexfn#1{\scriptscriptstyle #1}

In this section we derive the expected values of the number of
vertices and edges of Kallenberg exchangeable graphs restricted to
$[0,\nu]$, in terms of their underlying graphex. We focus on those 
graphexes where $\IsoF=0$ and $\StarF=0$ so we refer
to $W$ as the graphex without any risk of confusion.
Throughout this
section we implicitly assume $W$ is non-random; in the case of random
$W$ the results can be understood as conditional statements.

The intuition for the main proof idea is to find the distribution of
the degree of a single point in the latent Poisson process, write the
statistics of interest as sums of functions of the degrees of the
points and appeal to the linearity of expectation to evaluate these
expressions. For example, the number of edges in a graph is the sum of
the degrees of all of the vertices divided by 2. This perspective allows
the use of powerful techniques for computing expectations of sums over
point processes.

Because the $\theta$ labels of the graph carry no information it is
easiest to treat $G_{\nu}$ by projecting the latent Poisson process
$\PPnu$ along its second coordinate on to a random point set in
$\boldsymbol{\vartheta}\simeq\NNReals$ as $\PPnu^{P}=\{
\vartheta_{i}\suchthat(\theta_{i},\vartheta_{i})\in\PPnu\} $, which is
then a rate $\nu$ Poisson process.  For $\varphi$ a locally finite,
simple sequence and $\{z_{\{i,j\}}\}$ a sequence of values in $[0,1]$
such that $z_{ij}=z_{ji}$, then for $x\in\varphi$ define the degree
function:
\begin{equation}
  D(x,\varphi,\{z_{ij}\})=\sum_{p\in\varphi\backslash\{x\}}1[W(x,p)\ge z_{\indexfn{i}(x)\indexfn{i}(p)}]
  +2\cdot1[W(x,x)\ge z_{\indexfn{i}(x)\indexfn{i}(x)}]\label{eq:degree_function}
\end{equation}
where $\indexfn{i}(x)=\indexfn{i}(x,\varphi)$ gives the index of the
point $x\in\varphi$ with respect to the natural ordering on
$\NNReals$. Intuitively speaking, for a symmetric array
$\zeta_{\{i,j\}}$ of uniform $[0,1]$ random variables,
\[
D(\vartheta,\PPnu^{P},(\zeta_{\{i,j\}}))
\]
is the degree of a point $(\theta,\vartheta)\in\PPnu$ under a KEG
process, conditional on $(\theta,\vartheta)\in\PPnu$.

For any $\lambda\in\NNReals$ the probability that
$\lambda\in\PPnu^{P}$ is $0$ and so
$D(\lambda,\PPnu^{P},\zeta_{\{i,j\}})$ is ill defined. We wish to
derive the distribution of the degree of a point $\lambda$ under the
promise that it's in the point process. Because this is a measure 0
event the conditioning is in general somewhat tricky. The idea is
formalized by Palm theory, which for a measure $P$ on point sequences
defines a Palm measure $P_{\lambda}$ that behaves as the required
conditional distribution; see \cite{Chiu:Stoyan:Kendall:Mecke:2013}
for an accessible introduction.  The Slivnyak--Mecke theorem asserts
that a Poisson process $\PP$ with a promise $\lambda\in\PP$ (in the
Palm sense) is equal in distribution to $\PP\cup\{\lambda\}$, so the
correct object to work with is
$D(\lambda,\PPnu^{P}\cup\{\lambda\},\zeta_{\{i,j\}})$. Recalling the
graphex marginal $\mu_{W}(x)=\int_{\NNReals}W(x,y)\intd y$:
\begin{lem}
  \label{lem:marginal_degree_dist}Let $x\in\NNReals$. Then
  $D(\lambda,\PPnu^{P}\cup\{\lambda\},(\zeta_{\{i,j\}}))\equaldist
  D_{\mbox{ext}}+D_{\mbox{self}}$ where
  $D_{\mbox{ext}}\dist\poiDist(\nu\mu_{W}(\lambda))$ and $\frac 1 2
  D_{\mbox{self}} \dist \bern (W(\lambda,\lambda))$
  independently.\end{lem}
\begin{proof}
  With probability 1, $\lambda\notin\PPnu^{P}$ so
  \begin{align}
    D(\lambda,\PPnu^{P}\cup\{\lambda\},\zeta_{\{i,j\}}) & =
    \sum_{p\in\PPnu^{P}}1[W(\lambda,p)\ge\zeta_{\indexfn{i}(\lambda)\indexfn{i}(p)}]+2\cdot1[W(\lambda,\lambda)\ge\zeta_{\indexfn{i}(\lambda)\indexfn{i}(\lambda)}].
  \end{align}
  Since $\zeta_{{\indexfn i}(\lambda){\indexfn i}(\lambda)}\sim
  U[0,1]$ independent of everything else letting
  \[D_{\mbox{self}}=2\cdot 1[W(\lambda,\lambda)\ge\zeta_{{\indexfn
      i}(\lambda){\indexfn i}(\lambda)}]\] and
  \[D_{\mbox{ext}}=\sum_{p\in\PPnu^{P}}1[W(\lambda,p)\ge\zeta_{{\indexfn
      i}(\lambda){\indexfn i}(p)}]\] establishes the independence of
  the two terms and that $\frac 1 2 D_{\mbox{self}}\dist \bern
  (W(\lambda,\lambda))$.

We have that
  \[
  \int_{\NNReals}\int_{\left[0,1\right]}1\left[u\le
    W(\lambda,y)\right]\nu\intd y \intd u
  =\nu\int_{\NNReals}W(\lambda,y)\intd y < \infty\mbox{ a.s.},
  \]
where the a.s. finiteness is one of the defining conditions of the graphex $W$.  
It then follows by a version of Campbell's theorem
  \citep[][\S5.3]{Kingman:1993}, the characteristic function of
  $D_{\mbox{ext}}$ is
  \begin{align}
    \expect{\exp(itD_{\mbox{ext}})} & = \expect{\exp(it\sum_{p\in\PPnu^{P}}1\left[\zeta_{{\indexfn i}(\lambda){\indexfn i}(p)}\le W(\lambda,p)\right])}\\
    & = \exp\{ \int_{\NNReals}\int_{[0,1]}(1-e^{it1\left[u\le W(\lambda,y)\right]}\nu\intd u\intd y)\} \\
    & = \exp\{ \nu\sum_{n=1}^{\infty}\frac{(it)^{n}}{n!}\int_{\NNReals}\int_{[0,1]}1\left[u\le W(\lambda,y)\right]\intd u\intd y\} \\
    & = \exp\{ \nu\mu_{W}(\lambda)(e^{it}-1)\} .
  \end{align}
  Hence, $D_{\mbox{ext}}$ is a $\poiDist(\nu\mu_{W}(\lambda))$
  distributed random variable, completing the proof.
\end{proof}
We would now like to access the first moments of various graph
quantities by writing them as sums of (functions of) the degree and
exploiting the linearity of expectation to circumvent
dependencies. For example, the total number of edges of the graph is
\[
e_{\nu}\equaldist\frac{1}{2}\sum_{\vartheta\in\PPnu^{P}}D(\vartheta,\PPnu^{P},(\zeta_{\{i,j\}})),
\]
where the equality is in distribution (as opposed to almost sure)
because the indexing $\indexfn{i}(x)$ of the latent Poisson process
used by the degree function is not the same as the indexing used in
\cref{thm:graphex_rep_theorem}.

Standard point process formulas deal with computing expressions of the
form
\[
\expect{\sum_{\lambda\in\Gamma}h(\lambda,\Gamma)}
\]
where $\Gamma$ is a simple point process. Sums across the degrees of
points of the process do not immediately have this form because the
degree depends on the \iid uniform array $(\zeta_{\{i,j\}})$, so we
will need a slight extension. Let $\pointspace$ denote the family of
all sets of points $\varphi$ in $\NNReals$ that are both locally
finite and simple, then:
\begin{lem}[{Extended Slivnyak--Mecke}]
  \label{lem:Extended-Slivnyak-Mecke} Let $\Phi$ be a rate $\nu$
  Poisson process on $\NNReals$, $U$ an independent uniform random
  variable, and $f:\NNReals\times\pointspace\times[0,1]\to\NNReals$ a
  measurable non-negative function. Then
  \[
  \expect{\sum_{p\in\Phi}f(p,\Phi,U)}=\nu\int_{\NNReals}\expect{f(x,\Phi\cup\{x\},U)}\intd
  x.
  \]
\end{lem}
\begin{proof}
  By the independence of $U$ and $\Phi$, the non-negativity of $f$,
  and Tonelli's theorem, we have
  \[
  \expect{\sum_{p\in\Phi}f(p,\Phi,U)}=\int_{0}^{1}
  \expect{\sum_{p\in\Phi}f(p,\Phi,u)} \intd u.
  \]
  By the usual Palm calculus, the inner expectation satisfies
  \[
  \expect{\sum_{p\in\Phi}f(p,\Phi,u)} = \int_{\NNReals}
  \int_{\pointspace}f(x,\varphi,u)P_{x}(\intd\varphi)\nu\intd x,
  \]
  where $P_{x}$ is the local Palm distribution of a unit rate Poisson
  process. Letting $P$ be the distribution of a unit rate Poisson
  process, the Slivnyak--Mecke theorem gives:
  \[
  \int_{\pointspace}f(x,\varphi,u)P_{x}(\intd\varphi) =
  \int_{\pointspace}f(x,\varphi\cup\{x\},u)P(\intd\varphi).
  \]
  The result then follows by a second application of Tonelli's theorem
  to change the order of integration.
\end{proof}
The main results of this section now follow easily:
\begin{thm}
  \label{thm:expected_edges}The expected number of edges
  $e_{\nu}=|\edgeset{G_{\nu}}|$ is
  \[
  \expect{e_{\nu}}=\frac{1}{2}\nu^{2}\iint_{\NNReals^{2}}W(x,y)\intd
  x\intd y+\nu\int_{\NNReals}W(x,x)\intd x.
  \]
\end{thm}
\begin{proof}
  By \cref{lem:Extended-Slivnyak-Mecke,lem:marginal_degree_dist},
  \begin{align}
    \expect{e_{\nu}} & = \frac{1}{2}\expect{\sum_{\vartheta\in\PPnu^{P}}D(\vartheta,\PPnu^{P},(\zeta_{\{i,j\}}))}\\
    & = \frac{1}{2}\nu\int_{\NNReals}\expect{D(x,\PPnu^{P}\cup\{x\},(\zeta_{\{i,j\}}))}\intd x\\
    & = \frac{1}{2}\nu\int_{\NNReals}\nu\mu_{W}(x)+2W(x,x)\intd x
  \end{align}
  By assumption, $\norm{\mu_{W}}_1 = \norm {W}_1 <\infty$ and
  $\int_{\NNReals}W(\lambda,\lambda)\intd\lambda<\infty$, and so
  $\expect{e_{\nu}} < \infty$ and the result follows by the linearity
  of integration.
\end{proof}

\begin{thm}
  \label{thm:expected_vertices}The expected number of visible vertices
  $v_{\nu}=|\vertexset{G_{\nu}}|$ is
  \[
  \expect{v_{\nu}}=\nu\int_{\NNReals}(1-e^{-\nu\mu_{W}(x)})\intd x
  +\nu\int_{\NNReals}e^{-\nu\mu_{W}(x)}W(x,x)\intd x.
  \]
\end{thm}
\begin{proof}
  By \cref{lem:Extended-Slivnyak-Mecke,lem:marginal_degree_dist},
  \begin{align}
    \expect{v_{\nu}} & = \expect{\sum_{\vartheta\in\PPnu^{P}}1\left[D(\vartheta,\PPnu^{P},(\zeta_{\{i,j\}}))\ge1\right]}\\
    & = \nu\int_{\NNReals}\Pr(D(x,\PPnu^{P}\cup\{x\},(\zeta_{\{i,j\}}))\ge1)\intd x\\
    & = \nu\int_{\NNReals}1-\Pr(D_{\mbox{ext}}=0)\Pr(D_{\mbox{self}}=0)\intd x\\
    & = \nu\int_{\NNReals}1-e^{-\nu\mu_{W}(x)}(1-W(x,x))\intd x,
  \end{align}
  where $D_{\mbox{ext}}$ and $D_{\mbox{self}}$ are defined as in
  \cref{lem:marginal_degree_dist}. Splitting up the integral is
  justified since 
  $1-\exp(-\nu\mu_W(x)) \ge 0$ and
  $\exp(-\nu\mu_W(x))W(x,x) \ge 0$ for all $x$.
\end{proof}

A nearly identical argument can be used to find the expected number of
vertices of a specified degree. This result is interesting in its own
right and is used as a lemma in \cref{sec:Degree-distribution}.
\begin{thm}\label{thm:expected_deg_k}
  The expected number of vertices of degree $k$
  in $G_{\nu}$, $\numDegNu k$, is
  \[
  \begin{split}
    \expect{\numDegNu k}
    &=\nu^{k+1}\int_{\NNReals}\biggl[\frac{\mu_{W}(x)^{k}}{k!}e^{-\nu\mu_{W}(x)} \\
    &\qquad\qquad\quad
    +\frac{1}{\nu^{2}}\frac{\mu_{W}(x)^{k-2}}{(k-2)!}e^{-\nu\mu_{W}(x)}(1-\frac{(\nu\mu_{W}(x))^{2}}{k(k-1)})W(x,x)
    \biggr] \intd x
  \end{split}
  \]
\end{thm}
\begin{proof}
  The result follows from essentially the same argument as the
  previous two theorems and some straightforward algebraic
  manipulations.
\end{proof}
Notice that in the limit as $\nu\to\infty$ the contribution of self
edges ($W(\lambda,\lambda)\neq0$) is negligible in the sense that
terms due to the edges between distinct vertices dominate asymptotically
for \cref{thm:expected_edges,thm:expected_vertices,thm:expected_deg_k}.

We end this section by applying our results on the expected number of
vertices and edges to show that a KEG is dense \iff the generating
graphex is compactly supported.

\begin{thm}\label{thm:dense_iff_compact}
  Let $G$ be Kallenberg exchangeable graph with graphex $(0,0,W)$.
  If $W$ is compactly supported, then $G$ is dense with probability 1.
  Conversely, if $W$ is integrable and not compactly supported, then $G$ is sparse
  with probability 1.
\end{thm}
\begin{proof}
  We have already shown in \cref{sec:graphon_models} that if $W$ is
  compactly supported then the corresponding KEG is dense (or empty)
  with probability 1 because these models correspond exactly to
  graphon models.

  Conversely, suppose that the KEG $G$ generated by $W$ is dense with
  positive probability. This means that there are constants $c,p>0$
  such that
  \[
  \liminf_{\nu\to\infty}\Pr(e_\nu>cv_\nu^2)>p,
  \]
  where $e_\nu=\edgeset{G_\nu}$ and $v_\nu=\vertexset{G_\nu}$. With
  \[
  \expect{e_\nu} \ge \Pr(e_\nu>cv_\nu^2)\expect{cv_\nu^2}
  \]
  and Jensen's inequality, this implies
  $\expect{e_\nu}=\Omega(\expect{v_\nu}^2).$

  Now, by \cref{thm:expected_vertices},
  \[
  \expect{v_{\nu}}=\nu\int_{\NNReals} 1-e^{-\nu\mu_{W}(x)} \intd x
  +\nu\int_{\NNReals}e^{-\nu\mu_{W}(x)}W(x,x)\intd x,
  \]
  and monotone convergence shows $\int_{\NNReals} 1-e^{-\nu\mu_{W}(x)}
  \intd x \upto \infty$ \iff $\mu_W$ is not compactly supported. Thus
  for $G$ dense with positive probability and $W$ not compactly
  supported it holds that
  \[
  \expect{e_\nu}=\omega(\nu^2).
  \]
  However, by \cref{thm:expected_edges},
  $\expect{e_\nu}=\Theta(\nu^2).$ This contradiction completes the
  proof.
  
\end{proof}

\section{Degree Distribution in the Asymptotic Limit\label{sec:Degree-distribution}}
\global\long\def\PPnew{\Pi_{(1,\nu+1]}}

\global\long\def\PPrest{\Pi_{(1,\infty)}}

\global\long\def\PPnaught{\Pi_{0}}

\newcommandx\degNu[1][usedefault, addprefix=\global, 1=]{D_{#1,\nu}^
  {}}

\global\long\def\pseudoDeg{H}

\global\long\def\Ngreater#1{N_{>#1}^{{}_{(\nu)}}}

\global\long\def\ngreater#1{n_{>#1}^{{}_{(\nu)}}}

\global\long\def\tNgreater#1{\tilde{N}_{>#1}^{{}_{(\nu)}}}

One of the major advantage of KEGs over previous exchangeable graph models is
that they allow for sparse graphs of the kind typically seen in
application\NA{; in particular this means the KEG models should allow for a
variety of degree (scaling) behaviours.
Caron and Fox \cite{Caron_Fox_CRM_Graphs} characterized the degree
distribution in the large graph limit for the particular case of directed graphs
based on generalized gamma processes. 
We now describe the limiting degree
distribution of Kallenberg exchangeable graphs.} 
We focus on those graphexes where $\IsoF=\StarF=0$ so we refer
to $W$ as the graphex without any risk of confusion.
To formalize the notion of limiting degree distribution, let $G_{\nu}$
be a Kallenberg exchangeable graph on $[0,\nu)$ with graphex $W$, and
let $D_{\nu}$ be the degree of a vertex chosen uniformly at random
from $\vertexset{G_{\nu}}$. The central object of study is then the
random distribution function $k\mapsto P(D_{\nu}\le k\given G_{\nu})$
and its scaling limit. The primary aim of this section is to prove the
following theorem:
\begin{thm}
  \label{thm:deg_dist_theorem}Let $W$ be an integrable graphex such that
  \begin{enumerate}
  \item There exist some constants $C,T>0$ such that for all $\lambda$
    and $\omega>T$ it holds that $\int W(\lambda,x)W(\omega,x)\intd x
    \le C\mu_{W}(\lambda)\mu_{W}(\omega)$.
  \item $\mu_{W}$ is monotonically decreasing.
  \item $\mu_{W}$ is differentiable.
  \item There is some $\chi>0$ such that for all $x>\chi$ holds that
    $\frac{\mu_{W}(x)}{\mu_{W}'(x)}\frac{1}{x}\ge-1$.
  \end{enumerate}
  Let $k_\nu=o(\nu)$.
  Then,
  \[
  \Pr(D_{\nu}>k_{\nu}\given
  G_{\nu})\convPr\lim_{\nu\to\infty}\frac{\sum_{n=k_{\nu+1}}^{\infty}\int\frac{1}{n!}e^{-\nu\mu_{W}(x)}(\nu\mu_{W}(x))^{n}\intd
    x}{\int1-e^{-\nu\mu_{W}(x)}\intd x}.
  \]

\end{thm}

In the case $\mu_{W}(\lambda)=(1+\lambda)^{-2}$ the right hand side of
this expression is in $(0,1)$ for $k_{\nu}=k$ for any choice of
$k$. That is, even in the infinite graph limit a constant fraction of
the vertices will have degree $\le k$ for a fixed integer $k$. By
contrast, for $\mu_{W}(\lambda)=e^{-\lambda}$ the degree of a randomly
chosen vertex goes to $\infty$ so, for fixed $k$, $\Pr(D_{\nu}>k\given
G_{\nu})\convPr1$.  However, we saw that $\Pr(D_{\nu} >
\nu^\beta\given G_{\nu})\to 1-\beta$ for $\beta\in(0,1)$; i.e., taking
$k_{\nu}=\nu^\beta$ results in a non-trivial limit on the right hand
side. That is, this theorem can be understood intuitively as
characterizing the rate of growth of the degree of a typical
vertex. This scaling limit affords a precise notion of ``how dense'' the
graph associated to a particular graphex is.

Let $\ngreater{l}$ denote the number of vertices of $G_{\nu}$ with degree
greater than $l$. It is immediate that
\[
\Pr(D_{\nu}\ge k_{\nu}\given
G_{\nu})=\frac{\ngreater{k_{\nu}}}{\ngreater{0}},
\]
i.e., the probability of choosing a vertex of degree greater than
$k_{\nu}$ is the proportion of such vertices among all vertices. Notice
that, even for fixed $l$, the random variable $\ngreater{l}$ grows with
$\nu$.
Further notice that like $D_{\nu}$ the random variable
$\ngreater{k_{\nu}}/\ngreater{0}$ is ill defined for the event
$\ngreater{0}=0$; however this is a measure $0$ event in the limit
$\nu\to\infty$. The content of \cref{thm:deg_dist_theorem} can be
understood as saying that the limit of the ratio
$\frac{\ngreater{l}}{\ngreater{0}}$ is the limit of the ratio of the
expectations,
\[
\frac{\ngreater{l}}{\ngreater{0}}\convPr\lim_{\nu\to\infty}\frac{\expect{\ngreater{l}}}{\expect{\ngreater{0}}}\asympLim{\nu}.
\]

Reasoning about the degree of a randomly selected vertex is
substantially simplified by selecting only from those with label
$\theta\in[0,1]$ and ignoring the contribution of edges
$(\theta_{i},\theta_{j})$ with $\theta_{i},\theta_{j}\le1$. The reason
for this is that it allows us to eliminate one form of dependence
between the degrees of distinct points; namely the dependence arising
from the requirement that each terminus attached to a vertex has a
matching terminus attached to some other vertex in the set. Intuitively,
studying this simplification is valid because the $\theta$ labels of
the points of the latent Poisson process are independent of their
degrees and as the graph becomes large only a negligible number of
edges have both termini with labels $\theta\le1$. Let $\Ngreater{l}$
be the number of vertices of $G_{\nu}$ with label $\theta_{i}<1$ and
greater than $l$ neighbours $\{\theta_{j}\}$ where $\theta_{j}>1$. The
following lemma establishes the claimed equivalence:
\begin{lem}
  \label{lem:restric_approx_valid}The limiting distribution of
  $\ngreater{l}/\ngreater{0}$ is the same as the limiting distribution
  of the ratio that considers only vertices with label $\theta_{i}\le1$
  and counts only edges $(\theta_{i},\theta_{j})$ with $\theta_{j}>1$,

\[
\lim_{\nu\to\infty}\frac{\ngreater{l}}{\ngreater{0}}\equaldist\lim_{\nu\to\infty}\frac{\Ngreater{l}}{\Ngreater{0}}.
\]
\end{lem}
\begin{proof}
  The validity of this equality is a consequence of the following
  three observations:
  \begin{enumerate}
  \item $\lim_{\nu\to\infty}\Pr(\Ngreater{0}=0)=0$ so
    $\lim_{\nu\to\infty}\frac{\Ngreater{l}}{\Ngreater{0}}$ is well
    defined.
  \item The $\theta$ label of a point of the latent Poisson process is
    independent of its degree. Let $\tilde{D}_{\nu}$ be the degree of
    a vertex chosen uniformly at random from those members of
    $\vertexset{G_{\nu}}$ with label $\theta<1$ and let
    $\tNgreater{l}$ be the number of such vertices with degree greater
    than $l$. Because the degree of a point
    $(\theta_{i},\vartheta_{i})\in\PP$ is independent of the value of
    $\theta_{i}$ it holds that, conditional on $\{\tNgreater{0}>0\}$,
    \begin{align}
      \Pr(D_{\nu}>l\suchthat G_{\nu}) &
      \equaldist\Pr(\tilde{D}_{\nu}>l\suchthat G_{\nu}).
    \end{align}
    This immediately implies
    \begin{align}
      \frac{\ngreater{l}}{\ngreater{0}} &
      \equaldist\frac{\tNgreater{l}}{\tNgreater{0}}.
    \end{align}

  \item The number of edges $(\theta_{i},\theta_{j})$ with
    $\theta_{i},\theta_{j}\le1$ is almost surely finite and
    $\Ngreater{0}\upto\infty$ almost surely, so the probability of
    randomly choosing a vertex that participates in at least one of
    the neglected edges goes to $0$ as $\nu\to\infty$, thus
    \[
    \lim_{\nu\to\infty}\Pr(\tilde{D}_{\nu}>l\suchthat
    G_{\nu})\equalas\lim_{\nu\to\infty}\frac{\Ngreater{l}}{\Ngreater{0}}.
    \]

  \end{enumerate}
\end{proof}
To treat the limiting distribution of this ratio we introduce
\begin{align}
  \PPnaught & = \{\vartheta\suchthat(\theta,\vartheta)\in\mbox{\ensuremath{\PP}}_{\nu+1},\theta\le1\}\\
  \PPnew & =
  \{(\theta,\vartheta)\suchthat(\theta,\vartheta)\in\mbox{\ensuremath{\PP}}_{\nu+1},\theta>1\},
\end{align}
i.e., we break the latent Poisson process into the component with
$\theta\le1$ and the component with $\theta>1$ and then project out
the $\theta$ value of $\PPnaught$ since it contains no useful
information. Notice that $\PPnaught$ and $\PPnew$ are independent
Poisson processes.

For $x\in\NNReals$, $\bar{u}=(u_{i})$ a sequence of values in $[0,1]$
and $\{(\phi_{i},\varphi_{i})\}$ a locally finite, simple sequence
with elements in $(1,\infty)\times\NNReals$ we define
\[
\degNuFn(x,\bar{u},\{(\phi_{i},\varphi_{i})\})=\sum_{i}1[W(x,\varphi_{i})>u_{i}]1[\phi_{i}\le\nu+1].
\]
There exists a marking $(\lambda_{i},\bar{\zeta}_{i})$ of $\PPnaught$
where each $\bar{\zeta}_{i}=(\zeta_{j}^{i})$ is a sequence of
independent $U[0,1]$ random variables such that
\[
\degNuFn(\lambda,\bar{\zeta}_{i},\PPrest)
\]
is the degree of the point $\lambda\in\PPnaught$. Let
$\bar{U}_{j}=(U_{i}^{j})$ be independent sequences of independent
$U[0,1]$ random variables and define
\[
\degNu[j](x)=\degNuFn(x,\bar{U}_{j},\PPrest).
\]
These random variables will arise naturally in the course of the
proof.

It follows by mimicking the proof of \cref{lem:marginal_degree_dist}
that
\[
\degNu[j](x)\dist\poiDist(\nu\mu_{W}(x))
\]
marginally. The importance of
$\degNuFn(\lambda,\bar{\zeta}_{i},\PPrest)$ in the context of the
present section comes from the relation
\[
\Ngreater{l}=\sum_{i}1[\degNuFn(\lambda_{i},\bar{\zeta}_{i},\PPrest)>l].
\]
where $(U_{i})^{\lambda}$ is a marking of $\PPnaught$. We will make
heavy use of the observation that, by Campbell's formula,
\[
\expect{\Ngreater{l}}=\int\Pr(\degNu[1](x)>l)\intd x.
\]

The idea of the proof of \cref{thm:deg_dist_theorem} is to show that
\begin{equation}
  \Ngreater{k_{\nu}}/\expect{\Ngreater{0}}\convPr\lim_{\nu\to\infty}\expect{\Ngreater{k_{\nu}}}/\expect{\Ngreater{0}}\asympLim{\nu}.\label{eq:deg_dist_basic_proof_idea}
\end{equation}
The special case $k_{\nu}=0$ gives
$\Ngreater{0}/\expect{\Ngreater{0}}\convPr1$ and an application
Slutsky's theorem then establishes
\[
\frac{\Ngreater{k_{\nu}}}{\Ngreater{0}}\convPr\lim_{\nu\to\infty}\expect{\Ngreater{k_{\nu}}}/\expect{\Ngreater{0}}\asympLim{\nu}.
\]
Using Chebyshev's inequality, a sufficient condition for
\cref{eq:deg_dist_basic_proof_idea} to hold is
\[
\var{\Ngreater{k_{\nu}}}=o(\expect{\Ngreater{0}}^{2}).
\]
The majority of the proof is aimed at characterizing the growth rate
of $\var{\Ngreater{k_{\nu}}}$.

In order to do this, we will need to make an assumption about the
graphex $W$ that controls the average dependence between the degrees
of different vertices of $G_{\nu}$:
\begin{assumption}
  \label{asmp:limited_dependence_assumption}There exist some constants
  $C,T>0$ such that for all $\lambda$ and $\omega>T$ it holds that
  $\int W(\lambda,x)W(\omega,x)\intd x \le
  C\mu_{W}(\lambda)\mu_{W}(\omega)$.

\end{assumption}
We do not know of any examples of an integrable graphex that violates this
assumption, although $W(x,y)=1[xy<1]$ does.  To understand what the assumption means, let
$L(\lambda,\omega)$ be the number of common neighbours of points
$(l,\lambda),(w,\omega)\in\PPnu$ under $G_{\nu}$ and observe that for
a graphex $W$ that is $0$ on the diagonal (i.e., forbidding
self-edges),
\[
L(\lambda,\omega)\dist\poiDist(\nu\int W(\lambda,x)W(\omega,x)\intd
x),
\]
with respect to the Palm measure $P_{\lambda,\omega}$%
\footnote{Recall this is just the measure that guarantees that
  $\lambda,\omega$ are elements of the point process.%
}. This can be shown by an argument very similar to
\cref{lem:marginal_degree_dist}.  Thus the assumption can be
understood as requiring that the average number of common neighbours
between a pair of vertices is at most a constant factor larger than it
would be in the case $W(x,y)=\mu_{W}(x)\mu_{W}(y)$.

We further assume for simplicity that $\mu_{W}(x)$ is strictly
monotonically decreasing, differentiable and that there is some
$\chi>0$ such that for all $x>\chi$ holds that
$\frac{\mu_{W}(x)}{\mu_{W}'(x)}\frac{1}{x}\ge-1$.  It is not clear
which, if any, of these assumptions are necessary for the result to
hold. The last condition in particular may already be implied by the
other assumptions. Moreover, the result will hold automatically for a
graphex $W$ if there is some other graphex $W^{'}$ such that $W^{'}$
satisfies the conditions of the theorem and the KEGs corresponding to
$W$ and $W^{'}$ are equal in distribution.

Invertibility implies that $W$ does not have compact support; i.e.,
the graph is sparse (\cref{thm:dense_iff_compact}). A particular
consequence of this last assumption is that for any function
$l(\nu)\to0$ as $\nu\to\infty$ it holds that
$\mu_{W}^{-1}(l(\nu))\to\infty$, a fact that will be used heavily in
this section and the next.

Subject to these assumptions we may now begin the argument to bound
$\var{\Ngreater{k_{\nu}}}$.
\begin{lem}
  \label{lem:var_expression}Let $k_{\nu}=o(\nu)$, then
  \begin{align}
    \var{\Ngreater{k_{\nu}}} & = \expect{\Ngreater{k_{\nu}}} \\
    &+\iint\Pr(\degNu[1](x)>k_{\nu},\degNu[2](y)>k_{\nu})-\Pr(\degNu[1](x)>k_{\nu})\Pr(\degNu[2](y)>k_{\nu})\intd
    x\intd y
  \end{align}
\end{lem}
\begin{proof}
  Let $\{(\lambda_{i},\bar{\zeta}_{i})\}$ be a marking of $\PPnaught$
  such that each $\bar{\zeta}_{i}=(\zeta_{j}^{i})$ is a sequence of
  independent identically distributed $U[0,1]$ random variables and
  \[
  \degNuFn(\lambda_{i},\bar{\zeta}_{i},\PPrest)
  \]
  is the degree of point $\lambda$. Conditional on $\PPrest$ the
  degrees $\degNuFn(\lambda,\bar{\zeta}_{i},\PPrest)$ of each point
  $\lambda\in\PPnaught$ are a marking of $\PPnaught$ so
  \[
  \Ngreater{k_{\nu}}\given\PPrest\dist\poiDist(\expect{\Ngreater{k_{\nu}}\suchthat\PPrest}).
  \]
  Using this, the formula for conditional variance is
  \begin{align}
    \var{\Ngreater{k_{\nu}}} & = \expect{\var{\Ngreater{k_{\nu}}\suchthat\PPrest}}+\var{\expect{\Ngreater{k_{\nu}}\suchthat\PPrest}}\\
    & =
    \expect{\Ngreater{k_{\nu}}}+\var{\expect{\Ngreater{k_{\nu}}\suchthat\PPrest}}.
  \end{align}

  An application of Campbell's formula to the second term gives:
  \begin{align}
    \expect{\Ngreater{k_{\nu}}\given\PPrest} & = \int_{\NNReals}\expect{1[\degNuFn(x,\bar{U},\PPrest)>k_{\nu}]\given\PPrest}\intd x\\
    & = \int_{\NNReals}\Pr(\degNu[1](x)>k_{\nu}\given\PPrest)\intd x,
  \end{align}
  where $\bar{U}$ is a sequence of $U[0,1]$ random variables
  independent of $\PPrest$. Then
  $\expect{\Ngreater{k_{\nu}}\given\PPrest}^{2}$ is
  \[
  \iint_{\NNReals^{2}}\Pr(\degNu[1](x)>k_{\nu}\wedge\degNu[2](y)>k_{\nu}\given\PPrest)\intd
  x\intd y.
  \]
  By Tonelli's theorem,
  \begin{align}
    \expect{\expect{\Ngreater{k_{\nu}}\given\PPrest}^{2}} & =
    \iint_{\NNReals^{2}}\Pr(\degNu[1](x)>k_{\nu}\wedge\degNu[2](y)>k_{\nu})\intd
    x\intd y
  \end{align}
  whence
  \begin{align}
    \var{\expect{\Ngreater{k_{\nu}}\given\PPrest}} & = \iint_{\NNReals^{2}}\Pr(\degNu[1](x)>k_{\nu}\wedge\degNu[2](y)>k_{\nu})\intd x\intd y\\
    &
    -\iint_{\NNReals^{2}}\Pr(\degNu[1](x)>k_{\nu})\Pr(\degNu[2](y)>k_{\nu})\intd
    x\intd y
  \end{align}
  and the claimed result follows.
\end{proof}
Bounding the variance requires controlling the average dependence
between $\degNu[1](x)$ and $\degNu[2](y)$, as captured by the second
term in the lemma above. The degree of a point $\lambda$ gives
information about the degree of a point $\omega$ only through
$\PPnew$. Intuitively, as $\nu\to\infty$, the degree of $\lambda$
gives very little information about $\PPnew$ so the pairwise
dependence between degrees is weak and the variance of $\Ngreater{l}$
is small. Formalizing this intuition proves to be somewhat
tricky. Essentially, the strategy is to find a bound of the form
\[
\Pr(\degNu[1](x)>k_{\nu},\degNu[2](y)>k_{\nu})-\Pr(\degNu[1](x)>k_{\nu})\Pr(\degNu[2](y)>k_{\nu})\\
\le\Pr(\degNu[1](x)>k_{\nu})g(y)
\]
so that
\begin{align}
  \var{\Ngreater{k_{\nu}}} & \le\expect{\Ngreater{k_{\nu}}}+\iint\Pr(\degNu[1](x)>k_{\nu})g(y)\intd x\intd y\\
  & =\expect{\Ngreater{k_{\nu}}}(1+\int g(y)\intd y).
\end{align}
The goal is then to find a bounding function $g(y)$ such that $\int
g(y)\intd y$ is small. The next lemma provides such an expression.

\global\long\def\degDepLem#1{D_{#1}}

\begin{lem}
  \label{lem:dependence_bound} Let $T$ be a value such that for $y>T$
  it holds that
  \begin{align}
    \int W(x,z)W(y,z)\intd z & \le C\mu_{W}(x)\mu_{W}(y)
  \end{align}
  and
  \begin{align}
    2C\mu(y) & \le1-\log2.
  \end{align}
  Further, let $B(y)\dist\binDist(5k_{\nu},C\mu_{W}(y))$ independently
  of $\degNu[2](y)$ and define
  \begin{align}
    g(y) & =\begin{cases}
      \Pr(\degNu[2](y)\le k_{\nu}) & y\le T\\
      \Pr(\degNu[2](y)+B(y)>k_{\nu}\wedge\degNu[2](y)\le k_{\nu}) &
      y>T.
    \end{cases}
  \end{align}
  Then,
  \begin{align}
    \Pr(\degNu[1](x)>k_{\nu},\degNu[2](y)>k_{\nu})-\Pr(\degNu[1](x)>k_{\nu})\Pr(\degNu[2](y)>k_{\nu})\\
    \le\Pr(\degNu[1](x)>k_{\nu})g(y)
  \end{align}
\end{lem}
\begin{proof}
  Let $x,y\in\NNReals$ and define
  \begin{align}
    \degDepLem a & = \degNu[1](x)\\
    \degDepLem b & = \degNu[2](y).
  \end{align}
  It is conceptually helpful to think of $a,b$ as points of the latent
  Poisson process with $\vartheta$ values $x,y$ respectively, but the
  proof does not make formal use of this. The expression
  \[
  \Pr(\degDepLem a>k_{\nu},\degDepLem b>k_{\nu})=\Pr(\degDepLem
  a>k_{\nu})\Pr(\degDepLem b>k_{\nu}|\degDepLem a>k_{\nu}),
  \]
  makes it clear that $g(y)$ is a bound on $\Pr(\degDepLem
  b>k_{\nu}|\degDepLem a>k_{\nu})-\Pr(\degDepLem b>k_{\nu})$.  The
  focus will be on bounding $\Pr(\degDepLem b>k_{\nu}|\degDepLem
  a>k_{\nu})$.  To do this, introduce a marking
  $\{((\theta_{i},\vartheta_{i}),M_{i})\}$ of $\PPrest$ where
  \[
  M_{i}=1[W(x,\vartheta_{i})>U_{i}^{1}]
  \]
  indicates whether each point connects to $a$. This induces the
  obvious marking\footnote{the full marking is defined on $\PPrest$ for consistency
    of the indices of the points $(\theta_{i},\vartheta_{i}).$%
  }
  on $\PPnew$ that
  breaks $\PPnew$ into two independent sets:
  \[
  N_{a}=\{\vartheta_{i}\suchthat(\theta_{i},\vartheta_{i})\in\PPnew,\
  M_{i}=1\},
  \]
  the neighbours of $a$, and
  \[
  \bar{N}_{a}=\{\vartheta_{i}\suchthat(\theta_{i},\vartheta_{i})\in\PPnew,\
  M_{i}=0\},
  \]
  the non-neighbours of $a$. By construction $\abs{N_{a}}=D_{a}$ and
  the neighbours $N_{a}=\{\vartheta_{i}\}_{i=1}^{\degDepLem a}$ are,
  conditional on $\degDepLem a$, independently and identically
  distributed with probability density
  \begin{align}
    \vartheta_{i} & \distiid \frac{W(x,\vartheta_{i})}{\mu_{W}(x)}.
  \end{align}
  The non-neighbours $\bar{N}_{a}$ are a Poisson process on $\NNReals$
  with intensity $\nu(1-W(x,\vartheta))$. The degree of the point $b$
  may be written as the sum of its connections to the neighbours and
  non-neighbours of $a$,
  \[
  \degDepLem b=\degDepLem b^{(N_{a})}+\degDepLem b^{(\bar{N}_{a})},
  \]
  where, by an application of Campbell's theorem,
  \[
  \degDepLem b^{(\bar{N}_{a})}\dist\poiDist(\nu(\mu_{W}(y)-\int
  W(x,z)W(y,z)\intd z))
  \]
  and
  \[
  \degDepLem b^{(N_{a})}\given\degDepLem a\sim\binDist(\degDepLem
  a,p_{x,y})
  \]
  independently, with
  \[
  p_{x,y}=\frac{1}{\mu_{W}(x)}\int W(x,z)W(y,z)\intd z.
  \]
  It is now clear that the dependence of $\degDepLem b$ on $\degDepLem
  a$ comes in only through the number of trials of $\degDepLem
  b^{(N_{a})}\given\degDepLem a$.

  To treat $\degDepLem b^{(N_{a})}$ conditional on the event
  $\degDepLem a>k_{\nu}$ we introduce random variables $L_{1},L_{2}$
  such that on the event $\{\degDepLem a>k_{\nu}\}$
  \[
  L_{1}+L_{2}=\degDepLem a
  \]
  and implicitly specify the joint distribution of $L_{1},L_{2}$ by
  requiring $L_{1}$ to have marginal distribution
  \begin{align}
    L_{1} & \dist \poiDist(\nu\mu_{W}(x))
  \end{align}
  conditional on $\{\degDepLem a>k_{\nu}\}$. Intuitively, $L_{1}$ is
  the number of neighbours of $a$ that would exist without
  conditioning on $\degDepLem a>k_{\nu}$ and $L_{2}$ is the number of
  additional neighbours that are present as a result of the
  conditioning. Therefore on the event $\{\degDepLem a>k_{\nu}\}$
  there are random variables $B_{1},B_{2}$ such that:
  \begin{align}
    \degDepLem b^{(N_{a})} & = B_{1}+B_{2},
  \end{align}
  and
  \begin{align}
    B_{1}\given L_{1} & \dist \binDist(L_{1},p_{x,y})\\
    B_{2}\given L_{2} & \dist \binDist(L_{2},p_{x,y})
  \end{align}
  independently conditional on $L_{1},L_{2}$. The point of introducing
  these auxiliary random now becomes clear as:
  \begin{align}
    B_{1} & \dist\poiDist(\nu\int W(x,z)W(y,z)\intd z)
  \end{align}
  and so
  \begin{align}
    (\degDepLem b^{(\bar{N}_{a})}+B_{1})\given\{\degDepLem a>k_{\nu}\}
    & \dist\poiDist(\nu\mu_{W}(y)).
  \end{align}
  Intuitively, conditional on $\{\degDepLem a>k_{\nu}\}$, $\degDepLem
  b$ splits into a term
  \[
  \pseudoDeg=\degDepLem b^{(\bar{N}_{a})}+B_{1}
  \]
  with the unconditional distribution of $\degDepLem b$ plus a term
  $B_{2}$ that accounts for the 'extra' neighbours of $b$ that one
  expects to see as a result of learning that the degree of $a$ is
  large.

  As $\degDepLem b=\pseudoDeg+B_{2}$,

\[
\Pr(\degDepLem b>k_{\nu}|\degDepLem
a>k_{\nu})=\expect{\Pr(\pseudoDeg+B_{2}>k_{\nu}\given
  L_{1},L_{2})\given\degDepLem a>k_{\nu}}.
\]
Then,
\begin{align}
  \Pr(\pseudoDeg+B_{2}>k_{\nu}\given L_{1},L_{2})=\\
  \Pr(\pseudoDeg>k_{\nu}\given
  L_{1})+\Pr(\pseudoDeg+B_{2}>k_{\nu}\wedge\pseudoDeg\le k_{\nu}\given
  L_{1},L_{2}),
\end{align}
and $L_{1}$ has been defined so that
\[
\expect{\Pr(\pseudoDeg>k_{\nu}\given L_{1})\given\degDepLem
  a>k_{\nu}}=\Pr(\degDepLem b>k_{\nu}).
\]
We have now arrived at
\[
\Pr(\degDepLem a>k_{\nu},\degDepLem b>k_{\nu})=\Pr(\degDepLem
a>k_{\nu})[\Pr(\degDepLem b>k_{\nu})+R],
\]
where the remainder term is
\begin{align}
  R & = \expect{\Pr(\pseudoDeg+B_{2}>k_{\nu}\wedge\pseudoDeg\le k_{\nu}\given L_{1},L_{2})\given\degDepLem a>k_{\nu}}\\
  & = \Pr(\pseudoDeg+B_{2}>k_{\nu}\wedge\pseudoDeg\le
  k_{\nu}\given\degDepLem a>k_{\nu}).
\end{align}
Note that
\begin{align}
  \Pr(\degDepLem a>k_{\nu},\degDepLem b>k_{\nu})-\Pr(\degDepLem
  a>k_{\nu})\Pr(\degDepLem b>k_{\nu}) & =\Pr(\degDepLem a>k_{\nu})R
\end{align}
so that to complete the proof it remains to show that $R\le g(y)$.
For $\nu\mu_{W}(y)$ large the crude bound
\begin{align}
  R & \le \Pr(\pseudoDeg\le k_{\nu}\given\degDepLem a>k_{\nu})\\
  & = \Pr(\degNu[2](y)\le k_{\nu})
\end{align}
suffices. This establishes the claim for $y\le T$ in the lemma
statement.  The remaining task is to find a good bound in the regime
of $y$ where $\nu\mu_{W}(y)$ is not large. In particular, it suffices
to find a bound for $B_{2}$ independent of $\pseudoDeg$ with a
distribution that does not depend on $x$. To that end, let $b>0$ and
write
\[
\Pr(B_{2}>b\given\pseudoDeg)=\expect{\Pr(B_{2}>b\given
  L_{2})\given\pseudoDeg)}.
\]
As $B_{2}\given L_{2}\dist\binDist (L_{2},p_{x,y})$,
\[
\Pr(B_{2}>b\given L_{2})=1-\sum_{n=0}^{b}{L_{2} \choose
  n}p_{x,y}{}^{n}(1-p_{x,y})^{L_{2}-n}.
\]
The salient fact here is that $\varphi(l)=\sum_{n=0}^{b}{l \choose
  n}p(x,y)^{n}(1-p(x,y)^{l-n}$ is a convex function in $l$ and so by a
conditional Jensen's inequality
\begin{equation}
  \Pr(B_{2}>b\given\pseudoDeg)\le\Pr(\tilde{B}>b\given\pseudoDeg),\label{eq:binom_prob_bound}
\end{equation}
where $\tilde{B}\given\pseudoDeg\dist\binDist
(\expect{L_{2}\given\pseudoDeg},p_{x,y}).$ The task is then to find a
bound for the conditional expectation that is independent of
$\pseudoDeg$, which we accomplish by demonstrating a constant bound
$\expect{L_{2}\given\pseudoDeg}\le5k_{\nu}$ for $y$ sufficiently
large. $L_{2}$ is independent of $\pseudoDeg$ conditional on $L_{1}$
so bounding the conditional expectation can be accomplished by
understanding the distribution of $L_{2}\given L_{1}$ and
$L_{1}\given\pseudoDeg$. There exists $Q$ with
\[
Q\equaldist\degDepLem a\given\{\degDepLem a>k_{\nu}\}
\]
and $Q$ independent of $L_{1}$ such that
\begin{align}
  L_{2} & =1[L_{1}\le k_{\nu}](Q-L_{1})\nonumber \\
  \implies\expect{L_{2}\given\pseudoDeg} & \le\Pr(L_{1}\le
  k_{\nu}\given\pseudoDeg)\expect Q.\label{eq:L2_given_H_bound}
\end{align}
This can be understood as the following sampling scheme for a
truncated Poisson distribution:
\begin{enumerate}
\item Draw $l_{1}$ from the Poisson distribution. If $l_{1}>k_{\nu}$
  stop.
\item Otherwise sample $y$ from the truncated distribution, so that
  $l_{1}+(y-l_{1})$ is a trivially a correct sample.
\end{enumerate}
The definitions above can be used to derive:
\[
L_{1}\given\pseudoDeg\dist\binDist (\pseudoDeg,\frac{\int
  W(x,z)W(y,z)\intd x}{\mu_{W}(y)})+Z
\]
where $Z\dist\poiDist(\nu\mu_{W}(x)(1-p_{x,y}))$ is independent of the
first term. Thus,
\[
\Pr(L_{1}\le k_{\nu}\given\pseudoDeg)\expect Q\le\Pr(Z\le
k_{\nu})\expect Q.
\]
Further,
\[
\expect Q<k_{\nu}+\nu\mu_{W}(x),
\]
which can be seen by noting that there is some random variable $G$
such that
\begin{align}
  G & \dist \gammaDist (k_{\nu},1)\given G<\nu\mu_{W}(x)\\
  Q & = k_{\nu}+\poiDist(\nu\mu_{W}(x)-G).
\end{align}
For $\nu\mu_{W}(x)\le2k_{\nu}$, it immediately follows that
\begin{equation}
  \Pr(Z\le k_{\nu})\expect Q\le5k_{\nu}\label{eq:bound_for_le2k}
\end{equation}
For $\nu\mu_{W}(x)>2k_{\nu}$ the assumption $2C\mu(y)\le1-\log2$ for
large enough $y$ implies $\expect Z\ge k_{\nu}$ so a Poisson tail
bound \cite{Glynn_Poisson_Tail_bounds} may be applied to $Z$ to find
\begin{align}
  \Pr(Z\le k_{\nu})\expect Q & \le k_{\nu}+2\Pr(Z=k_{\nu})\nu\mu_{W}(x)\\
  & = k_{\nu}+2\frac{1}{k_{\nu}!}e^{-\nu\mu_{W}(x)(1-p_{x,y})}(\nu\mu_{W}(x)(1-p_{x,y}))^{k_{\nu}}(\nu\mu_{W}(x))\\
  & \le
  k_{\nu}+2\frac{1}{k_{\nu}!}e^{-\nu\mu_{W}(x)(1-C\mu_{W}(y))}(\nu\mu_{W}(x)(1-C\mu_{W}(y)))^{k_{\nu}}(\nu\mu_{W}(x))
\end{align}
The second term satisfies
\[
\frac{2}{k_{\nu}!}e^{-\nu\mu_{W}(x)(1-p_{x,y})}(\nu\mu_{W}(x)(1-p_{x,y}))^{k_{\nu}}(\nu\mu_{W}(x))\\
=2\frac{k_{\nu}+1}{1-C\mu_{W}(y)}\Pr(\tilde{Z}=k_{\nu}+1),
\]
where $\tilde{Z}\dist\poiDist(\nu\mu_{W}(x)(1-C\mu_{W}(y)))$.  This
term is maximized over $\nu\mu_{W}(x)\ge2k_{\nu}$ when
$\expect{\tilde{Z}}$ is minimal, i.e., when $\nu\mu_{W}(x)=2k_{\nu}$.
Subbing in,
\begin{align}
  2\frac{k_{\nu}+1}{1-C\mu_{W}(y)}\Pr(\tilde{Z}=k_{\nu}+1) & \le 2(1-C\mu_{W}(y))^{k_{\nu}}\frac{1}{k_{\nu}!}e^{-2k_{\nu}(1-C\mu_{W}(y))}(2k_{\nu})^{k_{\nu}+1}\\
  & \le 4k_{\nu}2^{k_{\nu}}e^{-k_{\nu}(1-2C\mu_{W}(y))}(\frac{1}{k_{\nu}!}k_{\nu}{}^{k_{\nu}}e^{-k_{\nu}})\\
  & \le 4k_{\nu},
\end{align}
where the final line uses $2C\mu_{W}(y)\le1-\log2$. It then follows
that
\begin{equation}
  \Pr(Z\le k_{\nu})\expect Q\le5k_{\nu}\label{eq:bound_for_g2k}
\end{equation}
for all values of $x$.

Putting together
\cref{eq:binom_prob_bound,eq:L2_given_H_bound,eq:bound_for_le2k,eq:bound_for_g2k}:
\[
\Pr(\pseudoDeg+B_{2}>k_{\nu}\wedge\pseudoDeg\le
k_{\nu}\given\degDepLem
a>k_{\nu})\le\Pr(\pseudoDeg+B(y)>k_{\nu}\wedge\pseudoDeg\le
k_{\nu}\given\degDepLem a>k_{\nu})
\]
where, conditional on $\degDepLem a>k_{\nu}$, $\pseudoDeg$ and $B(y)$
are independent with
\begin{align}
  \pseudoDeg\given\{\degDepLem a>k_{\nu}\} & \equaldist \degNu[1](y)\\
  B(y) & \dist \binDist (5k_{\nu},C\mu_{W}(x)).
\end{align}
This completes the proof of the lemma.
\end{proof}
Roughly speaking, the content of the previous two lemmas amounts to
\begin{align}
  \var{\Ngreater{k_{\nu}}} & \le\expect{\Ngreater{k_{\nu}}}(1+\int
  g(y)\intd y).
\end{align}
That is, the growth of the variance with $\nu$ is controlled by $\int
g(y)\intd y$.  Recalling that our aim is to show
$\var{\Ngreater{k_{\nu}}}=o(\expect{\Ngreater{k_{\nu}}}^{2})$ we must
establish that $\int g(y)\intd y=o(\expect{\Ngreater{0}})$. The
remainder of the proof is devoted to showing this. It turns out that
the appropriate way to do this depends on whether $k_{\nu}$ goes to
infinity.
\begin{lem}
  \label{lem:bounded_k_case}Let $g(y)$ be as in
  \cref{lem:dependence_bound} and suppose $W$ is integrable.  If the sequence $k_{\nu}$ is bounded
  then
  \[
  \int g(y)\intd y=o(\expect{\Ngreater{k_{\nu}}}).
  \]
\end{lem}
\begin{proof}
  Let $T_{\nu}=\sqrt{\expect{\Ngreater{k_{\nu}}}}$ so that by
  \cref{lem:dependence_bound} for $\nu$ large enough
  \[
  \int_{\NNReals}g(y)\intd y\le
  T_{\nu}+\int_{T_{\nu}}^{\infty}\Pr(\degNu[2](y)+B(y)>k_{\nu}\wedge\degNu[2](y)\le
  k_{\nu})\intd y.
  \]
  Moreover
  \[
  \Pr(\degNu[2](y)+B(y)>k_{\nu}\wedge\degNu[2](y)\le
  k_{\nu})\le\Pr(\tilde{B}(y)>1),
  \]
  where, letting $k=\lim_{\nu\to\infty}k_{\nu}$,
  $\tilde{B}(y)\dist\binDist (5k,C\mu_{W}(y))$.  By Markov's
  inequality
  \begin{align}
    \Pr(\tilde{B}(y)>1) & \le5kC\mu_{W}(y)
  \end{align}
  so that
  \begin{align}
    \int_{T_{\nu}}^{\infty}\Pr(\tilde{B}(y)>1)\intd y & \le5kC\int_{T_{\nu}}^{\infty}\mu_{W}(y)\intd y\\
    & =o(1),
  \end{align}
  where the final line follows by the integrability of $\mu_{W}$. Thus
  $\int_{\NNReals}g(y)\intd y=O(\sqrt{\expect{\Ngreater{k_{\nu}}}})$.
\end{proof}
The case $k\upto\infty$ is substantially trickier. Essentially the
strategy here is to break up to domain of $y$ into three components
and use a different tractable and reasonably tight bound on $g(y)$ in
each region, see \cref{tab:Upper-bounds-on-g}. An important
intermediate step is the observation
\begin{align}
  \expect{\Ngreater{0}} & =\Omega(\mu_{W}^{-1}(\frac{1}{\nu})),
\end{align}
which will eventually allow us to show $\int g(y)\intd
y=o(\expect{\Ngreater{0}})$ by establishing bounds on the integral in
terms of $\mu_{W}^{-1}(\frac{1}{\nu})$.  For instance, the next lemma
can be understood as establishing that
$\int_{0}^{\mu_{W}^{-1}((1+\epsilon)\frac{k_{\nu}}{\nu})}\Pr(\degNu[2](y)\le
k_{\nu})\intd y$ is at most an exponentially vanishing (in $k_{\nu}$)
fraction of $\expect{\Ngreater{0}}$.

\begin{table}
  \begin{centering}
    \begin{tabular}{|c|c|}
      \hline 
      Region of $\NNReals$ & Upper bound for $g(y)$\tabularnewline
      \hline 
      \hline 
      $[0,\mu^{-1}((1+\epsilon)\frac{k_\nu}{\nu})]$ & $\Pr(\degNu[2](y)\le k_{\nu})$\tabularnewline
      \hline 
      $(\mu^{-1}((1+\epsilon)\frac{k_\nu}{\nu}),\mu^{-1}((1-\epsilon)\frac{k}{\nu}))$ & $1$\tabularnewline
      \hline 
      $(\mu^{-1}((1-\epsilon)\frac{k_\nu}{\nu}),\infty)$ & $\Pr(B(y)>\frac{\epsilon}{2}k_{\nu})+\Pr(\degNu[2](y)>(1-\frac{\epsilon}{2})k_{\nu})$\tabularnewline
      \hline 
    \end{tabular}
    \par\end{centering}

  \caption{Upper bounds on $g(y)$\label{tab:Upper-bounds-on-g}}
\end{table}

\begin{lem}\label{lem:deg_bound_small_front}
For $0<\epsilon<1$,
  \[
  \int_{0}^{\mu_{W}^{-1}((1+\epsilon)\frac{k_{\nu}}{\nu})}\Pr(\degNu[2](y)\le
  k_{\nu})\intd
  y\le\frac{1+\epsilon}{\epsilon}(\frac{1+\epsilon}{e^{\epsilon}})^{k_{\nu}}\mu_{W}^{-1}(\frac{k_{\nu}}{\nu}).
  \]
\end{lem}
\begin{proof}
  Because $\Pr(\degNu[2](y)\le k_{\nu})$ is monotonically increasing
  in $y$ over the domain of integration, the integral is bounded by
  \begin{align}
    \mu_{W}^{-1}((1+\epsilon)\frac{k_{\nu}}{\nu})\Pr(\degNu[2](\mu_{W}^{-1}((1+\epsilon)\frac{k_{\nu}}{\nu}))\le
    k_{\nu}).
  \end{align}
  As
  $\expect{\degNu[2](\mu_{W}^{-1}((1+\epsilon)\frac{k_{\nu}}{\nu}))}=(1+\epsilon)k_{\nu}>k_{\nu}$
  a tail bound \cite{Glynn_Poisson_Tail_bounds} applies:
  \begin{align}
    \Pr(\degNu[2](\mu_{W}^{-1}((1+\epsilon)\frac{k_{\nu}}{\nu}))\le k_{\nu}) & \le(1+\frac{1}{\epsilon})\Pr(\degNu[2](y)=k_{\nu})\\
    & =(1+\frac{1}{\epsilon})\frac{1}{k_{\nu}!}((1+\epsilon)k_{\nu})^{k_{\nu}}e^{-(1+\epsilon)k_{\nu}}\\
    &
    \le\frac{1}{e}\frac{1+\epsilon}{\epsilon}(\frac{1+\epsilon}{e^{\epsilon}})^{k_{\nu}}.
  \end{align}
 
\end{proof}

For $y>\mu_{W}^{-1}((1-\epsilon)\frac{k_{\nu}}{\nu})$ we can bound
$g(y)$ (and thus
$\int_{\mu_{W}^{-1}((1-\epsilon)\frac{k_{\nu}}{\nu})}^\infty g(y)
\intd y$) by
\begin{align}
  &\Pr(\degNu[2](y)+B(y)>k_{\nu}\wedge\degNu[2](y)\le k_{\nu}) \\
  &\le\Pr(\degNu[2](y)+B(y)>k_{\nu}\wedge\degNu[2](y)\le(1-\frac{\epsilon}{2})k_{\nu})\\
  & +\Pr(\degNu[2](y)>(1-\frac{\epsilon}{2})k_{\nu})\\
  &
  \le\Pr(B(y)>\frac{\epsilon}{2}k_{\nu})+\Pr(\degNu[2](y)>(1-\frac{\epsilon}{2})k_{\nu}).
\end{align}
The next lemma controls the second term in this bound.
\begin{lem}
  \label{lem:deg_bound_small_tail}Suppose there is some $\chi>0$ such
  that for all $x>\chi$ it holds that
  \begin{align}
    \frac{\mu_{W}(x)}{x\mu_{W}'(x)} & \ge-1,
  \end{align}
  then, for $\nu$ sufficiently large such that
  $\frac{k_{\nu}}{\nu}\le\mu_{W}(\chi)$ and $\epsilon$ such that
  $0<\epsilon<1$,
  \[
  \int_{\mu_{W}^{-1}((1-\epsilon)\frac{k_{\nu}}{\nu})}^{\infty}\Pr(\degNu[2](y)>(1-\frac{\epsilon}{2})k_{\nu})\intd
  y\le\frac{2}{k_{\nu}\epsilon+2}\mu_{W}^{-1}(\frac{k_{\nu}}{\nu})
  \]
\end{lem}
\begin{proof}
  For $y\in[\mu_{W}^{-1}((1-\epsilon)\frac{k_{\nu}}{\nu}),\infty)$ it
  holds that $\expect{\degNu[2](y)}<(1-\epsilon/2)k_{\nu}$ so a tail
  bound \cite{Glynn_Poisson_Tail_bounds} applies:
  \[
  \Pr(\degNu[2](y)>(1-\frac{\epsilon}{2})k_{\nu})\le(\frac{1-\epsilon/2+1/k_{\nu}}{\epsilon/2+1/k_{\nu}})\frac{1}{\floor{(1-\frac{\epsilon}{2})k_{\nu}}!}e^{-\nu\mu_{W}(y)}(\nu\mu_{W}(y))^{(1-\frac{\epsilon}{2})k_{\nu}}.
  \]
  Because $\mu_{W}(y)$ is strictly monotonic the component of the
  bound that depends on $y$ may be integrated by substitution. For
  notational simplicity, let $f(x)=\mu_{W}^{-1}(y)$, then
  \begin{align}
    \int_{\mu_{W}^{-1}((1-\epsilon)\frac{k_{\nu}}{\nu})}^{\infty}e^{-\nu\mu_{W}(y)}(\nu\mu_{W}(y))^{(1-\frac{\epsilon}{2})k_{\nu}}\intd
    y &
    =-\int_{0}^{(1-\epsilon)k_{\nu}}e^{-x}x^{(1-\frac{\epsilon}{2})k_{\nu}}\frac{1}{\nu}f'(\frac{x}{\nu})\intd
    x.
  \end{align}
  Let $z=f(x)$ and write
  \begin{align}
    \frac{\mu_{W}(z)}{z\mu_{W}'(z)} & =\frac{f'(x)x}{f(x)}
  \end{align}
  so by assumption for $x\le\mu_{W}(\chi)$ holds that
  $x\frac{f'(x)}{f(x)}\ge-1$.  Thus for $\nu$ sufficiently large that
  $\frac{k_{\nu}}{\nu}\le\mu_{W}(\chi)$ it holds that
  \begin{align}
    -\int_{0}^{(1-\epsilon)k_{\nu}}e^{-x}x^{(1-\frac{\epsilon}{2})k_{\nu}}\frac{1}{\nu}f'(\frac{x}{\nu})\intd
    x &
    \le\int_{0}^{(1-\epsilon)k_{\nu}}e^{-x}x^{(1-\frac{\epsilon}{2})k_{\nu}-1}f(\frac{x}{\nu})\intd
    x.
  \end{align}
  Moreover, $xf(x)$ is a monotonically non-decreasing function on
  $x\le\mu_{W}(\chi)$, which may be established by:
  \begin{align}
    (xf(x))' & =f(x)+xf'(x)\\
    & =f(x)(1+x\frac{f'(x)}{f(x)})\\
    & \ge0.
  \end{align}
  This implies
  \begin{align}
    \int_{0}^{(1-\epsilon)k_{\nu}}e^{-x}x^{(1-\frac{\epsilon}{2})k_{\nu}-1}f(\frac{x}{\nu})\intd x & \le(1-\epsilon)k_{\nu}f(\frac{k_{\nu}}{\nu})\int_{0}^{(1-\epsilon)k_{\nu}}e^{-x}x^{(1-\frac{\epsilon}{2})k_{\nu}-2}\intd x\\
    & \le(1-\epsilon)k_{\nu}f(\frac{k_{\nu}}{\nu})\Gamma((1-\frac{\epsilon}{2})k_{\nu}-1)\\
    & =f(\frac{k_{\nu}}{\nu})\Gamma((1-\frac{\epsilon}{2})k_{\nu}).
  \end{align}
  This establishes
  \begin{align}
    \int_{\mu_{W}^{-1}((1-\epsilon)\frac{k_{\nu}}{\nu})}^{\infty}\Pr(\degNu[2](y)>(1-\frac{\epsilon}{2})k_{\nu})\intd y\le & (\frac{1-\epsilon/2+1/k_{\nu}}{\epsilon/2+1/k_{\nu}})f(\frac{k_{\nu}}{\nu})\frac{\Gamma((1-\frac{\epsilon}{2})k_{\nu})}{\Gamma((1-\frac{\epsilon}{2})k_{\nu}+1)}\\
    = &
    \frac{1}{\epsilon/2+1/k_{\nu}}\frac{1}{k_{\nu}}f(\frac{k_{\nu}}{\nu})
  \end{align}
  as claimed.
\end{proof}
The next lemma establishes the other half of the tail bound for
$g(y)$:
\begin{lem}
  \label{lem:g_bound_upper_tail_b_part}Suppose there is some $\chi>0$
  such that, for all $x>\chi$,
  \begin{align}
    \frac{\mu_{W}(x)}{x\mu_{W}'(x)} & \ge-1,
  \end{align}
  and let $B$ and $C$ be as in \cref{lem:dependence_bound}.  For $\nu$
  sufficiently large such that $\frac{k_{\nu}}{\nu}\le\mu_{W}(\chi)$
  and $\epsilon$ such that $10C\frac{k_{\nu}}{\nu}\le\epsilon<1$,

  \[
  \int_{\mu_{W}^{-1}((1-\epsilon)\frac{k_{\nu}}{\nu})}^{\infty}\Pr(B>\frac{\epsilon}{2}k_{\nu})\intd
  y
  \le(\frac{C}{10}\frac{\epsilon}{1-\epsilon}\frac{k_{\nu}}{\nu})^{\epsilon
    k_{\nu}/2}\frac{1}{\epsilon k_{\nu}/2-1}
  \mu_W^{-1}(\frac{(1-\epsilon)k_{\nu}}{\nu}).
  \]
\end{lem}
\begin{proof}
  The condition $10C\frac{k_{\nu}}{\nu}\le\epsilon$ ensures that
  \begin{align}
    C\mu_{W}(y) & <\frac{\epsilon/2k_{\nu}}{5k_{\nu}},
  \end{align}
  for $y>\mu_{W}^{-1}((1-\epsilon)\frac{k_{\nu}}{\nu})$. Recalling
  $B\dist\binDist(5k_{\nu},C\mu_{W}(y))$, this allows a large
  deviation bound
  \cite{Arratia_Gordon_binomial_large_deviation_bounds} to be applied:
  \begin{align}
    \Pr(B>\frac{\epsilon}{2}k_{\nu}) &
    \le\exp(-5k_{\nu}S(\frac{\epsilon/2k_{\nu}}{5k_{\nu}}\|C\mu_{W}(y))),
  \end{align}
  where $S(q\|p)=q\log\frac{q}{p}+(1-q)\log\frac{1-q}{1-p}$ is the
  relative entropy between $\bern(q)$ and $\bern(p)$.
  \begin{align}
    S(\frac{\epsilon}{10}\|C\mu_{W}(y)) &
    \ge\frac{\epsilon}{10}\log\frac{10}{C\epsilon}\frac{1}{\mu_{W}(y)},
  \end{align}
  whence
  \begin{align}
    \Pr(B>\frac{\epsilon}{2}k_{\nu}) &
    \le(\frac{C}{10}\epsilon)^{\epsilon
      k/2}\mu_{W}(y)^{\frac{\frak{\epsilon k}}{2}}.
  \end{align}
  It remains to integrate this bound. Let $f(x)=\mu_{W}^{-1}(x)$ then
  \begin{align}
    \int_{\mu_{W}^{-1}((1-\epsilon)\frac{k_{\nu}}{\nu})}^{\infty}\mu_{W}(y)^{\epsilon
      k_{\nu}/2}\intd y & =\nu^{-\epsilon
      k/2}\int_{0}^{(1-\epsilon)k_{\nu}}x^{\epsilon
      k_{\nu}/2}\frac{1}{\nu}f'(\frac{x}{\nu})\intd x.
  \end{align}
  Following the same reasoning as in the proof of
  \cref{lem:deg_bound_small_tail},
  \[x^{2}\frac{1}{\nu}f'(\frac{x}{\nu})\le(1-\epsilon)k_{\nu}f(\frac{(1-\epsilon)k_{\nu}}{\nu})\]
  on the domain of integration so,
  \begin{align}
    \nu^{-\epsilon k/2}\int_{0}^{(1-\epsilon)k_{\nu}}x^{\epsilon k/2}\frac{1}{\nu}f'(\frac{x}{\nu})\intd x & \le\nu^{-\epsilon k/2}(1-\epsilon)k_{\nu}f(\frac{(1-\epsilon)k_{\nu}}{\nu})[\frac{1}{\epsilon k_{\nu}/2-1}((1-\epsilon)k_{\nu})^{\epsilon k_{\nu}/2-1}]\\
    & =(\frac{k_{\nu}}{\nu})^{\epsilon
      k_{\nu}/2}(1-\epsilon)^{\epsilon k_{\nu}/2}\frac{1}{\epsilon
      k_{\nu}/2-1}f(\frac{(1-\epsilon)k_{\nu}}{\nu}).
  \end{align}

\end{proof}
In particular, the last several lemmas combine to show that for
$\epsilon_{\nu}\le1$ such that
$\epsilon_{\nu}=\omega(\frac{1}{k_{\nu}})$ and
$\epsilon_{\nu}=\omega(\frac{k_{\nu}}{\nu})$ it holds that
\begin{align}
  \int_{0}^{\mu_{W}^{-1}((1+\epsilon)\frac{k_{\nu}}{\nu})}g(y)\intd
  y+\int_{\mu_{W}^{-1}((1-\epsilon)\frac{k_{\nu}}{\nu})}^{\infty}g(y)\intd
  y & =o(\mu_{W}^{-1}(\frac{1}{\nu})).
\end{align}
With the observation that
$\expect{\Ngreater{0}}=\Omega(\mu_{W}^{-1}(\frac{1}{\nu}))$ this
leaves only the region
\[
(\mu_{W}^{-1}(1+\epsilon)\frac{k_{\nu}}{\nu},\mu_{W}^{-1}(1-\epsilon)\frac{k_{\nu}}{\nu})
\]
as a possible foil to $\int g(y)\intd y=o(\expect{\Ngreater{0}})$. In
this regime we expect
\begin{align}
  g(y) & =\Pr(\degNu[2](y)+B(y)>k_{\nu}\wedge\degNu[2](y)\le k_{\nu})
\end{align}
to be approximately constant because $\expect{\degNu[2](y)}\approx
k_{\nu}$ so we make due with the bound $g(y)\le1$.
\begin{lem}
  \label{lem:g_bound_middle_region}
  Suppose that $\mu_{W}$ is differentiable and that there is some
  $\chi>0$ such that for all $x>\chi$ it holds that
  \begin{align}
    \frac{\mu_{W}(x)}{x\mu_{W}'(x)} & \ge-1.
  \end{align}
  Then for $\epsilon>0$ and $\nu$ sufficiently large such that
  $(1+\epsilon)\frac{k_{\nu}}{\nu}\le\mu_{W}(\chi)$, it holds that
  \begin{align}
    \mu_{W}^{-1}((1-\epsilon)\frac{k_{\nu}}{\nu})-\mu_{W}^{-1}((1+\epsilon)\frac{k_{\nu}}{\nu})
    & \le
    2\frac{\epsilon}{1-\epsilon}\mu_{W}^{-1}((1-\epsilon)\frac{k_{\nu}}{\nu})
  \end{align}
\end{lem}
\begin{proof}
  Let $f(x)=\mu_{W}^{-1}(x)$. Since $\mu_{W}$ is differentiable so is
  $f$. By the mean value theorem there is some point
  $(1-\epsilon)\frac{k_{\nu}}{\nu}\le
  x^{*}\le(1+\epsilon)\frac{k_{\nu}}{\nu}$ such that
  \begin{align}
    f((1-\epsilon)\frac{k_{\nu}}{\nu})-f((1+\epsilon)\frac{k_{\nu}}{\nu}) & =-2\epsilon\frac{k_{\nu}}{\nu}f'(x^{*})\\
    & =-2\epsilon\frac{k_{\nu}}{\nu}\frac{1}{x^{*}}x^{*}f'(x^{*})\\
    &
    \le2\frac{\epsilon}{1-\epsilon}f((1-\epsilon)\frac{k_{\nu}}{\nu}),
  \end{align}
  where the final line follows as in \cref{lem:deg_bound_small_tail}.
\end{proof}
We can now complete our intermediate goal:
\begin{lem}
  \label{lem:g_bound}Let $g(y),\ T$ and $C$ be as in
  \cref{lem:dependence_bound}.  Suppose $k_{\nu}\upto\infty$ and
  $k_{\nu}=o(\nu).$ Suppose that $\mu_{W}$ is differentiable and that
  there is some $\chi>0$ such that for all $x>\chi$ it holds that
  \begin{align}
    \frac{\mu_{W}(x)}{x\mu_{W}'(x)} & \ge-1.
  \end{align}
  Then
  \begin{align}
    \int g(y)\intd y & =o(\expect{\Ngreater{0}})
  \end{align}
\end{lem}
\begin{proof}
  Let $\epsilon_{\nu}\downto0$ such that
  $\epsilon_{\nu}=\omega(\sqrt{\frac{1}{k_{\nu}}})$ and
  $\epsilon_{\nu}=\omega(\sqrt{\frac{k_{\nu}}{\nu}})$. Let
  \begin{align}
    h(y) & =\begin{cases}
      \Pr(\degNu[2](y)\le k_{\nu}) & y\le\mu_{W}^{-1}((1+\epsilon_{\nu})\frac{k_{\nu}}{\nu})\\
      1 & y\in(\mu_{W}^{-1}((1+\epsilon_{\nu})\frac{k_{\nu}}{\nu}),\mu_{W}^{-1}((1-\epsilon_{\nu})\frac{k_{\nu}}{\nu}))\\
      \Pr(B(y)>\frac{\epsilon}{2}k_{\nu})+\Pr(\degNu[2](y)>(1-\frac{\epsilon}{2})k_{\nu})
      & y\ge\mu_{W}^{-1}((1-\epsilon_{\nu})\frac{k_{\nu}}{\nu}).
    \end{cases}
  \end{align}
  Because $\mu_W$ is not compactly supported, for $\nu$ sufficiently
  large $\mu_{W}^{-1}((1+\epsilon_{\nu})\frac{k_{\nu}}{\nu})>T$ and in
  this regime it is immediate that
  \begin{align}
    g(y) & \le h(y).
  \end{align}
  Moreover, it is straightforward to verify that the conditions on
  $\epsilon_{\nu}$ with
  \cref{lem:deg_bound_small_front,lem:deg_bound_small_tail,lem:g_bound_upper_tail_b_part,lem:g_bound_middle_region}
  imply
  \begin{align}
    \int h(y)\intd y &
    =o(\mu_{W}^{-1}((1-\epsilon_{\nu})\frac{k_{\nu}}{\nu})).
  \end{align}
  (For \cref{lem:deg_bound_small_front} it suffices to consider the
  worst case $\epsilon_\nu=\sqrt{\frac 1 {k_{\nu}}}$.)

  Next,
  \begin{align}
    \expect{\Ngreater{0}} & =\int_{\NNReals}1-e^{-\nu\mu_{W}(y)}\intd y\\
    & \ge\int_{0}^{\mu_{w}^{-1}(\frac{1}{\nu})}1-e^{-1}\intd y\\
    & =\Omega(\mu_{W}^{-1}(\frac{1}{\nu})).
  \end{align}
  Thus $\expect{\Ngreater{0}}
  =\Omega(\mu_{W}^{-1}((1-\epsilon_{\nu})\frac{k_{\nu}}{\nu}))$,
  completing the proof.
\end{proof}

We are now equipped to give the proof of the main result:
\begin{proof}[Proof of \cref{thm:deg_dist_theorem}]
  By \cref{lem:restric_approx_valid} it suffices to show
  $\var{\Ngreater{k_{\nu}}}=o(\expect{\Ngreater{0}}^{2})$.  By
  \cref{lem:var_expression,lem:dependence_bound},
  \begin{align}
    \var{\Ngreater{k_{\nu}}} & \le\expect{\Ngreater{k_{\nu}}}(1+\int
    g(y)\intd y),
  \end{align}
  where $g(y)$ is as defined in
  \cref{lem:dependence_bound}. \cref{lem:bounded_k_case}, for bounded
  $k_\nu$, and \cref{lem:g_bound}, for $k_\nu\upto\infty$, establish
  \begin{align}
    \int g(y)\intd y & =o(\expect{\Ngreater{0}}),
  \end{align}
  completing the proof.
\end{proof}

\section{Connectivity for Separable KEGs\label{sec:Connectivity}}

A serious omission in the results presented thus far is that they give
virtually no information about the global structure of the KEGs.  In
particular, we have as yet made no statements about the connectivity
structure of these graphs. The sparse structure that we explore here
could, in principle, arise from graphs that consist of large numbers
of disconnected dense components. If this were to be the case then
these graphs would be uninteresting for physical applications. Our aim
in this section is to give a preliminary result showing that this is
not the case. \begin{defn} We call a KEG \defnphrase{separable} if the
  associated graphex has $\IsoF=\StarF=0$ and $W$ of the form 
  \[
  W(x,y)=\begin{cases}
    0 & x=y\\
    f(x)f(y) & \mbox{otherwise.}
  \end{cases}
  \]
  We prove that separable KEGs have an arbitrarily large fraction of
  the vertices contained in a single connected component in the large
  graph limit. (As usual, because there is no risk of confusion, we will use the term graphex to refer to the function $W$. )

\begin{rem}
Separability in combination with the graphex integrability conditions immediately implies that $f$ and hence $W$ is integrable
and thus that this result only applies for graphs that have a finite expected number of edges
when restricted to finite support $\nu$.
\end{rem}

  The main obstacle to the study of connectivity in the KEG setting is
  that the graphs are naturally defined in terms of the infinite
  collection of points in the latent Poisson process with only a
  finite number of these participating as points in a sampled
  graph. The difficulty is that traditional tools (e.g.\
  \cite{bollobas:2001}) for studying connectivity begin with a fixed
  set of vertices of the graph and examine how they become connected as
  edges are randomly introduced, an approach that is apparently futile
  in the present setting where we must specify the edge set in
  order to specify the vertex set. The tactic we use to circumvent
  this problem hinges on the division of the KEG into three parts
  based on the latent $\vartheta$ values of the vertices: the induced
  subgraph below some threshold value, the induced subgraph above this
  threshold and the bi-graph between them; see
  \cref{fig:sep_KEG_structure}.  The first piece intuition is that for
  fixed $\nu$ we can set the threshold $T_{\nu}$ such that nearly
  every point of the latent Poisson process with $\vartheta$ below
  $T_{\nu}$ will have an edge connected to it; because of this we can
  treat the connectivity of the below $T_{\nu}$ induced subgraph using
  the traditional random graph machinery.  The connectivity of vertices
  lying above $T_{\nu}$ that participate in at least one edge
  connecting below $T_{\nu}$ then follows straightforwardly.  This
  leaves only the vertices in the induced subgraph above $T_{\nu}$ that
  do not connect to a point below $T_{\nu}$ and it will turn out that
  these constitute a negligible fraction of the graph.

  \begin{figure}
    \begin{centering}
      \includegraphics[width=.25\linewidth]{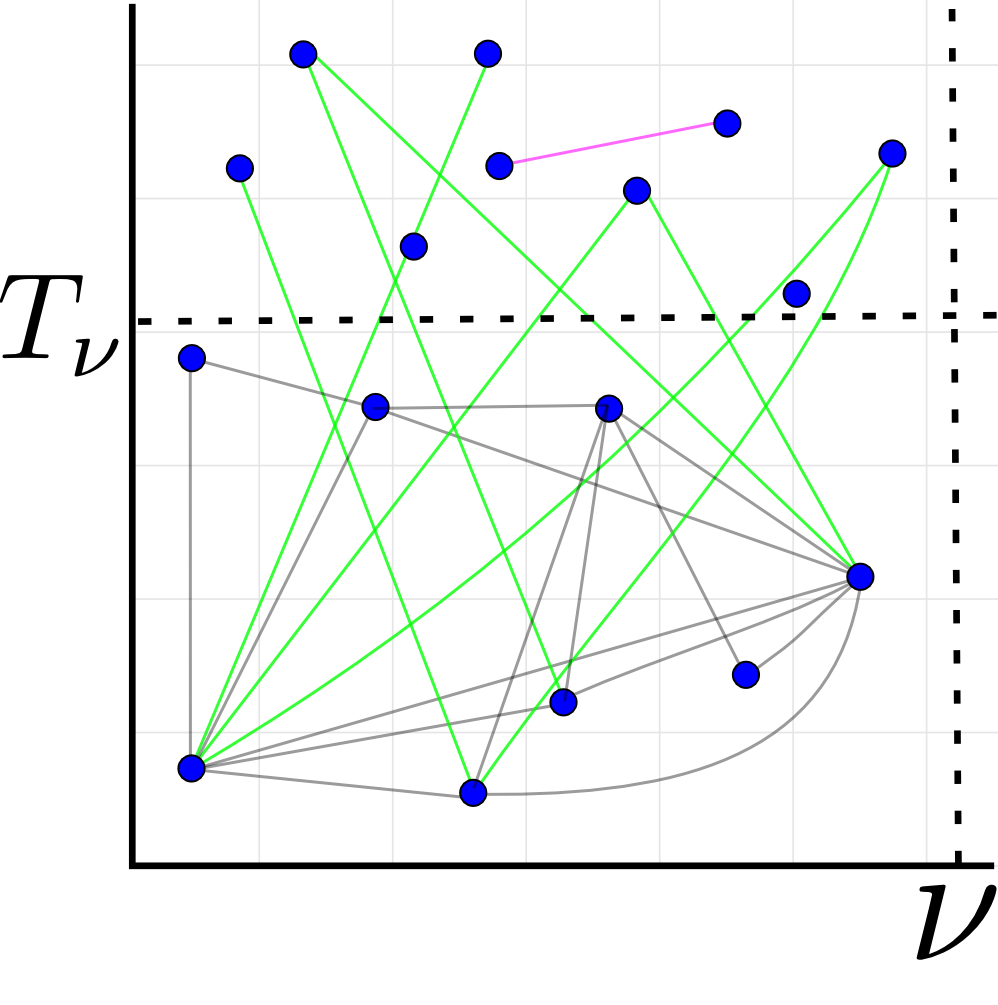}
      \par\end{centering}

    \caption{\label{fig:sep_KEG_structure}The basic structure of
      separable KEGs.  The induced subgraph below $T_{\nu}$ in gray is
      fully connected.  Above $T_{\nu}$ the vast majority of the vertices
      of the graph connect to the below threshold subgraph, in
      green. This leaves only the very small number of vertices connected
      only to vertices 
      that lie entirely above $T_{\nu}$, in magenta.}
  \end{figure}

  We fix some notation that we will need for the rest of this section:
  Let $\PP$ be the unit rate Poisson process on $\NNReals^{2}$ and let
  $\PP_{\nu}=\{
  (\theta_{i},\vartheta_{i})\in\Pi\given\theta_{i}\le\nu\} $ be the
  restriction of this process to label-space $\le\nu$. Let the Poisson
  process below a cutoff value $x$ in $\vartheta$ space be $\PPbelowth
  x=\{ (\theta_{i},\vartheta_{i})\in\PPnu\given\vartheta_{i}<x\} $ and
  let the process above the cutoff be $\PPaboveth x=\{
  (\theta_{i},\vartheta_{i})\in\PPnu\given\vartheta_{i}>x\} $.

  We begin by showing we can take $f(x)$ to be monotone decreasing
  without loss of generality:
  \begin{lem} Let $W(x,y)=f(x)f(y)1[x\neq y]$ be a separable graphex,
    then there is some other separable graphex $W^{'}=h(x)h(y)1[x\neq
    y]$ such that $h$ is monotone decreasing and the KEGs associated
    to $W$ and $W^{'}$ are equal in distribution.
  \end{lem}
  \begin{proof}
    Because the distribution of a KEG is invariant under measure
    preserving transformations of the generating graphon, it suffices
    to show that there are some measure preserving transformations
    $\tau , \varphi : \NNReals\to\NNReals$ and a monotonically
    decreasing function $h$ such that $f\circ\tau=h\circ \varphi$

    If $f(x)$ has bounded domain (i.e., $W$ is a graphon) then the
    result follows immediately from \citep[][Prop.A19]{Lovasz:2013:A},
    which shows that for any bounded $f$ with compact support there is
    some measure preserving transformation $\varphi$ on the domain of
    $f$ and monotone decreasing $h$ such that $f=h\circ\varphi$.

    Assume $f(x)$ has unbounded domain.
    Because $f$ is integrable and measurable the sets
    $A_k=\{x\suchthat f(x)\in [\frac 1 k, \frac{1}{k+1})\}$ for
    $k\in\Nats$ are Borel sets of finite measure. This means in
    particular (\citep[][Thm. A.20]{Kerr:2014}) that for $A_k$ with
    measure $c_k$ there is some measure preserving transformation
    $\tilde{\tau}$ such that $\tilde{\tau}(A_k)=[0,c_k]$. From this it
    immediately follows that there exists a measure preserving
    transformation $\tau$ such that $\tau(A_k)=[c_{k-1},c_k]$ with
    $c_0=0$.
    That is, $\tau$ imposes a pseudo-monotonicity where $f(\tau(x)) <
    \frac 1 k$ and $f(\tau(y)) \ge \frac 1 k$ implies
    $\tau(x)>\tau(y)$. By \citep[][Prop.A19]{Lovasz:2013:A} there is a
    measure preserving transformation $\varphi_k$ and a monotonically
    decreasing $h_k$ with support $\tau(A_k)$ such that
    $1_{\tau(A_k)}f\circ\tau=h_k\otimes\varphi_k$. Letting $\varphi =
    \bigotimes_i \varphi_i$ and $h=\bigotimes_i h_i$ completes the
    proof.
  \end{proof}

  We take $f$ to be monotone decreasing for the remainder of the
  section. Because the result is trivial for $f$ with bounded domain
  (the KEG is dense) we also take $f$ to have unbounded domain.
  Denote the left continuous inverse of $f$ by $f\linverse(t)=\inf\{
  \lambda\st f(\lambda)=t\} $. We will make frequent use of the
  observation that for $l_\nu\in o(1)$ it holds that
  $f\linverse(l_\nu)\in\omega(1)$. Let $G$ be a Kallenberg
  Exchangeable Graph associated with $W$ and let $G_{\nu}$ be the
  restriction to $\left[0,\nu\right]$.
\end{defn}

\begin{defn}
  Let $t_{\nu}$ be a function of $\nu$ such that $t_{\nu}\in o(1)$ and
  $t_\nu\in\omega(\frac 1 \nu)$ and define the \defnphrase{threshold}
  $\threshold=f\linverse(\frac{1}{\nu}+t_{\nu})$. \end{defn}
\begin{rem}
  This notation for the threshold suppresses the dependence on
  $t_{\nu}$, which should be thought of as going to $0$ as quickly as
  possible consistent with $t_{\nu}\in\omega(\frac{1}{\nu})$.
\end{rem}
The proof now proceeds roughly as follows:
\begin{enumerate}
\item We establish the existence of a connected core that we will show
  nearly every vertex of the graph connects to
  (\cref{lem:Erdos_Renyi_near_0})
\item We show that nearly every point of $\PPbelowth{\threshold}$
  participates in an edge connecting
  to the connected core
  (\cref{lem:all_below_threshold})
\item We lower bound the number of points of $\PPaboveth{\threshold}$
  that connect to the connected core (\cref{lem:visible_above})
\item We consider the induced subgraph of $G_{\nu}$ given by $\{
  \theta_{i}\in\vertexset{G_{\nu}}\given\vartheta_{i}>\threshold\} $
  and show that the number of points in this subgraph that fail to
  connect to the connected core is an arbitrarily small fraction of
  the number of vertices in the graph (\cref{lem:invisible_above})
\end{enumerate}
The first step of the proof is to show that there is an induced
subgraph $\popgraph$ that is both connected and very popular in the
sense that every other vertex of the graph will connect to it with high
probability.  The notion of popularity that we use is the that total
mass in the subgraph, $\sum_{p\in\popgraph}f(p)$, is an arbitrarily
large fraction of the total expected mass in the entire graph:
$\expect{\sum_{\vartheta_{i}\in\Pi_{\nu}}f(\vartheta_{i})}=\nu\norm
f_{1}$.  The critical fact for use in later parts of the argument
turns out to be that the mass of the popularity subgraph scales as
$\nu$.
\begin{lem}
  \label{lem:Erdos_Renyi_near_0} Suppose $f$ does not have compact
  support. Let $\popthreshold=f\linverse(\sqrt{\frac{\log\nu}{\nu}})$
  and let $\popgraph$ be the induced subgraph of $G_{\nu}$ given by
  including only vertices in $\Pi_{\nu,<T_{\nu,\mbox{pop}}}$, then:
  \begin{enumerate}
  \item Every element of $\PPbelowth{\popthreshold}$ connects to an
    edge;
    $\lim_{\nu\to\infty}\abs{\PPbelowth{\popthreshold}\backslash\vertexset{\popgraph}}=0\as$
  \item $\popgraph$ is almost surely connected; let $C(\popgraph)=1$
    if $\popgraph$ is connected and $0$ otherwise, then
    $\lim_{\nu\to\infty}C(\popgraph)=1\as$
  \item $\popgraph$ is ``ultra-popular'' almost surely; letting
    $S_{\nu}=\sum_{p\in\popgraph}f(p)$ we have for $\epsilon>0$ that
    $\lim_{\nu\to\infty}\frac{S_{\nu}}{\nu}\ge(1-\epsilon)\|f\|_{1}\as$
  \end{enumerate}
\end{lem}
\begin{proof}
  The key insight is that the connection probabilities below
  $\popthreshold$ are lower bounded by
  $p_{\nu}=f(\popthreshold)^{2}=\frac{\log\nu}{\nu}$ so that a
  sufficient condition for claims 1 and 2 is that the \ErdosRenyi\
  random graph $G(N_{\nu},p_{\nu})$ with
  $N_{\nu}\sim\poiDist(\nu\popthreshold)$ is almost surely connected
  in the limit. A sufficient condition \cite{bollobas:2001} for this
  is that there exists some $\delta>0$ such that
  \[
  \lim_{\nu\to\infty}\frac{p_{\nu}}{\log N_{\nu}/N_{\nu}}>1+\delta\as
  \]
  For arbitrary $\gamma>0$, it holds that
  $\lim_{\nu\to\infty}N_{\nu}/\nu\popthreshold\ge(1-\gamma)\mbox{
    a.s.}$ and so we have that:
  \begin{align}
    \lim_{\nu\to\infty}\frac{p_{\nu}}{\log N_{\nu}/N_{\nu}} & \ge \lim_{\nu\to\infty}\frac{\log\nu/\nu}{\log(1-\gamma)\nu\popthreshold/(1-\gamma)\nu\popthreshold}\mbox{ a.s.}\\
    & = \infty.
  \end{align}
  Thus in the limit as $\nu\to\infty$, the random graph with vertices
  $\PPbelowth{\popthreshold}$ and independent edge probabilities
  $f(\vartheta_{i})f(\vartheta_{j})$ is connected and, in particular,
  every vertex is contained in an edge, thereby establishing claims $1$
  and $2$.

  It remains to show that $S_{\nu}$ grows as claimed. For $\gamma>0$,
  by Hoeffding's inequality we have:
  \begin{align}
    \Pr(S_{\nu}<(1-\gamma)\expect{S_{\nu}\given N_{\nu}}\given N_{\nu}) & \le \Pr(\abs{S_{\nu}-\expect{S_{\nu}\given N_{\nu}}}<\gamma\expect{S_{\nu}\given N_{\nu}}\given N_{\nu})\\
    & \le 2\exp(-2\gamma^{2}\frac{\expect{S_{\nu}\given N_{\nu}}^{2}}{N_{\nu}})\\
    & = 2\exp(-2\gamma^{2}\frac{N_{\nu}}{\popthreshold^{2}}(\int_{0}^{\popthreshold}f(x)\intd x)^{2})\\
    & \le
    2\exp(-2\gamma^{2}\frac{N_{\nu}}{\popthreshold^{2}}(1-\gamma)^{2}\norm
    f_{1}^{2}),
  \end{align}
  for $\nu$ sufficiently large since $\popthreshold\to\infty\mbox{ as
  }\nu\to\infty$.  Whence,
  \begin{align}
    &\Pr(\frac{S_{\nu}}{\nu(1-\gamma)^{2}\norm f_{1}}<1-\gamma\given
    N_{\nu}\ge(1-\gamma)\nu\popthreshold) \\&\quad \le
    \Pr(\frac{S_{\nu}}{\expect{S_{\nu}\given N_{\nu}}}<1-\gamma\given
    N_{\nu}\ge(1-\gamma)\nu\popthreshold) \\&\quad \le
    2\exp(-2\gamma^{2}\frac{\nu}{\popthreshold}(1-\gamma)^{3}\norm
    f_{1}^{2}).
  \end{align}
  Using that $f(x)$ is monotonic and must be integrable we have that
  $f(x)=o(\frac{1}{x})$ so $\nu/\popthreshold\ge(\nu\log\nu)^{1/2}$
  and
  \[
  \Pr(\frac{S_{\nu}}{\nu(1-\gamma)^{2}\norm f_{1}}<1-\gamma\suchthat
  N_{\nu}\ge(1-\gamma)\nu\popthreshold)\le2\exp(-2\gamma^{2}(\nu\log\nu)^{1/2}(1-\gamma)^{3}\norm
  f_{1}^{2}).
  \]
  Finally, using
  $\lim_{\nu\to\infty}\frac{N_{\nu}}{\nu\popthreshold}\ge(1-\gamma)\mbox{
    a.s.}$ and the Borel--Cantelli lemma establishes
  \begin{align}
    \lim_{\nu\to\infty}\frac{S_{\floor{\nu}}}{\floor{\nu}+1} & \ge (1-\gamma)^{3}\norm f_{1}\mbox{ a.s.}\\
    \implies\lim_{\nu\to\infty}\frac{S_{\nu}}{\nu} & \ge
    (1-\gamma)^{3}\norm f_{1}\mbox{ a.s.}
  \end{align}
  and the result follows since $\gamma>0$ is arbitrary.
\end{proof}
We now have a promise that every point of the latent Poisson process
$\PPbelowth{\popthreshold}$ participates in the graph. We now
establish that, with high probability, as $\nu\to\infty$ an
arbitrarily large fraction of the points in $\PPbelowth{\threshold}$
connect to the popular connected core $\popgraph$. In particular, this
means an arbitrarily large fraction of the points of
$\PPbelowth{\threshold}$ participate in a single connected component
of $G_{\nu}$.
\begin{lem}
  \label{lem:all_below_threshold} Suppose $f$ does not have compact
  support. Let a point
  $(\theta_{i},\vartheta_{i})\in\PPbelowth{\threshold}$ be visible if
  $\theta_{i}\in\vertexset{G_{\nu}}$ and it participates in an edge
  connecting to $\popgraph$, and call a point invisible otherwise. Let
  $N_{\mbox{invis},\le\threshold}$ be the number of points in
  $\PPbelowth{T_{\nu}}$ that are invisible and let
  $N_{\mbox{vis},\le\threshold}$ be the number of points in
  $\PPbelowth{\threshold}$ that are visible, then for $\epsilon>0$
  \[
  \lim_{\nu\to\infty}P(N_{\mbox{invis},<T_{\nu}}>\epsilon
  N_{\mbox{vis},<T_{\nu}})=0.
  \]
\end{lem}
\begin{proof}
  By \cref{lem:Erdos_Renyi_near_0} it follows that as $\nu\to\infty$
  there are no invisible vertices below
  $\popthreshold=f^{-1}(\sqrt{\frac{\log\nu}{\nu}})$ so it suffices to
  bound the number of invisible vertices between $\popthreshold$ and
  $\threshold$. Conditional on $\popgraph$, each point
  $(\theta_{i},\vartheta_{i})\in\PPaboveth{\popthreshold}$ connects to
  $\popgraph$ independently with probability
  $1-\prod_{p\in\popgraph}(1-f(\vartheta_{i})f(p))\ge1-e^{-f(\vartheta_{i})S_{\nu}}$
  where $S_{\nu}=\sum_{p\in\popgraph}f(p)$. Since labeling each point
  of the Poisson process $\PPaboveth{\popthreshold}$ by whether or not
  it connects to $\popthreshold$ is, conditional on $\popgraph$, a
  marking of the Poisson process, we immediately have that the number
  of visible and invisible points in $\{
  (\theta_{i},\vartheta_{i})\in\PPnu\given\popthreshold<\threshold\} $
  are independent random variables and that there exists random
  variables $N_{\nu,\mbox{ub}}$ and $N_{\nu,\mbox{vis}}$ such that,
  \[
  N_{\nu,\mbox{ub}}\sim\poiDist(\nu\int_{\popthreshold}^{\threshold}e^{-f(x)S_{\nu}}\intd
  x)
  \]
  is a upper bound for $N_{\mbox{invis},<T_{\nu}}$ and
  \[
  N_{\nu,\mbox{vis}}\sim\poiDist(\nu\int_{\popthreshold}^{\threshold}1-e^{-f(x)S_{\nu}}\intd
  x)
  \]
  is an independent lower bound for $N_{\mbox{vis},<T_{\nu}})$.

  Thus a sufficient condition for the claim is
  $\Pr(\frac{N_{\nu,\mbox{ub}}}{N_{\nu,\mbox{vis}}}>\epsilon)\to0,\
  \nu\to\infty$.  Conditional on $S_{\nu}$, this is a ratio of
  independent Poisson random variables and this condition will hold if
  the ratio of their means goes to $0$:
  \begin{align}
    \frac{\nu\int_{\popthreshold}^{\threshold}e^{-f(x)S_{\nu}}\intd x}{\nu\int_{\popthreshold}^{\threshold}1-e^{-f(x)S_{\nu}}\intd x} & \le \frac{(\threshold-\popthreshold)e^{-f(T_{\nu})S_{\nu}}}{(\threshold-\popthreshold)}\\
    & \le e^{-(\frac{1}{\nu}+t_{\nu})S_{\nu}}.
  \end{align}
  Invoking
  $\lim_{\nu\to\infty}S_{\nu}/\nu\ge\frac{1}{2}\|f\|_{1}=1\as$ from
  \cref{lem:Erdos_Renyi_near_0} completes the result since this means
  $\lim_{\nu\to\infty}t_{\nu}S_{\nu}=\infty\as$
\end{proof}
The next step is to determine the total number of vertices above
$\threshold$ that connect to the popular connected core:
\begin{lem}
  \label{lem:visible_above}Suppose $f$ does not have compact
  support. Let
  \[
  N_{\nu,>T_{\nu}}=\abs{\{(\theta_{i},\vartheta_{i}) \in
    \PPaboveth{T_{\nu}} \given \exists p \in \vertexset{P_{\nu}}
    \mbox{ such that } \{ \theta_{i},p\} \in\edgeset{G_{\nu}}\} },
  \]
  be the number of points above $\threshold$ that connect to
  $\popgraph$.  Then there exists a random variable $N_{\nu,+}$ such
  that $N_{\nu,+}\le N_{\nu,>T_{\nu}}$ and
  \[
  N_{\nu,+}\given S_\nu
  \sim\poiDist(\nu\int_{\threshold}^{\infty}1-e^{-f(x)S_{\nu}}\intd x)
  \]
\end{lem}
\begin{proof}
  Conditional on $\popgraph$, each point
  $(\theta_{i},\vartheta_{i})\in\PPaboveth{\threshold}$ connects to
  $P_{\nu}$ independently with probability
  $1-\prod_{p\in\popgraph}(1-f(\vartheta_{i})f(p))\ge1-e^{-f(\vartheta_{i})S_{\nu}}$.
  This is a marking of the Poisson process so the random subset of
  $\PPaboveth{\threshold}$ that connects to $P_{\nu}$ is itself a
  Poisson process with rate
  $\nu(1-\prod_{p\in\popgraph}(1-f(\vartheta_{i})f(p)))$.  We may then
  further independently mark the points of this process such that the
  new random subset will be, conditional on $S_{\nu}$, a Poisson
  process with rate
  $\nu\int_{\threshold}^{\infty}1-e^{-f(x)S_{\nu}}\intd x$.  Let the
  number of points in this process be $N_{\nu,+}$ then it follows
  immediately that $N_{\nu,+}$ is a lower bound $N_{\nu,>T_{\nu}}$ and
  that $N_{\nu,+}\given S_\nu
  \sim\poiDist(\nu\int_{T_{\nu}}^{\infty}1-e^{-f(x)S_{\nu}}\intd x)$.
\end{proof}
The final step is to bound the number of vertices above $T_{\nu}$ that
will be neglected. These are the vertices that participate in edges lying
entirely above $T_{\nu}$ and have a minimum distance greater than $2$
to the popular subgraph $\popgraph$. Note that they may be part of the
giant component, but their contribution is negligible. We begin with a
small technical lemma:
\begin{lem}
  \label{lem:threshold_grows_slowly}Let $f:\NNReals\to[0,1]$ be
  monotonically decreasing and integrable, then
  $f\linverse(\frac{1}{t})=o(t)$.
\end{lem}
\begin{proof}
  Suppose otherwise so that $\exists c>0$ such that
  $f\linverse(\frac{1}{t})\ge ct$ infinitely often. Let $\{ t_{i}\}
  _{i=1}^{\infty}$ be a strictly increasing sequence of such $t$s,
  then for each $t_{i}$ there exists a box $B_{t_{i}}$ of area at
  least $c$ that lies under the graph: namely the box
  $\left[0,ct\right]\times\left[0,f(ct)\right]$.  For $\epsilon>0$ we
  may choose a subsequence $\{ \tilde{t}_{j}\}
  _{j=1}^{\infty}\subset\{ t_{i}\} _{i=1}^{\infty}$ such that
  $\abs{B_{t_{i}}\cap B_{t_{i+1}}}\le\epsilon$, so that the area below
  $f$ is bounded below by an infinite sum where each term has value at
  least $c-\epsilon>0$ thereby arriving at a contradiction.
\end{proof}
Following our interpretation of $T_{\nu}$ as a cutoff below which
every candidate vertex participates in the graph, the requirement
$\threshold=o(\nu)$ is obvious. Suppose otherwise, then there would be
$\Omega(\nu^{2})$ visible vertices in the graph and $\Theta(\nu^{2})$
expected edges, pushing the graph into the ultra-sparse regime where
$\abs{e(G_{\nu})}=O(\abs{\vertexset{G_{\nu}}})$.  The above lemma
shows that $\threshold=o(\nu)$ does indeed hold, since
$T_{\nu}=o(f\linverse(1/\nu))$ and $f\linverse(1/\nu)=o(\nu)$.  With
this result in hand,
\begin{lem}
  \label{lem:invisible_above}Suppose $f$ does not have compact
  support. Call a vertex $\theta_{i}\in\vertexset{G_{\nu}}$ ignored if
  $(\theta_{i},\vartheta_{i})\in\PPaboveth{\threshold}$ and its
  distance to $P_{\nu}$ is greater than $2$. Let $N_{\mbox{ignore}}$
  be the number of ignored vertices; then fixing $\epsilon>0$,
  \[
  \lim_{\nu\to\infty}\Pr(\frac{N_{\mbox{ignore}}}{\abs{\vertexset{G_{\nu}}}}>\epsilon)\to0.
  \]
\end{lem}
\begin{proof}
  \begin{figure}
    \begin{centering}
      \includegraphics[width=.25\linewidth]{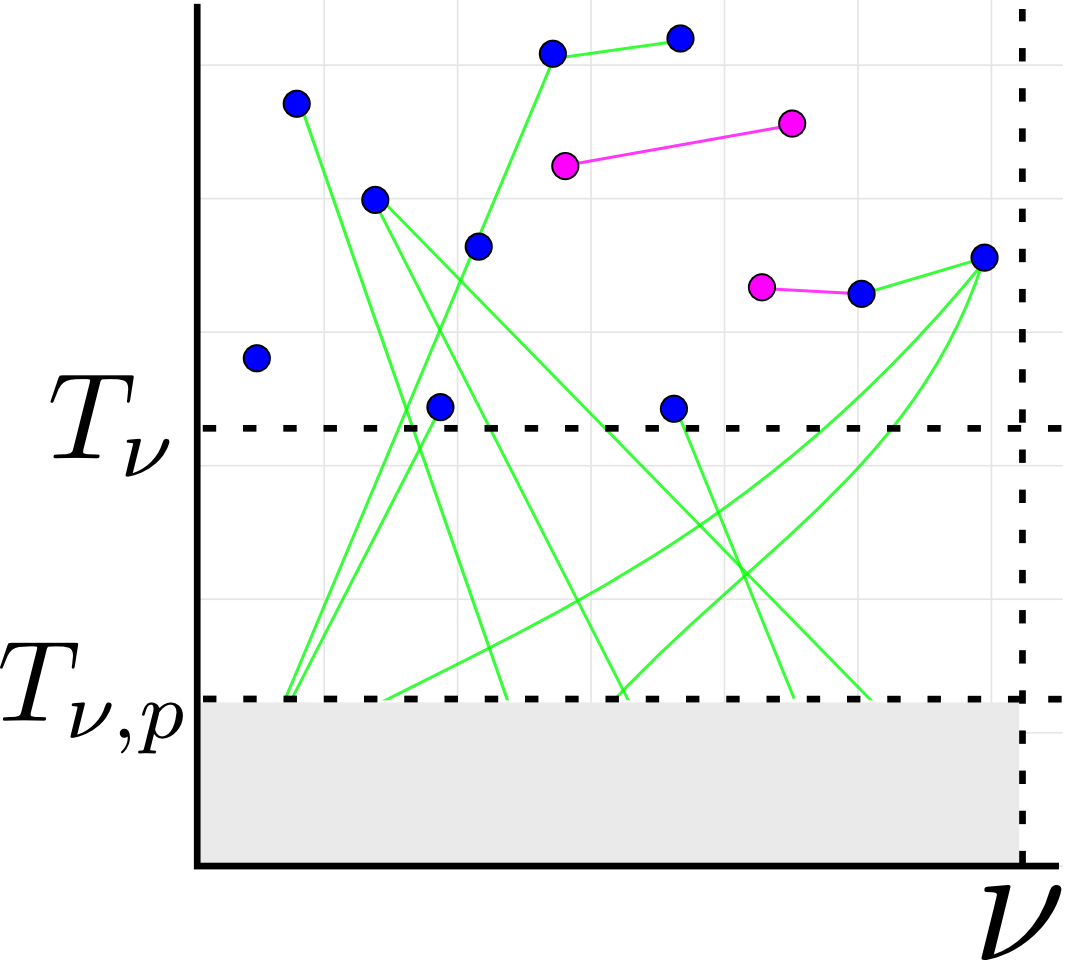}
      \par\end{centering}

    \caption{\label{fig:invis_above_figure}The structure of negligible
      vertices above $T_{\nu}$. Vertices with distance $>2$ to $P_{\nu}$
      (below $T_{\nu,p}$) are ignored, these are marked in magenta.}
  \end{figure}
  We mark each point in the Poisson process $\Pi_{\nu,>T_{\nu}}$ above
  $T_{\nu}$ by whether it participates in an edge with a terminus in
  $P_{\nu}$. As in \cref{lem:visible_above}, this forms a marking of
  the Poisson process conditional on $\popgraph$ so that the random
  subset of $\PPaboveth{\threshold}$ that is at distance one
  (\textbf{c}lose) to $\popgraph,$
  \[
  C_\nu=\{
  (\theta_{i},\vartheta_{i})\in\PPaboveth{T_{\nu}}\given\exists
  p\in\vertexset{P_{\nu}}\mbox{ such that }\{ \theta_{i},p\}
  \in\edgeset{G_{\nu}}\} ,
  \]
  and the remaining subset $\PPaboveth{T_{\nu}}\backslash C_\nu$ are
  independent Poisson processes conditional on $P_\nu$.

  Let $e_{\mbox{\ensuremath{\nu},ignore}}$ be the the number of edges
  in the induced subgraph of $G_{\nu}$ given by restricting the vertex
  set to $\PPaboveth{T_{\nu}}\backslash C_\nu$. It is immediate that
  $N_{\mbox{ignore}}\le 2e_{\mbox{\ensuremath{\nu},ignore}}$ (see
  \cref{fig:invis_above_figure}). Obviously
  $\abs{\vertexset{G_{\nu}}}>\abs{C_\nu}$ and by
  \cref{lem:visible_above} $\abs{C_\nu}>N_{\nu,+}$ so
  \[
  \Pr(\frac{N_{\mbox{ignore}}}{\abs{\vertexset{G_{\nu}}}}>\epsilon)\le
  \Pr(\frac{2e_{\mbox{\ensuremath{\nu},ignore}}}{N_{\nu,+}}>\epsilon),
  \]
  where in particular $e_{\mbox{\ensuremath{\nu},ignore}}$ and
  $N_{\nu,+}$ are independent conditional on $\PPbelowth{\threshold}$.

  We have very little distributional information about
  $e_{\mbox{\ensuremath{\nu},ignore}}$ so we use Markov's
  inequality. Since $\PPaboveth{T_{\nu}}\backslash C_\nu$ is a Poisson
  process with rate at most $\nu e^{-f(x)S_{\nu}}$ we may repeat the
  argument of \cref{thm:expected_edges} to bound
  $\expect{e_{\mbox{\ensuremath{\nu},ignore}}\given\popgraph}$ so that
  \begin{align}
    \expect{\frac{e_{\mbox{\ensuremath{\nu},ignore}}}{N_{\nu,+}}\given\PPbelowth{\threshold}} & = \frac{\expect{e_{\mbox{\ensuremath{\nu},ignore}}\given\PPbelowth{\threshold}}}{\expect{N_{\nu,+}\given\PPbelowth{\threshold}}}\\
    & \le
    2\frac{\nu^{2}(\int_{T_{\nu}}^{\infty}e^{-2S_{\nu}f(x)}f(x)\intd
      x)^{2}}{\nu\int_{T_{\nu}}^{\infty}1-e^{-f(x)S_{\nu}}\intd x}.
  \end{align}
  From this we see that the bound is $S_{\nu}$ measurable.  Taking
  $\gamma>0$ and working in the regime where
  $(1-\gamma)\le\frac{S_{\nu}}{\nu\|f\|_{1}}\le(1+\gamma)$ we have:
  \begin{align}
    \expect{\frac{e_{\mbox{\ensuremath{\nu},ignore}}}{N_{\nu,+}}\suchthat(1-\gamma)\le\frac{S_{\nu}}{\nu\|f\|_{1}}\le(1+\gamma)}
    & \le
    2\nu\frac{(\int_{T_{\nu}}^{\infty}e^{-2(1+\gamma)\|f\|_{1}\nu
        f(x)}f(x)\intd
      x)^{2}}{\int_{T_{\nu}}^{\infty}1-e^{-(1-\gamma)\|f\|_{1}\nu
        f(x)}\intd x}.
  \end{align}
  This can be treated by breaking up the integrals into the
  contributions above and below and upper threshold
  $\upperthreshold=f\linverse(\frac{1}{\nu})$. The numerator breaks up
  as,

\begin{align}
  \int_{\threshold}^{\upperthreshold}e^{-2(1+\gamma)\norm f_{1}\nu f(x)}f(x)\intd x+\int_{\upperthreshold}^{\infty}e^{-2(1+\gamma)\norm f_{1}\nu f(x)}f(x)\intd x\\
  \le O(\frac{\upperthreshold-\threshold}{\nu}) +
  O(\int_{\upperthreshold}^{\infty}f(x)\intd x),
\end{align}
where we have bounded the left term by the maximum of its
integrand. The denominator breaks up as,
\begin{align}
  \int_{\threshold}^{\upperthreshold}1-e^{-(1-\gamma)\norm f_{1}\nu f(x)}\intd x+\int_{\upperthreshold}^{\infty}1-e^{-(1-\gamma)\norm f_{1}\nu f(x)}\intd x\\
  \ge
  \Omega(\upperthreshold-\threshold)+\Omega(\nu\int_{\upperthreshold}^{\infty}f(x)\intd
  x),
\end{align}
where the bound on the right term follows from the fact that for
constant $c>0$ there exists $L_C$ depending only on $c$ such that
$1-e^{-cx}\ge L_c x$ for $x<1$.  Thus, in particular,
\begin{align}
  \expect{\frac{e_{\mbox{\ensuremath{\nu},ignore}}}{N_{\nu,+}}\suchthat(1-\gamma)\le\frac{S_{\nu}}{\nu\|f\|_{1}}\le(1+\gamma)}\\
  =O\left(\nu(\frac{\upperthreshold-\threshold}{\nu})^2\frac{1}{\upperthreshold-\threshold},
    \nu(\int_{\upperthreshold}^{\infty}f(x)\intd x)^2\frac{1}{\nu\int_{\upperthreshold}^{\infty}f(x)\intd x}\right),\\
\end{align}
and this goes to $0$ as $\nu\to\infty$; the left term because
$\upperthreshold=o(\nu)$ by \cref{lem:threshold_grows_slowly} and the
right term because $f$ is integrable and $\upperthreshold\to\infty$.

Putting all of this together and using that
$(1-\gamma)\le\lim_{\nu\to\infty}\frac{S_{\nu}}{\nu\|f\|_{1}}\le(1+\gamma)\as$
by \cref{lem:Erdos_Renyi_near_0} we have that:
\begin{align}
  \lim_{\nu\to\infty}\Pr(\frac{2e_{\mbox{\ensuremath{\nu},ignore}}}{N_{\nu,+}}>\epsilon) & = \lim_{\nu\to\infty}\Pr(\frac{2e_{\mbox{\ensuremath{\nu},ignore}}}{N_{\nu,+}}>\epsilon\suchthat(1-\gamma)\le\frac{S_{\nu}}{\nu\norm f_{1}}\le(1+\gamma))\\
  & \le \lim_{\nu\to\infty}2\expect{\frac{e_{\mbox{\ensuremath{\nu},ignore}}}{N_{\nu,+}}>\epsilon\suchthat(1-\gamma)\le\frac{S_{\nu}}{\nu\|f\|_{1}}\le(1+\gamma)}\\
  & = 0,
\end{align}
where the second line follows by Markov's inequality. This establishes
our claim.

\end{proof}
We can now put all of this together:
\begin{thm}
  Let $G$ be the KEG generated by $W=f(x)f(y)1[x\neq y]$, let
 $C_{1}(G_{\nu})$ be the largest connected component of $G_{\nu}$,
  and let $\epsilon>0$, then
  \[
  \lim_{\nu\to\infty}\Pr(\abs{C_{1}(G_{\nu})}>(1-\epsilon)\abs{\vertexset{G_{\nu}}})=1.
  \]
\end{thm}
\begin{proof}
  For $f$ with compact support this is a trivial consequence of
  \cref{thm:dense_iff_compact}, which shows that the graph is dense.
  For $f$ without compact support this is an immediate consequence of
  the lemmas of this section.
\end{proof}
A couple of concluding remarks are in order. Notice that the result
extends trivially to allow separable graphs that include self edges
because only a vanishing fraction of the vertices have a self edge. The
proofs in this section reveal some further interesting structure of
separable KEGs beyond connectivity, in particular:
\begin{enumerate}
\item If two points of a separable KEG are chosen at random there will
  be a very short path between them with high probability, even for
  very sparse random graphs. This is because both vertices very likely
  connect to the very dense subgraph $P_{\nu}$ by paths of length at
  most 2.
\item Although vertices of $G_{\nu}$ chosen uniformly at random are
  overwhelmingly likely to follow a degree distribution of the type
  given in \cref{thm:deg_dist_theorem} there are a vanishingly small
  fraction of the vertices (those in $P_{\nu}$) with much higher degree.
\end{enumerate}
Applied networks folk wisdom \cite{Newman:2009,Durrett:2006} holds
that real-world graphs often exhibit ``small world'' behaviour, with
very short paths between random vertices even for sparse graphs.
Similarly, it's common to observe that real-world graphs tend to
follow power law degree distribution except for the highest degree
vertices, which have much higher degree than would be expected from such
a law.  It's interesting that both of these features arise as emergent
behaviour of the simple random graph model considered in this section.

\section{Discussion}

This work was motivated by the need for a statistical framework 
for the analysis of the sparse graph structure of real-world networks.
The Kallenberg random graph model provides such a framework,
although the applicability and suitability of this framework---from either empirical or 
theoretical perspectives---is still to be determined.
Our work characterizing the limiting degree distribution and connectivity establish
that these models possess at least some of the properties of real-world networks 
we might hope to model.  The pioneering work of Caron and Fox yields further evidence.

The Kallenberg exchangeable graph model
is a natural generalization of the (dense) exchangeable graph model:
not only does the defining probabilistic symmetry still retain the
interpretation that the vertex labels do not carry any information
about the structure of the random graph, but graphons, which parametrize the exchangeable graphs, 
correspond with compactly-supported graphexes.
There are many deep results in the graphon theory for
which it is desirable to find sparse graph analogues. 
Several immediate goals worth pursuing are: identifying the sampling scheme that gives rise to KEGs; 
finding consistent estimators for a graphex, and identifying their properties;
and determining the graph limit theory corresponding to graphexes and its connection 
with existing graph limit theories for sparse graph sequences.
We now discuss these three directions in more detail.

A basic missing piece preventing us from confidently applying KEGs
to real-world network data
is a characterization of the processes that they model. In
particular, consider the problem of studying the properties of a very
large graph by sampling a small subgraph according to some random sampling design. 
Clearly any particular design licenses certain inferences and may even prevent others.
In this case the natural question is: what sampling schemes for 
subgraphs give rise to KEGs? 
It is well understood that a size-$n$ (dense) exchangeable graph model
corresponds to the process of observing the subgraph induced on $n$ vertices sampled 
uniformly at random from a large (even continuum-sized) graph.  One can see this interpretation in the work of Kallenberg~\cite{Kallenberg:1999} and the later independent work within graph theory, beginning with \citet{Lovasz:Szegedy:2006}.
The generative process for a KEG suggests the following 
sampling scheme for a finite graph $H$ corresponding to a KEG restricted to $[0,\nu]$: 
\begin{enumerate}
\item Sample a Poisson number $N$ of vertices uniformly at random with replacement 
         from $H$, where the mean of $N$ is $c\,\nu$.
\item Return the induced edge set, implicitly dropping isolated vertices.
\end{enumerate}
The corresponding graphex is $(0,0,c\cdot H)$ where $c \cdot H$ denotes the $c$-dilation of the empirical graphon associated with the finite graph $H$. (See \cref{sec:graphon_models}.)
The norm of the dilation is $\| c \cdot H\|_1 = c^2 \| H \|_1$, which we expect to approach zero as the 
graph $H$ becomes increasingly sparse.  This suggests normalizing, by taking the dilation $c$  
to be proportional to $\| H \|_1^{-\frac 1 2}$.  Such a renormalization bears some
resemblance to that of the $L^p$ theory discussed below, and is likely to feature in a graph limit theory.
This sampling scheme immediately suggests a notion of an empirical graphex, which one 
would expect to feature prominently in an estimation theory.
Identifying other sampling scheme(s) would provide both a sharp
understanding of the applicability of our models and substantive
guidance on how to subsample large networks.

In the absence of theoretical guidelines to the applicability of the KEG model, a pragmatic approach is to
simply fit KEG models to data and assess their appropriateness by
empirical evaluations, e.g., of their predictive performance. In practice, this entails
identifying classes of KEGs that both admit computationally
tractable inference procedures and are flexible enough to capture the
structure of real-world networks. The first step in this direction was
taken by Caron and Fox \cite{Caron_Fox_CRM_Graphs} with Bayesian
non-parametric models defined in terms of products of completely
random measures. The carefully crafted structure of their model 
allowed them to develop an efficient Markov Chain
Monte Carlo algorithm to fit their model to sparse graph data comprised of tens of thousands of vertices. 
More recently, 
\cite{Herlau:Schmidt:Morup:2015} have extended the work of Caron and
Fox to obtain an analogue of the well-known stochastic block model.  The analogue is easily seen to also be a KEG. 
Going forward, the close connection between
graphexes and graphons suggests that many of the existing models in the (dense) exchangeable graph framework 
will have natural analogues in KEG framework. This includes many popular
models in the literature, e.g.,\
\cite{Nowicki:Snijders:2001:A,Hoff2002,Airoldi2008,Miller2009,Lloyd:Orbanz:Ghahramani:Roy:2012};
see \cite{Orbanz_Roy_Exchangeable_Struct} for a review.

Finally, 
it is interesting to consider the connection with graph limit theory. 
There are at least two distinct contexts in which graphons arise: 
First, as we have already described in detail,
is as the structures characterizing the extreme elements among the exchangeable graphs. 
Second,  is as the limit objects for dense graph sequences~\cite{Lovasz:Szegedy:2006,Lovasz:Szegedy:2007,Lovasz:2013:A}.
The connection between the two perspectives is explained by
\cite{Diaconis:Janson:2007}.  The focus of the present paper is the
generalization of the first perspective to the sparse regime. 
Recent work~\cite{Borgs:Chayes:Cohn:Zhao:2014:sgc1,Borgs:Chayes:Cohn:Zhao:2014:sgc2}
has generalized the limit theory to the sparse regime by introducing a new notion of convergence and class
of limit objects called
$L^p$ graphons, which are symmetric integrable functions $W:[0,1]^2\to\NNReals$. 
The corresponding $W$-sparse random graph model is not projective, in contrast to the Kallenberg exchangeable graph model.  Understanding the link between the graphex theory and the $L^p$ graphon theory 
could provide new insights in both graph theory and the statistical analysis of networks.

\section*{Acknowledgements}
The authors would like to thank Nate Ackerman, Cameron Freer, Benson
Joeris, and Peter Orbanz for helpful discussions.  The authors would
also like to thank Mihai Nica for suggesting the proof of
\cref{lem:threshold_grows_slowly}.  This work was supported by
U.S. Air Force Office of Scientific Research grant \#FA9550-15-1-0074.

\printbibliography

\vfill 

\end{document}